\documentclass[11pt]{article}
\usepackage{amsmath,amssymb,amsthm,comment}
\usepackage{mathtools}
\usepackage{color}
\usepackage[nohead,margin=1.1in]{geometry}

\usepackage{booktabs}
\usepackage{threeparttable}
\usepackage{graphicx,psfrag,epstopdf}
\usepackage{algorithmic,algorithm}
\usepackage{url}
\graphicspath{{./pics/}}
\usepackage{tcolorbox}
\allowdisplaybreaks
\theoremstyle{plain}% default
\newtheorem{theorem}{Theorem}[section]
\newtheorem{remark}{Remark}[section]
\newtheorem{definition}{Definition}[section]
\newtheorem{lemma}{Lemma}[section]

\newtheorem{prop}{Proposition}[section]
\newtheorem{cor}{Corollary}[section]

\numberwithin{equation}{section}
\def\E{\mathbb{E}}

\usepackage{tikz}
\usetikzlibrary{positioning}
\begin{document}
\title{On the recovery of two function-valued coefficients in the Helmholtz equation for inverse scattering problems via neural networks}
\author{Zehui Zhou\thanks{Department of Mathematics,
Rutgers University, Piscataway, New Jersey, USA. (\texttt{zz569@math.rutgers.edu})}}
\date{}
\maketitle

\begin{abstract}
Recently, deep neural networks (DNNs) have become powerful tools for solving inverse scattering problems. 
However, the approximation and generalization rates of DNNs for solving these problems remain largely under-explored. 
In this work, we introduce two types of combined DNNs (uncompressed and compressed) to reconstruct {two function-valued coefficients} in the Helmholtz equation for inverse scattering problems from the scattering data at two different frequencies. 
An analysis of the approximation and generalization capabilities of the proposed neural networks for simulating the regularized pseudo-inverses of the linearized forward operators in direct scattering problems is provided. 
The results show that, with sufficient training data and parameters, the proposed neural networks can effectively approximate the inverse process with desirable generalization.
Preliminary numerical results show the feasibility of the proposed neural networks for recovering two types of isotropic inhomogeneous media. 
Furthermore, the trained neural network is capable of reconstructing the isotropic representation of certain types of anisotropic media.

\end{abstract}

\section{Introduction}\label{sec:intro}

%\subsection{Problem setup and main results}

The scattering of the time-harmonic incident plane wave $u^i(x):=e^{i\omega x\cdot \zeta}$ with the probing wave frequency $\omega$ and incident direction $\zeta\in\mathbb{S}$ (with $\mathbb{S}$ denoting the unit sphere) by an inhomogeneity occupying a bounded Lipschitz domain $\Omega \subset \mathbb{R}^2$ is mathematically formulated {as follows}: find the total wave field $u\in H_{loc}^1(\mathbb{R}^2)$ with $u=u^s+u^i$ such that
\begin{align}
&\nabla\cdot a(x) \nabla u+\omega^2 n(x) u=0  \;\;\mbox{in}\;\; \mathbb{R}^2,\label{eqn:prob}\\
&\lim_{r\to \infty} r^\frac12 \left(\dfrac{\partial u^s}{\partial r}-i\omega u^s\right)=0 \;\;\mbox{with}\;\;r=|x|  \label{eqn:src}
\end{align}
where 
\begin{align*}
a(x)=&\left\{
\begin{array}{cl}
    I & \mbox{in }\; \mathbb{R}^2\setminus \bar{\Omega}\\
    I+\gamma(x) & \mbox{in }\; \Omega
\end{array}
\right. \quad \mbox{ and }\quad 
n(x)=\left\{
\begin{array}{cl}
    1 & \mbox{in }\; \mathbb{R}^2\setminus \bar{\Omega}\\
    1+\eta(x) & \mbox{in }\; \Omega
\end{array}
\right.
\end{align*}
with $I$ being the identity function, $\gamma \in C^1(\bar{\Omega})$ (entry-wise in the case when $\gamma$ is a $2\times 2$ positive definite matrix-valued function) and $\eta \in L^\infty(\bar{\Omega})$, and the Sommerfeld radiation condition \eqref{eqn:src} on the scattered wave field $u^s$ is satisfied uniformly with respect to $\hat{x}=x/|x|$. 
{Problem \eqref{eqn:prob}--\eqref{eqn:src} models} the scattering of acoustic waves by a (possibly anisotropic) inhomogeneous medium with contrasts in both the sound speed and density \cite{colton-kress} in two dimensions, or the specially polarized electromagnetic waves by a (possibly orthotropic) inhomogeneity with contrasts in both electric permittivity and magnetic permeability \cite[Section 5.1]{CakoniColton:2005}. Thanks to Sommerfeld radiation condition \eqref{eqn:src}, the scattered filed $u^s$ assumes the following asymptotic behavior
$$u^s(x,\zeta)=\frac{e^{i\omega |x|}}{\sqrt{|x|}}\Lambda^\omega(\hat{x},\zeta)+\mathcal{O}\left(\frac{1}{|x|^{3/2}}\right),$$
where $\Lambda^\omega(\hat{x},\zeta)$, known as the far-field pattern of the scattered field, is a function of $\hat x \in {\mathbb S}$ for a fixed $\zeta$.
\begin{definition}[Scattering Data and the Inverse Problem] {\em The set of measured scattering data
$$\left \{\Lambda^\omega(\hat{x},\zeta): \, \hat x\in {\mathbb S}, \zeta\in{\mathbb S}\right\}$$
is called the {\it scattering data} at frequency $\omega$. The {\it inverse scattering problem} we consider here is to determine the contrasts $\gamma$ and $\eta$ from scattering data at two different frequencies $\omega_1\neq \omega_2$.}
\end{definition}
 
It is well-known that the uniqueness of the inverse problem holds in the case when $\gamma$ is a scalar function (i.e. isotropic media) \cite{uniq2} (see \cite{uni2, uni1} for coefficients less regular than assumed here). 
In the case of anisotropic media, i.e. when $\gamma$ is a $2\times 2$ symmetric matrix-valued function, there is a natural obstruction to the unique determination of $\gamma$ and $\eta$ (even if the scattering data is known for all frequencies). 
Indeed, let ${\mathcal O}$ be a simply connected region containing the support of the perturbations, i.e. $\overline\Omega \subset {\mathcal O}$, and let $\Phi$ be a sufficiently smooth diffeomorphism from ${\mathcal O}$ onto ${\mathcal O}$, with $\Phi(x)=x$ on $\partial {\mathcal O}$, and then define
$$
\tilde a =  \frac{D\Phi\,  a\,  D\Phi^{\top}}{|\det D \Phi |} \circ \Phi^{-1} ~~\hbox{ and }~~ \tilde n = \frac{n}{|\det D\Phi|}\circ \Phi^{-1} \qquad \mbox{in}\; {\mathcal O},
$$
where the $2\times2$ matrix-valued function $D\Phi$ is the Jacobian of $\Phi$ and $(\tilde a,\tilde n, {\mathcal O})$ is referred to as the pushforward of $(a, n, {\mathcal O})$ under $\Phi$. 
It is known that the operator $\nabla \cdot \tilde a\nabla + \omega^2 \tilde n$ has the same Dirichlet-to-Neumann data map (or Cauchy data) as $\nabla \cdot a\nabla  + \omega^2  n $ on $\partial {\mathcal O}$, hence the same scattering data, for any  frequency $\omega>0$ (see e.g. \cite{uniq2, KV}). The perturbations $\gamma, \eta$ and $(\tilde a-I), (\tilde n-1)$ yield the same scattering data at any frequency $\omega$, and thus cannot be identified from the scattering  data. 
In ${\mathbb R}^2$, owing to the isothermal coordinates, it is established that the set of all possible pushforwards $(a, n, {\mathcal O})$ contains a unique isotropic coefficient $(\tilde a,\tilde n, {\mathcal O})$, i.e., where $\tilde a$ is a scalar function \cite{john2, john1}.  
In this paper, although the neural network is built and trained to determine scalar functions $\gamma$ and $\eta$, we also present the results of numerical examples where the scattering data corresponds to anisotropic media.

\begin{remark}
{\em In practice, we have access to the scattering data for only a limited number of incident and observation directions.
Thanks to the analyticity of the far-field pattern, the above uniqueness results for the inverse problem also hold when the scattering data is considered for $\hat x \in {\mathbb S}_r\subset {\mathbb S}$ and $\zeta \in {\mathbb S}_t\subset {\mathbb S}$, where ${\mathbb S}_r$ and ${\mathbb S}_t$ are open subsets (possibly the same) of the unit sphere. 
However, the resolution of the reconstructions decreases for small spatial apertures ${\mathbb S}_r$ and ${\mathbb S}_t$.}
\end{remark}

By noting that the incident wave field $u^i(x)=e^{i\omega x\cdot \zeta}$ solves the Helmholtz equation with the homogeneous media (i.e. $\gamma\equiv0$ and $\eta\equiv0$) in ${\mathbb R}^2$,  the perturbed Helmholtz equation \eqref{eqn:prob} for the scattered field takes the form 
\begin{align}
\nabla\cdot a(x) \nabla u^s+\omega^2 n(x)u^s=- \nabla\cdot \gamma(x)\nabla u^i-\omega^2 \eta(x) u^i \;\;\;\mbox{for}\;\; \; x\in \mathbb{R}^2 \label{eqn:prob1}
\end{align}
or equivalenty
\begin{equation}
\Delta u^s+\omega^2 u^s=- \nabla\cdot \gamma(x)\nabla u-\omega^2 \eta(x) u \;\;\;\mbox{for}\;\;\; x\in  \mathbb{R}^2. \label{eqn:prob2}
\end{equation}
Together with the Sommerfeld radiation condition \eqref{eqn:src}, it is easy to verify that $u^s$ satisfies the following Lippmann-Schwinger volume integral equation,
\begin{equation}\label{lsvi}
u^s(x,\zeta)= \nabla_x\cdot \int_{\Omega} G(x-y)\gamma(y)\nabla_y (u^s+u^i)(y) {\rm d}y+\int_{\Omega} G(x-y)\omega^2 \eta(y) (u^s+u^i)(y){\rm d}y,
\end{equation}
with $G(x-y)=i/4 H_0^{(1)}(\omega |x-y|)$ being the full-space fundamental solution of the Helmholtz equation  in ${\mathbb R}^2$, where $H_0^{(1)}$ is the Hankel function of the first kind of order zero. 
The volume integral equation \eqref{lsvi} allows $\gamma$ to have jumps across $\partial \Omega$, and we refer the reader to  \cite{Costabel:2015} for the sufficient conditions on its solvability. 
Due to the dependency of $u^s$ on $\gamma$ and $\eta$, equation \eqref{lsvi} is nonlinear with respect to $\gamma$ and $\eta$. 
Limited by the structure of neural networks, we consider the linearized version of \eqref{lsvi} to reconstruct the perturbations $\gamma$ and $\eta$:
 \begin{eqnarray*}
&&u^s(x,\zeta)\approx \nabla_x\cdot \int_{\Omega} G(x-y)\gamma(y)\nabla_y u^i(y) {\rm d}y+\int_{\Omega} G(x-y)\omega^2 \eta(y)u^i(y){\rm d}y \\
&& \qquad \qquad  =  \int_{\Omega} \nabla_x\cdot  G(x-y) \gamma(y) i\omega \zeta e^{i\omega y\cdot \zeta} {\rm d}y+\int_{\Omega} G(x-y)\omega^2 \eta(y)e^{i\omega x\cdot \zeta}{\rm d}y\\
&&  \qquad \qquad  =  \int_{\Omega}i \omega \nabla_x G(x-y)\cdot \zeta\,  \gamma(y)e^{i\omega y\cdot \zeta} {\rm d}y+\int_{\Omega} G(x-y)\omega^2 \eta(y)e^{i\omega x\cdot \zeta}{\rm d}y.
\end{eqnarray*}
This approximation of the scattered field can be justified for small perturbation and can be viewed as resulting from differentiating the scattered field with respect to a small parameter $\epsilon>0$ for given perturbations $\epsilon \gamma$ and $\epsilon \eta$ of the homogeneous media. 
The asymptotic expressions 
\begin{align*}
G(x-y)=&\frac{e^{i\frac{\pi}{4}}}{\sqrt{8\pi \omega}}\frac{e^{i\omega |x|}}{\sqrt{|x|}}e^{-i \omega \hat{x}\cdot y}+{O}\left(\frac{1}{|x|^{3/2}}\right)\qquad \qquad \mbox{and}\\
\nabla_x G(x-y)=&-\nabla_y G(x-y)=\frac{e^{i\frac{\pi}{4}}}{\sqrt{8\pi \omega}}\frac{e^{i\omega |x|}}{\sqrt{|x|}} i\omega \hat x e^{-i \omega \hat{x}\cdot y}+{O}\left(\frac{1}{|x|^{3/2}}\right),
\end{align*}
together with the definition of the far-field pattern $\Lambda^\omega(\hat{x},\zeta)$ (i.e. the scattering data), imply that
\begin{align*}
\Lambda^\omega(\hat{x},\zeta)\approx \frac{e^{i\frac{\pi}{4}}}{\sqrt{8\pi \omega}}   \left(\int_{\Omega} -\omega^2(\hat x\cdot \zeta) e^{-i \omega \hat{x}\cdot y} \gamma(y)  e^{i\omega y\cdot \zeta} {\rm d}y+\int_{\Omega} e^{-i \omega \hat{x}\cdot y}\omega^2 \eta(y)e^{i\omega y\cdot \zeta}  {\rm d}y\right).
\end{align*}
The support $\Omega$ of the perturbations is part of the unknown. {By noting} that $\gamma$ and $\eta$ are both zero outside $\Omega$, we consider a sufficiently large ball $\mathcal{B}_R$ of radius $R$, centered at the origin, which we know a priori contains the perturbations, i.e., $\overline \Omega\subset \mathcal{B}_R$. Then {by extending} $\gamma$ and $\eta$ by zero in $\mathcal{B}_R\setminus \overline{\Omega}$, we arrive at the following (approximate) expression  for the scattering data
\begin{align}
\Lambda^\omega(\hat{x},\zeta)=
&\frac{e^{i\frac{\pi}{4}}\omega^2}{\sqrt{8\pi \omega}}  \left(\int_{\mathcal{B}_R} - (\hat{x}\cdot \zeta)   e^{i\omega (\zeta-\hat{x})\cdot y}\gamma(y){\rm d}y+\int_{\mathcal{B}_R}  e^{i\omega (\zeta-\hat{x})\cdot y}\eta(y){\rm d}y\right)\nonumber\\
:=&c_\omega \big([F_1^\omega\gamma](\hat{x},\zeta)+[F_2^\omega\eta](\hat{x},\zeta)\big),\label{eqn:F1F2}
\end{align}
where the constant $c_\omega=\frac{e^{i\frac{\pi}{4}}\omega^2}{\sqrt{8\pi \omega}}$ and $F_1^\omega$, $F_2^\omega$ are viewed as linear operators on $\gamma$ and $\eta$ {such that 
\begin{align*}
[F_1^\omega\gamma](\hat{x},\zeta)=\int_{\mathcal{B}_R} - (\hat{x}\cdot \zeta)   e^{i\omega (\zeta-\hat{x})\cdot y}\gamma(y){\rm d}y \quad \mbox{and}\quad [F_2^\omega\eta](\hat{x},\zeta)=\int_{\mathcal{B}_R}  e^{i\omega (\zeta-\hat{x})\cdot y}\eta(y){\rm d}y.
\end{align*}} 
Note that although we have assumed that $\gamma\in C^1(\overline{\Omega})$, in general,  $\gamma$ is not necessarily zero on $\partial \Omega$ hence its extension by zero is only bounded. 
In order to reconstruct the perturbations $\gamma$ and $\eta$, we need to solve the linearized ill-posed inverse problem with scattering data generated for two different interrogating frequencies $\omega_1$ and $\omega_2$, namely 
\begin{equation}\label{eqn:prob_linear_sys}
F\begin{bmatrix}
\gamma\\
\eta    
\end{bmatrix}=\Lambda   
, \quad\mbox{where} \quad \Lambda := \begin{bmatrix}
\Lambda^{\omega_1}\\
\Lambda^{\omega_2}    
\end{bmatrix} \quad \mbox{and}\quad 
F:=\begin{bmatrix}
F_1^{\omega_1} &F_2^{\omega_1}\\
F_1^{\omega_2} &F_2^{\omega_2}  
\end{bmatrix}.
\end{equation}
That is, we seek to invert (in practice a discretized version of) $F: (L^\infty(\mathcal{B}_R))^2\to (L^2({\mathbb S}\times {\mathbb S}))^2$. 
Throughout the paper, we assume that $\omega_1<\omega_2$. 
The benefit of the linearization \eqref{lsvi} is that the obtained operators  $F_1^\omega$ and $F_2^\omega$ become of Fourier convolution type, which is suitable for neural network inversion.  {Due to the compactness of $F$ \cite{colton-kress},} the problem \eqref{eqn:prob_linear_sys} is ill-posed, i.e., a solution may not exist and even if it does exist, the solution is highly unstable with respect to the noisy data $\Lambda=\left(\Lambda^{\omega_1}(\hat{x},\zeta), \Lambda^{\omega_2}(\hat{x},\zeta)\right)^T$, thus regularization is  required for the stable and accurate numerical solutions. 
To this end, it is convenient to work on  Hilbert spaces, hence we view $F$ as an operator $(L^2(\mathcal{B}_R))^2\to (L^2({\mathbb S}\times {\mathbb S}))^2$, and denote by $F^*: (L^2({\mathbb S}\times {\mathbb S}))^2\to (L^2(\mathcal{B}_R))^2$ its $L^2$ adjoint. 
In this work, we apply the Tikhonov regularization technique \cite{EnglHankeNeubauer:1996}, and thus considered the regularized minimization problem with regularization parameter $\alpha>0$:
$$\min_{\gamma,\eta}\bigg\|
F\begin{bmatrix}
\gamma\\
\eta    
\end{bmatrix}-\Lambda   
\bigg\|^2_{\mathbb{S}\times \mathbb{S}}+\alpha\bigg\|
\begin{bmatrix}
\gamma\\
\eta    
\end{bmatrix}\bigg\|^2_{\mathcal{B}_R},$$
where generically for a vector-valued function $v:=\left(v_1, v_2\right)\in L^2({\mathcal O})\times L^2({\mathcal O})$ we denote by $\|v\|_{{\mathcal O}}=(\int_{\mathcal O}\|v(y)\|_2^2{\rm d}y)^\frac12$ with $\|\cdot\|_2$ denoting the standard Euclidean norm of a vector (and also the spectral norm of a matrix later in the paper). 
The unique minimizer to this problem, i.e.,
\begin{equation}
(\gamma_\alpha,\eta_\alpha)^T = F_\alpha^\dag \Lambda \quad \mbox{with}\quad F_\alpha^\dag=(F^*F+\alpha I)^{-1}F^*,
\end{equation} 
defines the regularized pseudo-inverse $F_\alpha^\dag$ of the forward linear operator $F$. 
The main goal of this work is to construct an approximation of $F_\alpha^\dag$ using neural networks.

Inspired by the works \cite{LiDemanetZepeda:2022, ZhangZepeda-NunezLi:2024}, we construct two types of combined (fully connected followed by convolutional) neural networks to approximate the regularized pseudo-inverse $F_\alpha^\dag=(F^*F+\alpha I)^{-1}F^*$.
{The explicit construction of the proposed neural networks and their approximation and generalization analyses together represent the main technical contributions of this work.}  
In specific, the Fourier integral operator $F^*$ can be decomposed into {double integral} operators with related integral kernels, followed by suitable linear radial scaling operators (see Proposition \ref{prop:F*_formula}), where the integral kernels can be further approximated by low-rank operators due to its complementary low-rank property (see Theorem \ref{thm:low_rank}). 
Thus, after discretization, $F^*$ can be approximated by the fully connected neural networks \cite{DeVoreHaninPetrova:2021} of shared weights layers on the polar coordinates (see Theorem \ref{thm:BF_0_app}), where the weights can be further approximated by the butterfly factorization \cite{CandesDemanetYing:2009,LiYangMartin:2015} (see Theorem \ref{thm:BF_1_app}). % due to its complementary low-rank property (see Theorem \ref{thm:low_rank}). 
Furthermore, the convolution-type pseudo-differential operator $F^*F+\alpha I$ and its inverse can be simulated by the convolutional neural networks \cite{FanFeliuLin:2019} with scaled residual layers \cite{HeZhangRen:2016} on the Cartesian coordinates (see Theorem \ref{thm:CNN_app}) in view of their translation invariant property (see Proposition \ref{prop:F*F_formula}). 
In this work, we shall discuss the neural networks without and with the butterfly structure (middle level butterfly factorization), named uncompressed and compressed neural networks in \cite{ZhangZepeda-NunezLi:2024}, respectively.
{It is worth noting that the case we consider in this work, involving two unknown function-valued coefficients, is fundamentally distinct from the single-coefficient scenario \cite{LiDemanetZepeda:2022, ZhangZepeda-NunezLi:2024}. 
This distinction introduces challenges for both theoretical and numerical analyses.}

The rest of the paper is organized as follows. 
Properties of the Fourier integral operator $F^*$ and pseudo-differential operator $F^*F+\alpha I$ are studied in Section \ref{sec:prop}, and the main theoretical results %of this paper 
concerning the approximation analysis for the combined neural network {are presented} in Section \ref{sec:NN_approx} while the generalization analysis %for the proposed neural networks 
is provided in Section \ref{sec:gene}. 
For sake of the reader convenience, we recall that the approximation error is defined as the $L^\infty$ distance between the true coefficients $\gamma, \eta$ and the reconstruction produced by the neural network for a given input. 
On the other hand, the generalization error, loosely speaking, measures the performance of the neural network when {it is trained} on a specific set of data but applied to a larger one. 
More specifically, Theorem \ref{thm:BF_0_app} and Theorem \ref{thm:BF_1_app} prove the approximation error for the fully connected parts of uncompressed and compressed neural networks, respectively, while Theorem \ref{thm:CNN_app} establishes the approximation error for the convolutional part of both proposed neural networks. 
{Theorem \ref{thm:gen_0} provides} generalization results for uncompressed neural networks whereas Theorem \ref{thm:gen_1} proves generalization results for compressed neural networks. Preliminary {numerical examples} of the performance of the proposed neural networks are provided in Section \ref{num}. 
It is worth noting that although the neural networks proposed in this work are based on the regularized pseudo-inverse of the linearized problem, in numerical experiments, we can increase the nonlinearity of the structures to improve their expressiveness, thereby obtaining better reconstructions than simply solving the linearized problem using traditional methods.

%%%%%%%%%%%%%%%%%%%%%%%%%%%%%%%%%%%%%%%%%%%%%%%%%%%%%%%%%%%%%%%%%%

\section{Properties of the regularized pseudo-inverse $F_\alpha^\dag$}\label{sec:prop}
In this section, we investigate the properties of the two components of the regularized pseudo-inverse $F_\alpha^\dag=(F^*F+\alpha I)^{-1}F^*$, namely $F^*:\left(L^2({\mathbb S}\times {\mathbb S})\right)^2\to \left(L^2(\mathcal{B}_R)\right)^2$ and $(F^*F+\alpha I)^{-1}:\left(L^2(\mathcal{B}_R)\right)^2\to \left(L^2(\mathcal{B}_R)\right)^2$ , in order to effectively utilize these features for the network design in Section \ref{sec:NN_approx}. First, we derive {an explicit formula} for $F^*$ on the polar coordinates in Propositions \ref{prop:F*_formula}. 
\begin{prop}\label{prop:F*_formula}
Let $F_1^\omega$, $F_2^\omega$ and $F$ be the operators defined in \eqref{eqn:F1F2} and \eqref{eqn:prob_linear_sys} respectively. Then their adjoint operators over the scattering data on the polar coordinates $(\theta,\rho)\in [0,2\pi]\times [0,\infty)$ are given by
\begin{align*}
F^*\Lambda=\begin{bmatrix}
F_1^{\omega_1 *} &F_1^{\omega_2 *}\\
F_2^{\omega_1 *} &F_2^{\omega_2 *} 
\end{bmatrix}\begin{bmatrix}
\Lambda^{\omega_1}\\
\Lambda^{\omega_2}    
\end{bmatrix}=
\sum_{i=1}^2\begin{bmatrix}
F_1^{\omega_i *}\Lambda^{\omega_i}\\
F_2^{\omega_i *}\Lambda^{\omega_i}
\end{bmatrix}
\end{align*}
and
\begin{align*}
[F_1^{\omega_i *}\Lambda^{\omega_i}](\theta,\rho)
=&-c_{\omega_i}\mathcal{S}_{\omega_i}\int_0^{2\pi} \overline{k(\theta_{\hat{x}},\rho)}\int_0^{2\pi} k(\theta_{\zeta},\rho)\cos(\theta_{\hat{x}}-\theta_{\zeta})  \Lambda_\theta^{\omega_i}(\theta_{\hat{x}},\theta_{\zeta}){\rm d}\theta_{\zeta}{\rm d}\theta_{\hat{x}},\\
[F_2^{\omega_i *}\Lambda^{\omega_i}](\theta,\rho)=&c_{\omega_i}\mathcal{S}_{\omega_i}\int_0^{2\pi} \overline{k(\theta_{\hat{x}},\rho)} \int_0^{2\pi} k(\theta_{\zeta},\rho)\Lambda_\theta^{\omega_i}(\theta_{\hat{x}},\theta_{\zeta}){\rm d}\theta_{\zeta}{\rm d}\theta_{\hat{x}},
\end{align*}
for $i=1,2$, with the linear radial scaling operator $\mathcal{S}_\omega f(\theta,\rho)=f(\theta,\frac{\omega}{\omega_1}\rho)$, the integral kernel $k(t,\rho)=e^{-i\omega_1 \rho \cos t}$ and the shifted far-field pattern $\Lambda_\theta^{\omega}(\theta_{\hat{x}},\theta_{\zeta})=\Lambda^{\omega}(\theta_{\hat{x}}+\theta,\theta_{\zeta}+\theta)$.
\end{prop}
\begin{proof}
{
By the definition of the adjoint operator, there holds
\begin{align*}
F^*=\begin{bmatrix}
F_1^{\omega_1 *} &F_1^{\omega_2 *}\\
F_2^{\omega_1 *} &F_2^{\omega_2 *} 
\end{bmatrix}.
\end{align*}}
For any $\omega$, it {follows} directly from the definitions of $F_1^\omega$ and $F_2^\omega$ in \eqref{eqn:F1F2} that
\begin{align*}
\langle \gamma, F_1^{\omega*}\Lambda^\omega \rangle_{\mathcal{B}_R}
=\langle F_1^{\omega}\gamma,\Lambda^\omega \rangle_{\mathbb{S}\times\mathbb{S}}
=&\int_{\mathbb{S}\times\mathbb{S}} \overline{\Lambda^{\omega}(\hat {x} ,\zeta)}\bigg(c_{\omega} \int_{\mathcal{B}_R} -(\hat{x}\cdot \zeta)  e^{i\omega (\zeta-\hat{x})\cdot y}\gamma(y){\rm d}y\bigg){\rm d}\hat{x}\, {\rm d}\zeta\\
=&\int_{\mathcal{B}_R}\gamma(y)\bigg(\overline{\int_{\mathbb{S}\times\mathbb{S}} -c_{\omega}(\hat{x}\cdot \zeta)  e^{i\omega (\hat{x}-\zeta)\cdot y}\Lambda^{\omega}(\hat {x} ,\zeta){\rm d}\hat{x}\, {\rm d}\zeta}\bigg){\rm d}y,\\
\langle \eta, F_2^{\omega*}\Lambda^\omega \rangle_{\mathcal{B}_R}
=\langle F_2^{\omega}\eta,\Lambda^\omega \rangle_{\mathbb{S}\times\mathbb{S}}
=&\int_{\mathbb{S}\times\mathbb{S}} \overline{\Lambda^{\omega}(\hat {x} ,\zeta)}\bigg(c_{\omega} \int_{\mathcal{B}_R}  e^{i\omega (\zeta-\hat{x})\cdot y}\eta(y){\rm d}y\bigg){\rm d}\hat{x}\, {\rm d}\zeta\\
=&\int_{\mathcal{B}_R}\eta(y)\bigg(\overline{\int_{\mathbb{S}\times\mathbb{S}} c_{\omega} e^{i\omega (\hat{x}-\zeta)\cdot y}\Lambda^{\omega}(\hat {x} ,\zeta){\rm d}\hat{x}\, {\rm d}\zeta}\bigg){\rm d}y,
\end{align*}
{which imply}
\begin{align}
[F_1^{\omega *}\Lambda^{\omega}](y)=&-c_{\omega}\int_{\mathbb{S}\times\mathbb{S}} (\hat{x}\cdot \zeta)   e^{i\omega (\hat{x}-\zeta)\cdot y}\Lambda^{\omega}(\hat{x},\zeta){\rm d}\hat{x}\, {\rm d}\zeta,\label{eqn:F1*Cart}\\
[F_2^{\omega *}\Lambda^{\omega}](y)=&c_{\omega}\int_{\mathbb{S}\times\mathbb{S}} e^{i\omega (\hat{x}-\zeta)\cdot y}\Lambda^{\omega}(\hat{x},\zeta){\rm d}\hat{x}\, {\rm d}\zeta, \label{eqn:F2*Cart}
\end{align}
for any $y\in \mathcal{B}_R$. 
%Similarly, by the definition of $F$ in \eqref{eqn:prob_linear_sys}, there holds
%\begin{align*}
%&\bigg \langle \begin{bmatrix}
%\gamma\\
%\eta    
%\end{bmatrix}, F^*\begin{bmatrix}
%\Lambda^{\omega_1}\\
%\Lambda^{\omega_2}    
%\end{bmatrix} \bigg \rangle_{\mathcal{B}_R}
%=\bigg \langle F\begin{bmatrix}
%\gamma\\
%\eta    
%\end{bmatrix}, \begin{bmatrix}
%\Lambda^{\omega_1}\\
%\Lambda^{\omega_2}    
%\end{bmatrix} \bigg \rangle_{\mathbb{S}\times\mathbb{S}}\\
%=&\sum_{i=1}^2 \int_{\mathbb{S}\times\mathbb{S}} \overline{\Lambda^{\omega_i}(\hat {x} ,\zeta)}\bigg(c_{\omega_i} \int_{\mathcal{B}_R} -(\hat{x}\cdot \zeta)  e^{i\omega_i (\zeta-\hat{x})\cdot y}\gamma(y)+ e^{i\omega_i (\zeta-\hat{x})\cdot y}\eta(y){\rm d}y\bigg){\rm d}\hat{x}\, {\rm d}\zeta\\
%=&\sum_{i=1}^2 \int_{\mathcal{B}_R}\gamma(y)\overline{[F_1^{\omega_i *}\Lambda^{\omega_i}](y)}{\rm d}y
%+\sum_{i=1}^2 \int_{\mathcal{B}_R}\eta(y)\overline{[F_2^{\omega_i *}\Lambda^{\omega_i}](y)}{\rm d}y,
%\end{align*}
%which indicates
%\begin{align*}
%F^*=\begin{bmatrix}
%F_1^{\omega_1 *} &F_1^{\omega_2 *}\\
%F_2^{\omega_1 *} &F_2^{\omega_2 *} 
%\end{bmatrix}.
%\end{align*}
Now, we convert the Cartesian coordinates to the polar coordinates with $\hat{x}=(\cos \theta_{\hat{x}} ,\sin \theta_{\hat{x}})\in \mathbb{S}$, $\zeta=(\cos \theta_{\zeta},\sin \theta_{\zeta})\in \mathbb{S}$ and $y=(\rho \cos \theta,\rho \sin \theta)\in \mathcal{B}_R$. Together with the trigonometric identities, we have
\begin{equation}\label{eqn:tri}
\begin{array}{ll}
 \hat{x}\cdot \zeta=\cos \theta_{\hat{x}}\cos \theta_{\zeta}+\sin \theta_{\hat{x}}\sin \theta_{\zeta}=\cos(\theta_{\hat{x}}-\theta_{\zeta}),\\
(\hat{x}-\zeta)\cdot y=\hat{x}\cdot y-\zeta\cdot y=\rho\big(\cos(\theta_{\hat{x}}-\theta)-\cos(\theta_{\zeta}-\theta)\big).
\end{array}    
\end{equation}
{Then,} we derive from \eqref{eqn:F1*Cart} and \eqref{eqn:F2*Cart} that
\begin{align*}
[F_1^{\omega *}\Lambda^{\omega}](\theta,\rho)=&-c_{\omega}\int_{[0,2\pi]^2} \cos(\theta_{\hat{x}}-\theta_{\zeta})   e^{i\omega\rho(\cos(\theta_{\hat{x}}-\theta)-\cos(\theta_{\zeta}-\theta))}\Lambda^{\omega}(\theta_{\hat{x}},\theta_{\zeta}){\rm d}\theta_{\hat{x}}{\rm d}\theta_{\zeta}\\
=&-c_{\omega}\int_0^{2\pi} e^{i\omega\rho\cos(\theta_{\hat{x}}-\theta)} \int_0^{2\pi} e^{-i\omega\rho\cos(\theta_{\zeta}-\theta)}\cos(\theta_{\hat{x}}-\theta_{\zeta})  \Lambda^{\omega}(\theta_{\hat{x}},\theta_{\zeta}){\rm d}\theta_{\zeta}{\rm d}\theta_{\hat{x}}\\
=&-c_{\omega}\int_0^{2\pi} e^{i\omega\rho\cos\theta_{\hat{x}}} \int_0^{2\pi} e^{-i\omega\rho\cos\theta_{\zeta}}\cos(\theta_{\hat{x}}-\theta_{\zeta})  \Lambda^{\omega}(\theta_{\hat{x}}+\theta,\theta_{\zeta}+\theta){\rm d}\theta_{\zeta}{\rm d}\theta_{\hat{x}},\\
[F_2^{\omega *}\Lambda^{\omega}](\theta,\rho)=&c_{\omega}\int_0^{2\pi} e^{i\omega\rho\cos\theta_{\hat{x}}} \int_0^{2\pi} e^{-i\omega\rho\cos\theta_{\zeta}}\Lambda^{\omega}(\theta_{\hat{x}}+\theta,\theta_{\zeta}+\theta){\rm d}\theta_{\zeta}{\rm d}\theta_{\hat{x}}.
\end{align*}
Then, we define $\Lambda_\theta^{\omega}(\theta_{\hat{x}},\theta_{\zeta}):=\Lambda^{\omega}(\theta_{\hat{x}}+\theta,\theta_{\zeta}+\theta)$ and the kernel $k(t,\rho)=e^{-i\omega_1 \rho \cos t}$ which gives 
\begin{align*}
[F_1^{\omega_1 *}\Lambda^{\omega_1}](\theta,\rho)
=&-c_{\omega_1}\int_0^{2\pi} \overline{k(\theta_{\hat{x}},\rho)}\int_0^{2\pi} k(\theta_{\zeta},\rho)\cos(\theta_{\hat{x}}-\theta_{\zeta})  \Lambda_\theta^{\omega_1}(\theta_{\hat{x}},\theta_{\zeta}){\rm d}\theta_{\zeta}{\rm d}\theta_{\hat{x}},\\
[F_1^{\omega_2 *}\Lambda^{\omega_2}](\theta,\frac{\omega_1}{\omega_2}\rho)=&-c_{\omega_2}\int_0^{2\pi} e^{i\omega_2\frac{\omega_1}{\omega_2}\rho\cos\theta_{\hat{x}}} \int_0^{2\pi} e^{-i\omega_2\frac{\omega_1}{\omega_2}\rho\cos\theta_{\zeta}}\cos(\theta_{\hat{x}}-\theta_{\zeta})\Lambda_\theta^{\omega_2}(\theta_{\hat{x}},\theta_{\zeta}){\rm d}\theta_{\zeta}{\rm d}\theta_{\hat{x}}\\
=&-c_{\omega_2}\int_0^{2\pi} \overline{k(\theta_{\hat{x}},\rho)} \int_0^{2\pi} k(\theta_{\zeta},\rho)\cos(\theta_{\hat{x}}-\theta_{\zeta})\Lambda_\theta^{\omega_2}(\theta_{\hat{x}},\theta_{\zeta}){\rm d}\theta_{\zeta}{\rm d}\theta_{\hat{x}}.
\end{align*}
Similarly, there hold
\begin{align*}
[F_2^{\omega_1 *}\Lambda^{\omega_1}](\theta,\rho)
=&c_{\omega_1}\int_0^{2\pi} \overline{k(\theta_{\hat{x}},\rho)}\int_0^{2\pi} k(\theta_{\zeta},\rho)\Lambda_\theta^{\omega_1}(\theta_{\hat{x}},\theta_{\zeta}){\rm d}\theta_{\zeta}{\rm d}\theta_{\hat{x}},\\
[F_2^{\omega_2 *}\Lambda^{\omega_2}](\theta,\frac{\omega_1}{\omega_2}\rho)=&c_{\omega_2}\int_0^{2\pi} \overline{k(\theta_{\hat{x}},\rho)} \int_0^{2\pi} k(\theta_{\zeta},\rho)\Lambda_\theta^{\omega_2}(\theta_{\hat{x}},\theta_{\zeta}){\rm d}\theta_{\zeta}{\rm d}\theta_{\hat{x}}.
\end{align*}
Finally, we define a linear radial scaling operator $\mathcal{S}_\omega$ such that $\mathcal{S}_\omega f(\theta,\rho):=f(\theta,\frac{\omega}{\omega_1}\rho)${. Then there holds} for $i,j=1,2$
\begin{align*}
[F_j^{\omega_i *}\Lambda^{\omega_i}](\theta,\rho)=&\mathcal{S}_{\omega_i}[F_j^{\omega_i *}\Lambda^{\omega_i}](\theta,\frac{\omega_1}{\omega_i}\rho).
\end{align*}
This completes the proof.
\end{proof}

\begin{remark}\label{rem:F*_cos}
From Proposition \ref{prop:F*_formula}, 
we have
\begin{align*}
[F_1^{\omega *}\Lambda^{\omega}](\theta,\rho)
=&-c_{\omega}\mathcal{S}_{\omega}\int_0^{2\pi} \overline{k(\theta_{\hat{x}},\rho)}\int_0^{2\pi} k(\theta_{\zeta},\rho)\cos(\theta_{\hat{x}}-\theta_{\zeta})  \Lambda_\theta^{\omega}(\theta_{\hat{x}},\theta_{\zeta}){\rm d}\theta_{\zeta}{\rm d}\theta_{\hat{x}}.
\end{align*}
We can further isolate the variables $\theta_{\hat{x}}$ and $\theta_{\zeta}$ from the interactive integrand $\cos(\theta_{\hat{x}}-\theta_{\zeta})$ by
\begin{align}\label{eqn:cos}
\cos(\theta_{\hat{x}}-\theta_{\zeta})=\frac12 (e^{i(\theta_{\hat{x}}-\theta_{\zeta})}+e^{-i(\theta_{\hat{x}}-\theta_{\zeta})})=\frac12 e^{i\theta_{\hat{x}}}e^{-i\theta_{\zeta}}+\frac12 e^{i\theta_{\zeta}}e^{-i\theta_{\hat{x}}},
\end{align}
which {converts} the operator $F_1^{\omega *}\Lambda^{\omega}$ into the average of two {double integral} operators with integral kernels $k^{\pm}(t,\rho)
:=e^{-i(\omega_1\rho\cos(t)\pm t)}$ and their complex conjugates $\overline{k^{\pm}(t,\rho)}=e^{i(\omega_1\rho\cos(t)\pm t)}$, {i.e., 
\begin{align*}
[F_1^{\omega *}\Lambda^{\omega}](\theta,\rho)
=-\frac12 c_{\omega}\mathcal{S}_{\omega}\bigg(&\int_0^{2\pi} \overline{k^+(\theta_{\hat{x}},\rho)}\int_0^{2\pi} k^+(\theta_{\zeta},\rho) \Lambda_\theta^{\omega}(\theta_{\hat{x}},\theta_{\zeta}){\rm d}\theta_{\zeta}{\rm d}\theta_{\hat{x}}\\
&+\int_0^{2\pi} \overline{k^-(\theta_{\hat{x}},\rho)}\int_0^{2\pi} k^-(\theta_{\zeta},\rho)\Lambda_\theta^{\omega}(\theta_{\hat{x}},\theta_{\zeta}){\rm d}\theta_{\zeta}{\rm d}\theta_{\hat{x}}\bigg).
\end{align*}} 
This decomposition will be used to analyze the generalization accuracy of the proposed neural networks in Section \ref{sec:gene} 
\end{remark}
\noindent
Next, we state the complementary low-rank property of the integral kernels $k$, $k^{+}$ and $k^{-}$ and establish the approximation capability of the low-rank operators to these integral kernels in the following theorem. A similar analysis can be found in \cite{CandesDemanetYing:2009}.
\begin{theorem}\label{thm:low_rank}
Let $k(t,\rho)=e^{-i\omega_1 \rho \cos t}$ and $k^{\pm}(t,\rho)=e^{-i(\omega_1\rho\cos(t)\pm t)}$.
Then for any $\epsilon>0$, $r\in\mathbb{N}$ and $c_{r,\epsilon}\in(0,e^{-1}r (e\epsilon)^{\frac{1}{r}})$, the operator $k_{r}(t,\rho)$ of rank $r$ defined as
\begin{align}\label{eqn:low_rank_func}
k_{r}(t,\rho)=\sum_{l=0}^{r-1}\frac{(-i\omega_1)^l}{l !}(\cos t-\cos t_0)^l e^{-i\omega_1\rho_0(\cos t-\cos t_0)}(\rho-\rho_0)^l e^{-i\omega_1\rho \cos t_0},
\end{align}
with any $t_0\in[0,2\pi]$ and $\rho_0\in[0,\infty)$, satisfies
\begin{align*}
|k(t,\rho)-k_{r}(t,\rho)|<\epsilon  \quad\mbox{and}\quad |k^\pm(t,\rho)-k_{r}(t,\rho)e^{\mp i t}|<\epsilon
\end{align*}
for any $(t,\rho)\in\{(t,\rho): \omega_1|\rho-\rho_0||t-t_0|\leq c_{r,\epsilon}\}$.
\end{theorem}
\begin{proof}
In order to approximate the kernel $k(t,\rho)=e^{-i\omega_1 \rho \cos t}$ by {low-rank} operators, i.e., a sum of separable operators, we decompose $k(t,\rho)$ into three components by
\begin{align*}
k(t,\rho)=e^{-i\omega_1 \rho \cos t_0}e^{-i\omega_1 \rho_0 (\cos t-\cos t_0)}e^{-i\omega_1 (\rho-\rho_0) (\cos t-\cos t_0)}.
\end{align*}
With the notation $x_{t,\rho}=-\omega_1 (\rho-\rho_0) (\cos t-\cos t_0)\in \mathbb{R}$, we apply $(r-1)$th order Taylor expansion on the last component $e^{ix_{t,\rho}}$ at the origin and derive that
\begin{align*}
e^{ix_{t,\rho}}=\sum_{l=0}^{r-1}\frac{(ix_{t,\rho})^l}{l !}+\frac{(ix_{t,\rho})^r}{r!}e^{ix^*_{t,\rho}}
\end{align*}
with some $x^*_{t,\rho}\in \mathbb{R}$.
Let $k_r(t,\rho)$ be the operator defined in \eqref{eqn:low_rank_func}{. Then there holds}
\begin{align*}
|k(t,\rho)-k_{r}(t,\rho)|=&|e^{-i\omega_1 \rho \cos t_0}e^{-i\omega_1 \rho_0 (\cos t-\cos t_0)}(e^{ix_{t,\rho}}-\sum_{l=0}^{r-1}\frac{(ix_{t,\rho})^l}{l !})|
\leq|\frac{(ix_{t,\rho})^r}{r!}e^{ix^*_{t,\rho}}|\leq \frac{|x_{t,\rho}|^r}{r!}.
\end{align*}
For any $t$ and $\rho$ such that $\omega_1|\rho-\rho_0||t-t_0|\leq c_{r,\epsilon}$, we bound $|x_{t,\rho}|$ by $$|x_{t,\rho}|=|\omega_1 (\rho-\rho_0) (\cos t-\cos t_0)|\leq \omega_1 |\rho-\rho_0||t-t_0|\leq c_{r,\epsilon}.$$
Then, together with the assumption $c_{r,\epsilon}<e^{-1}r (e\epsilon)^{\frac{1}{r}}$ and the inequality $e(\frac{r}{e})^r\leq r!$,
we obtain 
\begin{align}\label{eqn:low_rank_c}
|k(t,\rho)-k_{r}(t,\rho)|\leq\frac{|x_{t,\rho}|^r}{r!}
\leq\frac{c_{r,\epsilon}^r}{r!}
<\frac{e^{1-r}r^r}{r!}\epsilon\leq \epsilon.
\end{align}
Finally, the inequality
\begin{align*}
|k^\pm(t,\rho)-k_{r}(t,\rho)e^{\mp i t}|\leq|k(t,\rho)-k_{r}(t,\rho)|<\epsilon
\end{align*}
completes the proof.
\end{proof}
\noindent
Now, with the formula for $F^*$ in Proposition \ref{prop:F*_formula}, we {can express} $F^*F+\alpha I$ and its inverse in the following proposition.
\begin{prop}\label{prop:F*F_formula}
Let $F_1^\omega$, $F_2^\omega$ and $F$ be the operators defined in \eqref{eqn:F1F2} and \eqref{eqn:prob_linear_sys} respectively. Then, for any vector-valued function ${\bf f}\in (L^2(\mathcal{B}_R))^2$, the operator $F^*F+\alpha I$ and its inverse are  translation invariant and given by
\begin{align*}
{\bf h}:=(F^*F+\alpha I){\bf f}={\bf g}*{\bf f} \quad \mbox{and}\quad (F^*F+\alpha I)^{-1}{\bf h}=\mathcal{F}^{-1}\big(\mathcal{F}({\bf g})^{-1}\big)*{\bf h}
\end{align*}
where $\mathcal{F}$ denotes the element-wise Fourier transform and the symmetric convolution kernel
\begin{align*}
{\bf g}=&\sum_{i=1}^2\begin{bmatrix}
g_{11}^{\omega_i}+\frac\alpha2\delta & g_{12}^{\omega_i} \\
g_{21}^{\omega_i} & g_{22}^{\omega_i}+\frac\alpha2\delta
\end{bmatrix}
\end{align*}
with the Dirac delta distribution $\delta$ and the real-valued radially symmetric filters {$$g_{j{j'}}^\omega(y)=(-1)^{j+{j'}}c_\omega^2 \int_{\mathbb{S}\times\mathbb{S}} (\hat{x}\cdot \zeta)^{4-j-{j'}} \cos(\omega (\hat{x}-\zeta)\cdot y) {\rm d}\hat{x}\, {\rm d}\zeta, \quad j,j'=1,2.$$}
Furthermore, the kernel of $(F^*F+\alpha I)^{-1}$ can be written as a convolution
\begin{align*}
\mathcal{F}^{-1}\big(\mathcal{F}({\bf g})^{-1}\big)
=(g_2 I_2)*\sum_{i=1}^2\begin{bmatrix}
g_{22}^{\omega_i}+\frac\alpha2\delta &-g_{12}^{\omega_i} \\
-g_{21}^{\omega_i} &g_{11}^{\omega_i}+\frac\alpha2\delta
\end{bmatrix},
\end{align*}
where $I_2$ is the $2\times 2$ indentity matrix and $g_2:=\mathcal{F}^{-1}\big(\det(\mathcal{F}({\bf g}))^{-1}\big):=\mathcal{F}^{-1}\big(\mathcal{F}(g_{de})^{-1}\big)$ is the kernel for the deconvolution of the convolution with kernel $g_{de}:=\sum_{i,i'=1}^2(g_{11}^{\omega_i}+\frac\alpha2\delta)*(g_{22}^{\omega_{i'}}+\frac\alpha2\delta) -g_{12}^{\omega_i}*g_{12}^{\omega_{i'}}$.
\end{prop}   
\begin{proof}
Combining the definitions of the operators $F_1^\omega$ and $F_2^\omega$ in \eqref{eqn:F1F2} with the formulas for their adjoint operators obtained in \eqref{eqn:F1*Cart} and \eqref{eqn:F2*Cart} directly implies that
\begin{align*}
[F_1^{\omega *}F_1^{\omega}\gamma](y)
=&(-c_\omega)^2 \int_{\mathbb{S}\times\mathbb{S}} (\hat{x}\cdot \zeta)   e^{i\omega (\hat{x}-\zeta)\cdot y}\Big(\int_{\mathcal{B}_R}  (\hat{x}\cdot \zeta)   e^{i\omega (\zeta-\hat{x})\cdot z}\gamma(z){\rm d}z\Big){\rm d}\hat{x}\, {\rm d}\zeta\\
=&c_\omega^2    \int_{\mathcal{B}_R}  \Big(\int_{\mathbb{S}\times\mathbb{S}} (\hat{x}\cdot \zeta)^2 e^{i\omega (\hat{x}-\zeta)\cdot (y-z)} {\rm d}\hat{x}\, {\rm d}\zeta\Big)\gamma(z){\rm d}z:=    
%\int_{\mathcal{B}_R} g_{11}^\omega(y-z)\gamma(z){\rm d}z
[g_{11}^\omega*\gamma](y),\\
[F_2^{\omega *}F_1^{\omega}\gamma](y)
=&-c_\omega^2 \int_{\mathbb{S}\times\mathbb{S}}   e^{i\omega (\hat{x}-\zeta)\cdot y}\Big(\int_{\mathcal{B}_R}  (\hat{x}\cdot \zeta)   e^{i\omega (\zeta-\hat{x})\cdot z}\gamma(z){\rm d}z\Big){\rm d}\hat{x}\, {\rm d}\zeta\\
=&-c_\omega^2    \int_{\mathcal{B}_R}  \Big(\int_{\mathbb{S}\times\mathbb{S}} (\hat{x}\cdot \zeta) e^{i\omega (\hat{x}-\zeta)\cdot (y-z)} {\rm d}\hat{x}\, {\rm d}\zeta\Big)\gamma(z){\rm d}z:=   
%\int_{\mathcal{B}_R} g_{21}^\omega(y-z)\gamma(z){\rm d}z
[g_{21}^\omega*\gamma](y),\\
[F_1^{\omega *}F_2^{\omega}\eta](y)
=&-c_\omega^2 \int_{\mathbb{S}\times\mathbb{S}} (\hat{x}\cdot \zeta)   e^{i\omega (\hat{x}-\zeta)\cdot y}\Big(\int_{\mathcal{B}_R}    e^{i\omega (\zeta-\hat{x})\cdot z}\eta(z){\rm d}z\Big){\rm d}\hat{x}\, {\rm d}\zeta\\
=&-c_\omega^2    \int_{\mathcal{B}_R}  \Big(\int_{\mathbb{S}\times\mathbb{S}} (\hat{x}\cdot \zeta) e^{i\omega (\hat{x}-\zeta)\cdot (y-z)} {\rm d}\hat{x}\, {\rm d}\zeta\Big)\eta(z){\rm d}z:=    
%\int_{\mathcal{B}_R} g_{12}^\omega(y-z)\eta(z){\rm d}z,
[g_{12}^\omega*\eta](y)\\
[F_2^{\omega *}F_2^{\omega}\eta](y)
=&c_\omega^2 \int_{\mathbb{S}\times\mathbb{S}}    e^{i\omega (\hat{x}-\zeta)\cdot y}\Big(\int_{\mathcal{B}_R}e^{i\omega (\zeta-\hat{x})\cdot z}\eta(z){\rm d}z\Big){\rm d}\hat{x}\, {\rm d}\zeta\\
=&c_\omega^2    \int_{\mathcal{B}_R}  \Big(\int_{\mathbb{S}\times\mathbb{S}}  e^{i\omega (\hat{x}-\zeta)\cdot (y-z)} {\rm d}\hat{x}\, {\rm d}\zeta\Big)\eta(z){\rm d}z:=    %\int_{\mathcal{B}_R} g_{22}^\omega(y-z)\eta(z){\rm d}
[g_{22}^\omega*\eta](y),
\end{align*}
where the real-valued convolution kernels
\begin{align*}
g_{11}^\omega(y)=&c_\omega^2\int_{\mathbb{S}\times\mathbb{S}} (\hat{x}\cdot \zeta)^2 e^{i\omega (\hat{x}-\zeta)\cdot y} {\rm d}\hat{x}\, {\rm d}\zeta=c_\omega^2\int_{\mathbb{S}\times\mathbb{S}} (\hat{x}\cdot \zeta)^2 \cos(\omega (\hat{x}-\zeta)\cdot y){\rm d}\hat{x}\, {\rm d}\zeta\\
g_{22}^\omega(y)=&c_\omega^2\int_{\mathbb{S}\times\mathbb{S}}  e^{i\omega (\hat{x}-\zeta)\cdot y}{\rm d}\hat{x}\, {\rm d}\zeta=c_\omega^2\int_{\mathbb{S}\times\mathbb{S}} \cos(\omega (\hat{x}-\zeta)\cdot y){\rm d}\hat{x}\, {\rm d}\zeta,\\
\mbox{and}\quad g_{12}^\omega(y)=&g_{21}^\omega(y)=-c_\omega^2\int_{\mathbb{S}\times\mathbb{S}} (\hat{x}\cdot \zeta) e^{i\omega (\hat{x}-\zeta)\cdot y} {\rm d}\hat{x}\, {\rm d}\zeta=-c_\omega^2\int_{\mathbb{S}\times\mathbb{S}} (\hat{x}\cdot \zeta) \cos(\omega (\hat{x}-\zeta)\cdot y){\rm d}\hat{x}\, {\rm d}\zeta,
\end{align*}
which is proved to be radially symmetric in
Lemma \ref{lem:filter}.
Then by the definition of $F$ in \eqref{eqn:prob_linear_sys} and its adjoint $F^*$ given in Proposition \ref{prop:F*_formula}, there holds 
\begin{align*}
F^*F+\alpha I=\sum_{i=1}^2\begin{bmatrix}
F_1^{\omega_i *}F_1^{\omega_i}+\frac\alpha2 I &F_1^{\omega_i *}F_2^{\omega_i}\\
F_2^{\omega_i *}F_1^{\omega_i} &F_2^{\omega_i *}F_2^{\omega_i}+\frac\alpha2 I
\end{bmatrix}.
\end{align*}
Now, we rewrite the above operator applied on a vector-valued function ${\bf f}=\begin{bmatrix}
f_1\\f_2
\end{bmatrix}$ where $f_1,f_2\in L^2(\mathcal{B}_R)$ as a convolution
\begin{align*}
{\bf h}:=(F^*F+\alpha I)
{\bf f}=&\sum_{i=1}^2\begin{bmatrix}
(F_1^{\omega_i *}F_1^{\omega_i}+\frac\alpha2 I)f_1 +F_1^{\omega_i *}F_2^{\omega_i}f_2\\
F_2^{\omega_i *}F_1^{\omega_i}f_1+(F_2^{\omega_i *}F_2^{\omega_i}+\frac\alpha2 I)f_2
\end{bmatrix}
=\sum_{i=1}^2\begin{bmatrix}
(g_{11}^{\omega_i}+\frac\alpha2 \delta) *f_1+g_{12}^{\omega_i} * f_2\\
g_{21}^{\omega_i}*f_1+(g_{22}^{\omega_i}+\frac\alpha2 \delta)* f_2
\end{bmatrix}\\
=&\sum_{i=1}^2\begin{bmatrix}
g_{11}^{\omega_i}+\frac\alpha2\delta & g_{12}^{\omega_i} \\
g_{21}^{\omega_i} & g_{22}^{\omega_i}+\frac\alpha 2 \delta
\end{bmatrix}*\begin{bmatrix}
f_1\\f_2
\end{bmatrix}:={\bf g}*{\bf f},
\end{align*}
where $\delta$ is the Dirac delta distribution.
Next, we compute the inverse operator $(F^*F+\alpha I)^{-1}$ applied on ${\bf h}$. We denote the element-wise Fourier transform {by $\mathcal{F}$. Then} with the properties of the Fourier transform and the convolution, there holds
\begin{align*}
\mathcal{F}({\bf h})=\mathcal{F}({\bf g}*{\bf f})=&\sum_{i=1}^2\begin{bmatrix}
\mathcal{F}(g_{11}^{\omega_i}+\frac\alpha2\delta) \mathcal{F}(f_1)+\mathcal{F}(g_{12}^{\omega_i}) \mathcal{F}(f_2)\\
\mathcal{F}(g_{21}^{\omega_i})\mathcal{F}(f_1)+\mathcal{F}(g_{22}^{\omega_i}+\frac\alpha2\delta)\mathcal{F}(f_2)
\end{bmatrix}
=\mathcal{F}({\bf g})\mathcal{F}({\bf f}),
\end{align*}
which implies $\mathcal{F}({\bf g})^{-1}\mathcal{F}({\bf h})=\mathcal{F}({\bf f})$ where the well-defined inverse matrix-valued function
\begin{align*}
\mathcal{F}({\bf g})^{-1}
%=&\det(\mathcal{F}({\bf g}))^{-1}\begin{bmatrix}
%\sum_{i=1}^2\mathcal{F}(g_{22}^{\omega_i}+\frac\alpha2\delta) &-\sum_{i=1}^2\mathcal{F}(g_{12}^{\omega_i}) \\
%-\sum_{i=1}^2\mathcal{F}(g_{21}^{\omega_i}) &\sum_{i=1}^2\mathcal{F}(g_{11}^{\omega_i}+\frac\alpha2\delta)
%\end{bmatrix}\\
=&\det(\mathcal{F}({\bf g}))^{-1}\sum_{i=1}^2\begin{bmatrix}
\mathcal{F}(g_{22}^{\omega_i}+\frac\alpha2\delta) &-\mathcal{F}(g_{12}^{\omega_i}) \\
-\mathcal{F}(g_{21}^{\omega_i}) &\mathcal{F}(g_{11}^{\omega_i}+\frac\alpha2\delta)
\end{bmatrix},
\end{align*}
with the positive determinant $\det(\mathcal{F}({\bf g}))\geq \alpha^2$ (proved in Lemma \ref{lem:det}).
Thus, we have $(F^*F+\alpha I)^{-1}{\bf h}={\bf f}=\mathcal{F}^{-1}\Big(\mathcal{F}({\bf g})^{-1}\mathcal{F}({\bf h})\Big)=\mathcal{F}^{-1}\big(\mathcal{F}({\bf g})^{-1}\big)*{\bf h}$ and
\begin{align*}
\mathcal{F}^{-1}\big(\mathcal{F}({\bf g})^{-1}\big)=&\Big(\mathcal{F}^{-1}\big(\det(\mathcal{F}({\bf g}))^{-1}\big)I_2\Big)*\mathcal{F}^{-1}\big(\sum_{i=1}^2\begin{bmatrix}
\mathcal{F}(g_{22}^{\omega_i}+\frac\alpha2\delta) &-\mathcal{F}(g_{12}^{\omega_i}) \\
-\mathcal{F}(g_{21}^{\omega_i}) &\mathcal{F}(g_{11}^{\omega_i}+\frac\alpha2\delta)
\end{bmatrix}\big)\\
=&\Big(\mathcal{F}^{-1}\big(\det(\mathcal{F}({\bf g}))^{-1}\big)I_2\Big)*\sum_{i=1}^2\begin{bmatrix}
g_{22}^{\omega_i}+\frac\alpha2\delta &-g_{12}^{\omega_i} \\
-g_{21}^{\omega_i} &g_{11}^{\omega_i}+\frac\alpha2\delta
\end{bmatrix}.
\end{align*}
%where $I_2$ is the $2\times 2$ indentity matrix and $\mathcal{F}^{-1}\big(\det(\mathcal{F}({\bf g}))^{-1}\big)$ is the kernel for the deconvolution of the convolution with kernel $g_{de}=\sum_{i,i'=1}^2(g_{11}^{\omega_i}+\frac\alpha2\delta)*(g_{22}^{\omega_{i'}}+\frac\alpha2\delta) -g_{12}^{\omega_i}*g_{12}^{\omega_{i'}}$.
Furthermore, the determinant $\det(\mathcal{F}({\bf g}))$ can be rewritten as
\begin{align*}
\det(\mathcal{F}({\bf g}))=&\sum_{i=1}^2\mathcal{F}(g_{11}^{\omega_i}+\frac\alpha2\delta) \sum_{i=1}^2\mathcal{F}(g_{22}^{\omega_i}+\frac\alpha2\delta) -(\sum_{i=1}^2\mathcal{F}(g_{12}^{\omega_i}))^2\\
=&\sum_{i,i'=1}^2\mathcal{F}\big((g_{11}^{\omega_i}+\frac\alpha2\delta)*(g_{22}^{\omega_{i'}}+\frac\alpha2\delta)\big) -\mathcal{F}(g_{12}^{\omega_i}*g_{12}^{\omega_{i'}})\\
=&\mathcal{F}\Big(\sum_{i,i'=1}^2(g_{11}^{\omega_i}+\frac\alpha2\delta)*(g_{22}^{\omega_{i'}}+\frac\alpha2\delta) -g_{12}^{\omega_i}*g_{12}^{\omega_{i'}}\Big).
\end{align*}
This completes the proof.
\end{proof}
\begin{remark}\label{rem:real_out}
By the formula in Proposition \ref{prop:F*F_formula} (when $\alpha=0$) that $F^*F{\bf f}={\bf g}*{\bf f}$ where ${\bf g}$ is a real-valued function, we obtain that $F^*\Lambda=F^*F \begin{bmatrix}
\gamma\\
\eta    
\end{bmatrix}={\bf g}*\begin{bmatrix}
\gamma\\
\eta    
\end{bmatrix}$ is real-valued.
\end{remark}

%%%%%%%%%%%%%%%%%%%%%%%%%%%%%%%%%%%%%%%%%%%%%%%%%%%%%%%%%%%%
\section{Combined neural networks}\label{sec:NN_approx}
In this section, by noting the properties of $F^*$ and $(F^*F+\alpha I)^{-1}$ {in Section \ref{sec:prop}}, we construct a combined neural network to approximate the regularized pseudo-inverse $$F_\alpha^\dag=(F^*F+\alpha I)^{-1}F^*: \left(L^2({\mathbb S}\times {\mathbb S})\right)^2\to \left(L^2(\mathcal{B}_R)\right)^2.$$

\subsection{Butterfly neural networks (BFNNs) for $F^*$}
Now, we construct a fully connected neural network with or without the butterfly factorization to approximate the adjoint operator $F^*:\left(L^2({\mathbb S}\times {\mathbb S})\right)^2\to \left(L^2(\mathcal{B}_R)\right)^2$ . 
We shall start with the $0$-level BFNN (uncompressed, no butterfly structure) that discretizes the integrals in Proposition \ref{prop:F*_formula} as matrix multiplications. 
After converting the Cartesian coordinates to the polar coordinates, we discretize the domain of the input $\Lambda(\theta_{\hat{x}},\theta_{\zeta})$, i.e., $[0,2\pi]^2$, uniformly with $\{(\theta_{\hat{x}_i},\theta_{\zeta_j})=(\frac{2\pi i}{n_\theta},\frac{2\pi j}{n_\theta})\}_{i,j=1}^{n_\theta}$. 
We consider the compact domains $[0,2\pi]\times[0,\frac{\omega_2}{\omega_1}R]$ of the integral kernel $k(t,\rho)$ and $[0,2\pi]\times[0,R]$ of the output $[F^*\Lambda](\theta,\rho)$ such that $\bar{\Omega} \subset [0,2\pi]\times[0,R]$. 
Without loss of generality, we assume that $R=1$ and $\frac{\omega_2}{\omega_1}\in \mathbb{N}$, and then discretize $[0,2\pi]\times[0,\frac{\omega_2}{\omega_1}]$ and $[0,2\pi]\times[0,1]$ uniformly with $\{(t_i,\rho_j)=(\frac{2\pi i}{n_\theta},\frac{\omega_2}{\omega_1}\frac{j}{n_\rho})\}_{i,j=1}^{i=n_\theta,j=n_\rho}$ and $\{(\theta_i,\rho_j)=(\frac{2\pi i}{n_\theta},\frac{\omega_2}{\omega_1}\frac{j}{n_\rho})\}_{i,j=1}^{i=n_\theta,j=\frac{\omega_1}{\omega_2}n_\rho}$, respectively, where $\frac{\omega_1}{\omega_2}n_\rho\in\mathbb{N}$ for some $n_\theta\leq n_\rho\in\mathbb{N}$. 
{In the rest of this section, the input $\Lambda=\begin{bmatrix}
\Lambda^{\omega_1}\\
\Lambda^{\omega_2}    
\end{bmatrix}$ and output $F^*\Lambda=
\sum_{i=1}^2\begin{bmatrix}
F_1^{\omega_i *}\Lambda^{\omega_i}\\
F_2^{\omega_i *}\Lambda^{\omega_i}
\end{bmatrix}$ refer to either functions or their discretization matrices with
\begin{align*}
\Lambda^{\omega}=&\big(\Lambda^\omega[i,j]\big)_{i,j=1}^{n_\theta}:=\big(\Lambda^{\omega}(\theta_{\hat{x}_i},\theta_{\zeta_j})\big)_{i,j=1}^{n_\theta}\in\mathbb{C}^{n_\theta\times n_\theta} \quad \mbox{and}\\
F_{\tilde{i}}^{\omega *}\Lambda^{\omega}=&\big((F_{\tilde{i}}^{\omega *}\Lambda^{\omega})[i,j]\big)_{i,j=1}^{i=n_\theta,j=\frac{\omega_1}{\omega_2}n_\rho}:=\big([F_{\tilde{i}}^{\omega *}\Lambda^{\omega}](\theta_i,\rho_j)\big)_{i,j=1}^{i=n_\theta,j=\frac{\omega_1}{\omega_2}n_\rho} \in \mathbb{R}^{n_\theta\times\frac{\omega_1}{\omega_2}n_\rho}, \quad \mbox{for } \tilde{i}=1, 2,
\end{align*}
depending on the context.}
With the above discretization, we define the discretized integrands as matrices given by
\begin{align}
\Lambda_\theta^\omega=&\big(\Lambda_\theta^\omega[i,j]\big)_{i,j=1}^{n_\theta}:=\big(\Lambda_\theta^{\omega}(\theta_{\hat{x}_i},\theta_{\zeta_j})\big)_{i,j=1}^{n_\theta}\in\mathbb{C}^{n_\theta\times n_\theta},\nonumber\\
K=&\big(K[i,j]\big)_{i,j=1}^{i=n_\theta,j=n_\rho}:=\big(k(t_i,\rho_j)\big)_{i,j=1}^{i=n_\theta,j=n_\rho}\in\mathbb{C}^{n_\theta\times n_\rho},\label{eqn:K}\\
C=&\big(C[i,j]\big)_{i,j=1}^{n_\theta}:=\big(-\cos(\theta_{\hat{x}_i}-\theta_{\zeta_j})\big)_{i,j=1}^{n_\theta}=\big(-\cos(\tfrac{2\pi(i-j)}{n_\theta})\big)_{i,j=1}^{n_\theta}\in\mathbb{R}^{n_\theta\times n_\theta},\label{eqn:C}
\end{align}
where the discretized shifted far-field pattern $\Lambda^\omega_\theta$ satisfies $$\Lambda^\omega_{\theta_m}[i,j]=\Lambda^{\omega}(\theta_{\hat{x}_i}+\theta_m,\theta_{\zeta_j}+\theta_m)=\Lambda^\omega[i+m,j+m],\quad \forall \theta_m=\frac{2\pi m}{n_\theta}, \;m=1,\cdots,n_\theta$$
and $C$ is a circulant matrix.
Further, we denote the discretization of $[F^*\Lambda](\theta,\rho)$ on the domain $[0,2\pi]\times[0,1]$ by $\Phi$.
Then, we construct of a neural network $\phi^0$ with 4 channels to approximate $F^*$ in the following theorem. 
\begin{theorem}\label{thm:BF_0_app}
Let $F^*$ be the adjoint operator of the linearized forward operator $F$ defined in \eqref{eqn:prob_linear_sys}. Then, there exists a neural network $\phi^0$ on the input $\Lambda$ given by
\begin{align*}
\phi^0(\Lambda)
=\sum_{i=1}^2 (P_u^T (\phi_1^0 R_1^{\omega_i}(\Lambda),\cdots,\phi_1^0 R_{n_\theta}^{\omega_i}(\Lambda) )^T P_{\omega_i}+ P_l^T (\phi_2^0 R_1^{\omega_i}(\Lambda),\cdots,\phi_2^0 R_{n_\theta}^{\omega_i}(\Lambda))^T P_{\omega_i})
\end{align*}
with the shifted inputs $\{R_m^{\omega}(\Lambda)\}_{m=1}^{n_\theta}$ defined in \eqref{eqn:rotaton}, explicit matrix transformations $P_u$, $P_l$, $P_{\omega}$ defined in {\eqref{eqn:P_ul} and \eqref{eqn:P_omega}, and linear operators $\phi_1^0$, $\phi_2^0$ defined in \eqref{eqn:phi_1^0} and \eqref{eqn:phi_2^0}}, such that
\begin{align*}
|\Phi-\phi^0(\Lambda)|_{\infty}\leq c_{{\bf bf}_0} n_\theta^{-1},
\end{align*}
with the constant $c_{{\bf bf}_0}=8\pi^3 (4\omega_2 +3)(c_{\omega_1}^2+c_{\omega_2}^2)(\|\eta(y)\|_{L^1(\mathcal{B}_R)}+ \|\gamma(y)\|_{L^1(\mathcal{B}_R)})$ and {$|\cdot|_{\infty}$ denoting the maximum magnitude of the entries in the vector or matrix.}
\end{theorem}
\begin{proof}
We decompose the operator $F^*$ applied on the input $\Lambda$ into three parts:
\begin{align*}
\Lambda=\begin{bmatrix}
\Lambda^{\omega_1}\\
\Lambda^{\omega_2}    
\end{bmatrix} \xrightarrow{\;(i)\;} \Big\{
\Lambda_{\theta_m}^{\omega_1},
\Lambda_{{\theta_m}}^{\omega_2}    
\Big\}_{m=1}^{n_\theta} \xrightarrow{\;(ii)\;} 
\Big\{
F_j^{\omega_i *}\Lambda^{\omega_i}  
\Big\}_{i,j=1}^{2}  \xrightarrow{\;(iii)\;}  F^*\Lambda=
\sum_{i=1}^2\begin{bmatrix}
F_1^{\omega_i *}\Lambda^{\omega_i}\\
F_2^{\omega_i *}\Lambda^{\omega_i}
\end{bmatrix}
\end{align*}
and simulate these three processes one by one.

\vspace{5mm}
\noindent$\textbf{(i)}$
By the discretized integrands defined above, we first define the permutation matrices for any $m\in\{1,\cdots,n_\theta\}$ as
\begin{align}\label{eqn:P_theta}
P_{1,\theta_m}=
\begin{bmatrix}
O_{(n_\theta-m)\times m}&I_{n_\theta-m}\\
I_{m}&O_{m\times(n_\theta-m)}
\end{bmatrix}\quad \mbox{and}\quad 
P_{2,\theta_m}=
\begin{bmatrix}
O_{m\times(n_\theta-m)}&I_{m}\\
I_{n_\theta-m}&O_{(n_\theta-m)\times m}
\end{bmatrix}
\end{align}
with $O$ and $I$ denoting the zero matrix and identity matrix respectively, such that $$\Lambda_{\theta_m}^\omega=P_{1,\theta_m}\Lambda^\omega P_{2,\theta_m}.$$
{Then by letting}
\begin{align}\label{eqn:P_ul}
P_u=[I_{n_\theta} \; O_{n_\theta\times n_\theta}]\quad \mbox{and}\quad P_l=[O_{n_\theta\times n_\theta}\;I_{n_\theta}]
\end{align}
we have 
\begin{align}\label{eqn:rotaton}
\Lambda_{\theta_m}^{\omega_1}%=P_{1,\theta_m}\Lambda^\omega P_{2,\theta_m}
=P_{1,\theta_m}P_u\Lambda P_{2,\theta_m}:= R_m^{\omega_1}(\Lambda)\quad \mbox{and}\quad  \Lambda_{\theta_m}^{\omega_2}=P_{1,\theta_m}P_l\Lambda P_{2,\theta_m}:= R_m^{\omega_2}(\Lambda).
\end{align}

\vspace{5mm}
\noindent$\textbf{(ii)}$
Next, we approximate the {double integral} operator $F_1^{\omega *}\Lambda^{\omega}$ at $(\theta_m,\frac{\omega_1}{\omega}\rho_n)$ for any $m\in\{1,\cdots,n_\theta\}$ and $n\in\{1,\cdots,n_\rho\}$ given in Proposition \ref{prop:F*_formula}, i.e.,
\begin{align*}
[F_1^{\omega *}\Lambda^{\omega}](\theta_m,\frac{\omega_1}{\omega}\rho_n)=-c_{\omega}\int_0^{2\pi} \overline{k(\theta_{\hat{x}},\rho_n)}\int_0^{2\pi} k(\theta_{\zeta},\rho_n)\cos(\theta_{\hat{x}}-\theta_{\zeta})  \Lambda_{\theta_m}^{\omega}(\theta_{\hat{x}},\theta_{\zeta}){\rm d}\theta_{\zeta}{\rm d}\theta_{\hat{x}}
\end{align*} 
by
\begin{align}\label{eqn:approx_F*}
(\phi_1^0 \Lambda^\omega_{\theta_m})[n]:=c_{\omega}(\frac{2\pi}{n_\theta})^2\sum_{i=1}^{n_\theta} \overline{K[i,n]}\sum_{j=1}^{n_\theta} K[j,n] C[i,j]\Lambda_{\theta_m}^{\omega}[i,j],
\end{align}
{with matrices $K$ and $C$ defined in \eqref{eqn:K} and \eqref{eqn:C}, respectively.}
It follows directly with $\odot$ denoting the element-wise Hadamard multiplication that
\begin{align*}
(\phi_1^0 \Lambda^\omega_{\theta_m})[n]=& \frac{4\pi^2 c_{\omega}}{n_\theta^2}\sum_{i=1}^{n_\theta} \overline{K[i,n]} \big((C\odot\Lambda_{\theta_m}^{\omega}) K\big)_{in}
=\frac{4\pi^2 c_{\omega}}{n_\theta^2}\sum_{i=1}^{n_\theta} \big(\overline{K} \odot((C\odot\Lambda_{\theta_m}^{\omega}) K)\big)_{in}\\
=&\frac{4\pi^2 c_{\omega}}{n_\theta^2}\Big(b^T\big(\overline{K} \odot((C\odot\Lambda_{\theta_m}^{\omega}) K)\big)\Big)_{1n},
\end{align*}
where $b=[1,\cdots,1]^T\in \mathbb{R}^{n_\theta}$. For further decomposition and approximation analysis (see Theorem \ref{thm:BF_1_app}), we rewrite the above entries {as}
\begin{align*}
(\phi_1^0 \Lambda^\omega_{\theta_m})[n]
=&\frac{4\pi^2 c_{\omega}}{n_\theta^2}\big(K^* (C\odot\Lambda_{\theta_m}^{\omega}) K\big)_{nn},
\end{align*}
with $K^*=\overline{K}^T$ denoting the conjugate transpose of $K$.
Now, we define the intermediate output matrix $O_{1,\omega}^0=\big(O_{1,\omega}^0[m,n]\big)_{m,n=1}^{m=n_\theta,n=\frac{\omega_1}{\omega_2}n_\rho}\in\mathbb{C}^{n_\theta\times(\frac{\omega_1}{\omega_2}n_\rho)}$, by selecting corresponding columns from the matrix $\big((\phi_1^0 \Lambda^\omega_{\theta_m})[n]\big)_{m,n=1}^{m=n_\theta,n=n_\rho}$ with the identity $\frac{\omega}{\omega_1}\rho_n=\rho_{\frac{\omega}{\omega_1}n}$, as
\begin{align*}
O_{1,\omega}^0[m,n]= (\phi_1^0 \Lambda^\omega_{\theta_m})[\frac{\omega}{\omega_1}n]\approx [F_1^{\omega *}\Lambda^{\omega}](\theta_m,\frac{\omega_1}{\omega}\rho_{\frac{\omega}{\omega_1}n})=[F_1^{\omega *}\Lambda^{\omega}](\theta_m,\rho_{n}).
\end{align*}
This selection process can be represented by the matrix multiplication 
\begin{align*}
O_{1,\omega}^0=(\phi_1^0\Lambda_{\theta_1}^\omega,\cdots,\phi_1^0\Lambda_{\theta_{n_\theta}}^\omega)^T P_\omega,
\end{align*}
where the matrix $P_\omega\in\mathbb{R}^{n_\rho \times (\frac{\omega_1}{\omega_2}n_\rho)}$ is defined as
\begin{align}\label{eqn:P_omega}
P_\omega[i,j]=\left\{\begin{array}{ll}
    1 &  \qquad\mbox{if } i=\frac{\omega}{\omega_1}j,\;j=1,\cdots,\frac{\omega_1}{\omega_2}n_\rho\\
    0 & \qquad\mbox{otherwise}
\end{array}
\right.
\end{align}
and the row vectors $\phi_1^0\Lambda_{\theta_m}^\omega=((\phi_1^0 \Lambda^\omega_{\theta_m})[1],\cdots,(\phi_1^0 \Lambda^\omega_{\theta_m})[n_\rho])\in \mathbb{R}^{1\times n_\rho}$ is given by
\begin{align}\label{eqn:phi_1^0}
\phi_1^0\Lambda_{\theta_m}^\omega=\frac{4\pi^2 c_{\omega}}{n_\theta^2}b^T\big(\overline{K} \odot((C\odot\Lambda_{\theta_m}^{\omega}) K)\big) \quad \mbox{or} \quad \phi_1^0\Lambda_{\theta_m}^\omega=\frac{4\pi^2 c_{\omega}}{n_\theta^2}{\rm diag}\big(K^* (C\odot\Lambda_{\theta_m}^{\omega}) K\big)
\end{align}
with ${\rm diag}(\cdot)$ denoting the row vector generated sequentially by the diagonal entries.

Similarly, by replacing $C$ with the $n_\theta \times n_\theta$ matrix in which all entries are ones, we define 
\begin{align*}
O_{2,\omega}^0=(\phi_2^0\Lambda_{\theta_1}^\omega,\cdots,\phi_2^0\Lambda_{\theta_{n_\theta}}^\omega)^T P_\omega,
\end{align*}
where $\phi_2^0\Lambda_{\theta_m}^\omega\in \mathbb{R}^{1\times n_\rho}$ is given by
\begin{align}\label{eqn:phi_2^0}
\phi_2^0\Lambda_{\theta_m}^\omega=\frac{4\pi^2 c_{\omega}}{n_\theta^2}b^T\big(\overline{K} \odot(\Lambda_{\theta_m}^{\omega} K)\big) \quad \mbox{or} \quad \phi_2^0\Lambda_{\theta_m}^\omega=\frac{4\pi^2 c_{\omega}}{n_\theta^2}{\rm diag}\big(K^* \Lambda_{\theta_m}^{\omega} K\big)
\end{align}
to approximate $\big([F_2^{\omega *}\Lambda^{\omega}](\theta_m,\rho_n)\big)_{m,n=1}^{m=n_\theta,n=\frac{\omega_1}{\omega_2}n_\rho}$.

\vspace{5mm}
\noindent$\textbf{(iii)}$
Finally, using the intermediate output matrices, we derive a discretization of $F^*\Lambda=
\sum_{i=1}^2\begin{bmatrix}
F_1^{\omega_i *}\Lambda^{\omega_i}\\
F_2^{\omega_i *}\Lambda^{\omega_i}
\end{bmatrix}$
by
\begin{align}\label{eqn:approx_F*_2}
\sum_{i=1}^2\begin{bmatrix}
\big([F_1^{\omega_i *}\Lambda^{\omega_i}](\theta_m,\rho_n)\big)_{m,n=1}^{m=n_\theta,n=\frac{\omega_1}{\omega_2}n_\rho}\vspace{3mm}\\
\big([F_2^{\omega_i *}\Lambda^{\omega_i}](\theta_m,\rho_n)\big)_{m,n=1}^{m=n_\theta,n=\frac{\omega_1}{\omega_2}n_\rho}\vspace{2mm}
\end{bmatrix}
\approx 
\sum_{i=1}^2\begin{bmatrix}
O_{1,\omega_i}^0\\
O_{2,\omega_i}^0
\end{bmatrix}=\sum_{i=1}^2 (P_u^T O_{1,\omega_i}^0+ P_l^T O_{2,\omega_i}^0).
\end{align}

\vspace{5mm}
Now, we shall summarize the structure of the uncompressed $0$-level BFNN $\phi^0$ and discuss the approximation error  of $\phi^0$ in $|\cdot|_\infty$ norm.
Combining the above approximations gives 
\begin{align*}
\phi^0(\Lambda)=&\sum_{i=1}^2 (P_u^T O_{1,\omega_i}^0+ P_l^T O_{2,\omega_i}^0)
=\sum_{i=1}^2 (P_u^T (\phi_1^0\Lambda_{\theta_1}^{\omega_i},\cdots,\phi_1^0\Lambda_{\theta_{n_\theta}}^{\omega_i})^T P_{\omega_i}+ P_l^T (\phi_2^0\Lambda_{\theta_1}^{\omega_i},\cdots,\phi_2^0\Lambda_{\theta_{n_\theta}}^{\omega_i})^T P_{\omega_i})\\
=&\sum_{i=1}^2 (P_u^T (\phi_1^0 R_1^{\omega_i}(\Lambda),\cdots,\phi_1^0 R_{n_\theta}^{\omega_i}(\Lambda) )^T P_{\omega_i}+ P_l^T (\phi_2^0 R_1^{\omega_i}(\Lambda),\cdots,\phi_2^0 R_{n_\theta}^{\omega_i}(\Lambda))^T P_{\omega_i})
\end{align*}
which represents a neural network with 4 channels (corresponding to $\omega_i$ and $\phi_j^0$ for $i,j=1,2$) and 3 or 4 weight layers per channel on $2n_\theta$ shifted inputs $\{R_{m}^{\omega_1}(\Lambda),R_{m}^{\omega_2}(\Lambda)\}_{m=1}^{n_\theta}$ followed by explicit matrix transformations and {an added layer}.
With the approximations in \eqref{eqn:approx_F*} and \eqref{eqn:approx_F*_2}, we bound the element-wise difference between the discretization of $F^*\Lambda$ (i.e., $\Phi$) and the output $\phi^0(\Lambda)$ by
\begin{align}
|\Phi-\phi^0(\Lambda)|_{\infty}:=&\max_{m,n}\big|\Phi[m,n]-\phi^0(\Lambda)[m,n]\big|\leq \max_{m,n,\tilde{j}}\sum_{\tilde{i}=1}^2\big|[F_{\tilde{j}}^{\omega_{\tilde{i}} *}\Lambda^{\omega_{\tilde{i}}}](\theta_m,\rho_n)-O_{\tilde{j},\omega_{\tilde{i}}}^0[m,n]\big|\nonumber\\
\leq & \max_{m,n,\tilde{j}}\sum_{\tilde{i}=1}^2\big|[F_{\tilde{j}}^{\omega_{\tilde{i}} *}\Lambda^{\omega_{\tilde{i}}}](\theta_m,\frac{\omega_1}{\omega_{\tilde{i}}}\rho_n)-(\phi_1^0 \Lambda^{\omega_{\tilde{i}}}_{\theta_m})[n]\big|:=\max_{m,n,\tilde{j}}\sum_{\tilde{i}=1}^2 {\rm I}_{\tilde{i}\tilde{j}}^0.\label{eqn:bf_bd0}
\end{align}
We shall bound ${\rm I}_{\tilde{i}1}^0$ and ${\rm I}_{\tilde{i}2}^0$ one by one. First, by the formulas of $F_{\tilde{j}}^{\omega_{\tilde{i}} *}\Lambda^{\omega_{\tilde{i}}}$ in Proposition \ref{prop:F*_formula} and the definition of $(\phi_1^0 \Lambda^{\omega_{\tilde{i}}}_{\theta_m})[n]$ in \eqref{eqn:approx_F*}, we have
\begin{align*}
{\rm I}_{\tilde{i}1}^0=&\big|[F_1^{\omega_{\tilde{i}} *}\Lambda^{\omega_{\tilde{i}}}](\theta_m,\frac{\omega_1}{\omega_{\tilde{i}}}\rho_n)-(\phi_1^0 \Lambda^{\omega_{\tilde{i}}}_{\theta_m})[n]\big|\\
\leq& c_{\omega_{\tilde{i}}}\big|-\int_0^{2\pi} \overline{k(\theta_{\hat{x}},\rho_n)}\int_0^{2\pi} k(\theta_{\zeta},\rho_n)\cos(\theta_{\hat{x}}-\theta_{\zeta})  \Lambda_{\theta_m}^{\omega_{\tilde{i}}}(\theta_{\hat{x}},\theta_{\zeta}){\rm d}\theta_{\zeta}{\rm d}\theta_{\hat{x}}\\
&\quad\;\;-\frac{4\pi^2}{n_\theta^2}\sum_{i=1}^{n_\theta} \overline{K[i,n]}\sum_{j=1}^{n_\theta} K[j,n] C[i,j]\Lambda_{\theta_m}^{\omega_{\tilde{i}}}[i,j]\big|\\
\leq& c_{\omega_{\tilde{i}}}\sum_{i,j=1}^{n_\theta} \big|\int_{\frac{2\pi (i-1)}{n_\theta}}^{\frac{2\pi i}{n_\theta}}\int_{\frac{2\pi (j-1)}{n_\theta}}^{\frac{2\pi j}{n_\theta}} \big(\overline{k(\theta_{\hat{x}},\rho_n)}  k(\theta_{\zeta},\rho_n)\cos(\theta_{\hat{x}}-\theta_{\zeta})  \Lambda_{\theta_m}^{\omega_{\tilde{i}}}(\theta_{\hat{x}},\theta_{\zeta})\\
&\quad\qquad\qquad\qquad\qquad\qquad\qquad\qquad+\overline{K[i,n]} K[j,n] C[i,j]\Lambda_{\theta_m}^{\omega_{\tilde{i}}}[i,j]\big){\rm d}\theta_{\zeta}{\rm d}\theta_{\hat{x}}\big|\\
\leq& c_{\omega_{\tilde{i}}}\sum_{i,j=1}^{n_\theta} \int_{\frac{2\pi (i-1)}{n_\theta}}^{\frac{2\pi i}{n_\theta}}\int_{\frac{2\pi (j-1)}{n_\theta}}^{\frac{2\pi j}{n_\theta}} \big|\overline{k(\theta_{\hat{x}},\rho_n)}  k(\theta_{\zeta},\rho_n)\cos(\theta_{\hat{x}}-\theta_{\zeta})  \Lambda_{\theta_m}^{\omega_{\tilde{i}}}(\theta_{\hat{x}},\theta_{\zeta})\\
&\qquad\qquad\qquad\qquad\qquad-\overline{k(\tfrac{2\pi i}{n_\theta},\rho_n)} k(\tfrac{2\pi j}{n_\theta},\rho_n) \cos(\tfrac{2\pi i}{n_\theta}-\tfrac{2\pi j}{n_\theta})  \Lambda_{\theta_m}^{\omega_{\tilde{i}}}(\tfrac{2\pi i}{n_\theta},\tfrac{2\pi j}{n_\theta})\big|{\rm d}\theta_{\zeta}{\rm d}\theta_{\hat{x}}.
\end{align*}
By decomposing the integrand into
\begin{align*}
&\big|\overline{k(\theta_{\hat{x}},\rho_n)}  k(\theta_{\zeta},\rho_n)\cos(\theta_{\hat{x}}-\theta_{\zeta})  \Lambda_{\theta_m}^{\omega_{\tilde{i}}}(\theta_{\hat{x}},\theta_{\zeta})-\overline{k(\tfrac{2\pi i}{n_\theta},\rho_n)} k(\tfrac{2\pi j}{n_\theta},\rho_n) \cos(\tfrac{2\pi i}{n_\theta}-\tfrac{2\pi j}{n_\theta})  \Lambda_{\theta_m}^{\omega_{\tilde{i}}}(\tfrac{2\pi i}{n_\theta},\tfrac{2\pi j}{n_\theta})\big|\\
\leq&\big|(\overline{k(\theta_{\hat{x}},\rho_n)} -\overline{k(\tfrac{2\pi i}{n_\theta},\rho_n)} ) k(\theta_{\zeta},\rho_n)\cos(\theta_{\hat{x}}-\theta_{\zeta})  \Lambda_{\theta_m}^{\omega_{\tilde{i}}}(\theta_{\hat{x}},\theta_{\zeta})\big|\\
&+\big|\overline{k(\tfrac{2\pi i}{n_\theta},\rho_n)}(k(\theta_{\zeta},\rho_n)-k(\tfrac{2\pi j}{n_\theta},\rho_n))\cos(\theta_{\hat{x}}-\theta_{\zeta})  \Lambda_{\theta_m}^{\omega_{\tilde{i}}}(\theta_{\hat{x}},\theta_{\zeta})\big|\\
&+\big|\overline{k(\tfrac{2\pi i}{n_\theta},\rho_n)}k(\tfrac{2\pi j}{n_\theta},\rho_n)(\cos(\theta_{\hat{x}}-\theta_{\zeta})-\cos(\tfrac{2\pi i}{n_\theta}-\tfrac{2\pi j}{n_\theta}) )  \Lambda_{\theta_m}^{\omega_{\tilde{i}}}(\theta_{\hat{x}},\theta_{\zeta})\big|\\
&+\big|\overline{k(\tfrac{2\pi i}{n_\theta},\rho_n)}k(\tfrac{2\pi j}{n_\theta},\rho_n)\cos(\tfrac{2\pi i}{n_\theta}-\tfrac{2\pi j}{n_\theta}) (\Lambda_{\theta_m}^{\omega_{\tilde{i}}}(\theta_{\hat{x}},\theta_{\zeta})
- \Lambda_{\theta_m}^{\omega_{\tilde{i}}}(\tfrac{2\pi i}{n_\theta},\tfrac{2\pi j}{n_\theta}))\big|,
\end{align*}
we further bound ${\rm I}_{\tilde{i}1}^0$ by
\begin{align*}
{\rm I}_{\tilde{i}1}^0\leq& c_{\omega_{\tilde{i}}}\frac{4\pi^2}{n_\theta^2}\sum_{i,j=1}^{n_\theta} \big(w_{\overline{k(\cdot,\rho_n)}}(\tfrac{2\pi}{n_\theta})\|k\|_{L^{\infty}}\|\Lambda^{\omega_{\tilde{i}}}\|_{L^{\infty}}+w_{k(\cdot,\rho_n)}(\tfrac{2\pi}{n_\theta})\|\bar{k}\|_{L^{\infty}}\|\Lambda^{\omega_{\tilde{i}}}\|_{L^{\infty}}\\
&+w_{\cos(\cdot)}(\tfrac{2\pi}{n_\theta})\|\bar{k}\|_{L^{\infty}}\|k\|_{L^{\infty}}\|\Lambda^{\omega_{\tilde{i}}}\|_{L^{\infty}}+\|\bar{k}\|_{L^{\infty}}\|k\|_{L^{\infty}} w_{\Lambda^{\omega_{\tilde{i}}}}(\tfrac{2\pi}{n_\theta})\big),
\end{align*}
where $w_{f}$ is the modulus of continuity of a function $f$ given by $w_{f}(s)=\sup_{|x-y|_{\infty}\leq s}|f(x)-f(y)|$.
Then, together with the estimates in Lemma \ref{lem:est1} (when $R=1$) and the inequalities $\rho_n\leq \frac{\omega_2}{\omega_1}$ and $$w_{\cos(\cdot)}(s)\leq\sup_{|t-t'|\leq s}|\cos(t)-\cos(t')|\leq \sup_{|t-t'|\leq s}|t-t'|\leq s,$$ 
there holds
\begin{align*}
{\rm I}_{\tilde{i}1}^0\leq& c_{\omega_{\tilde{i}}}\frac{4\pi^2}{n_\theta^2}\sum_{i,j=1}^{n_\theta} \big(2w_{k(\cdot,\rho_n)}(\tfrac{2\pi}{n_\theta})\|\Lambda^{\omega_{\tilde{i}}}\|_{L^{\infty}}+w_{\cos(\cdot)}(\tfrac{2\pi}{n_\theta})\|\Lambda^{\omega_{\tilde{i}}}\|_{L^{\infty}}+ w_{\Lambda^{\omega_{\tilde{i}}}}(\tfrac{2\pi}{n_\theta})\big)\\
\leq &c_{\omega_{\tilde{i}}}\frac{4\pi^2}{n_\theta^2}\sum_{i,j=1}^{n_\theta} \big((2\omega_1 \rho_n+1+2\omega_{\tilde{i}}) \|\eta(y)\|_{L^1(\mathcal{B}_R)}+(2\omega_1 \rho_n+3+2\omega_{\tilde{i}}) \|\gamma(y)\|_{L^1(\mathcal{B}_R)}\big)c_{\omega_{\tilde{i}}}\frac{2\pi}{n_\theta}\\
\leq &8\pi^3 c_{\omega_{\tilde{i}}}^2\big((2\omega_2 +1+2\omega_{\tilde{i}}) \|\eta(y)\|_{L^1(\mathcal{B}_R)}+(2\omega_2 +3+2\omega_{\tilde{i}}) \|\gamma(y)\|_{L^1(\mathcal{B}_R)}\big) n_\theta^{-1}.
\end{align*}
Similarly, for the second term ${\rm I}_{\tilde{i}2}^0$ where the function $-\cos(\theta_{\hat{x}}-\theta_{\zeta})$ in the integrand of $F_{\tilde{j}}^{\omega_{\tilde{i}} *}\Lambda^{\omega_{\tilde{i}}}$ and the matrix $C$ in the summand of $(\phi_1^0 \Lambda^{\omega_{\tilde{i}}}_{\theta_m})[n]$ are replaced by $1$ and an $n_\theta \times n_\theta$ matrix in which all entries are ones respectively, we obtain that
\begin{align*}
{\rm I}_{\tilde{i}2}^0\leq& c_{\omega_{\tilde{i}}}\frac{4\pi^2}{n_\theta^2}\sum_{i,j=1}^{n_\theta} \big(2w_{k(\cdot,\rho_n)}(\tfrac{2\pi}{n_\theta})\|\Lambda^{\omega_{\tilde{i}}}\|_{L^{\infty}}+w_{\Lambda^{\omega_{\tilde{i}}}}(\tfrac{2\pi}{n_\theta})\big)\\
\leq &8\pi^3 c_{\omega_{\tilde{i}}}^2 \big((2\omega_2+2\omega_{\tilde{i}}) \|\eta(y)\|_{L^1(\mathcal{B}_R)}+(2\omega_2+2+2\omega_{\tilde{i}})\|\gamma(y)\|_{L^1(\mathcal{B}_R)})\big)n_\theta^{-1}.
\end{align*}
Finally, taking all above estimates and the assumption $\omega_1\leq\omega_2$ into account, we derive from \eqref{eqn:bf_bd0} that
\begin{align*}
|\Phi-\phi^0(\Lambda)|_{\infty}\leq\max_{m,n,\tilde{j}}\sum_{\tilde{i}=1}^2 {\rm I}_{\tilde{i}\tilde{j}}^0\leq c_{{\bf bf}_0} n_\theta^{-1}
\end{align*}
with the constant $c_{{\bf bf}_0}=8\pi^3 (4\omega_2 +3)(c_{\omega_1}^2+c_{\omega_2}^2)(\|\eta(y)\|_{L^1(\mathcal{B}_R)}+ \|\gamma(y)\|_{L^1(\mathcal{B}_R)})$.
This completes the proof.
\end{proof}

\begin{remark}
The proposed uncompressed $0$-level BFNN $\phi^0$ is a neural network of 4 channels with 3 or 4 weight layers per channel on $2n_\theta$ shifted inputs followed by explicit matrix transformations and {an added layer}. Most weight layers (14 in total) share the same weights. In fact, there are 5 different weights including $C\in\mathbb{R}^{n_\theta\times n_\theta}$, $K$, $\bar{K}\in\mathbb{C}^{n_\theta\times n_\rho}$ (or $K^*\in\mathbb{C}^{n_\rho\times n_\theta}$) and $\frac{4\pi^2 c_{\omega_1}}{n_\theta^2}b^T,\frac{4\pi^2 c_{\omega_2}}{n_\theta^2}b^T\in\mathbb{C}^{n_\theta}$ with $n_\theta$ and $2$ free parameters in $C$ and $\frac{4\pi^2 c_{\omega_1}}{n_\theta^2}b^T\&\frac{4\pi^2 c_{\omega_2}}{n_\theta^2}b^T$ respectively. Furthermore, by noting the {symmetry} of the kernel $k(t,\rho)=e^{-i\omega_1 \rho \cos t}$ along the $\rho$-axis, there are $\frac12 n_\theta n_\rho$ free parameters in {$K$, $\bar{K}$ and $K^*$}. 
\end{remark}

By noting the complementary low-rank property of the integral kernel $k$ stated in Theorem \ref{thm:low_rank}, we can further approximate $k$ by a low-rank operator and simulate this process by compressing the matrix $K$ via the butterfly factorization to reduce the computational complexity of the neural networks. 
We give the approximation analysis of $K$ by an $n_r \times n_r$ low-rank (of rank $r$) block matrix $K_r$ in the following lemma.

\begin{lemma}\label{lem:K_r}
Let $K$ be the discretization of the integral kernel $k(t,\rho)=e^{-i\omega_1 \rho \cos t}$ on the gird $\{(t_i,\rho_j)=(\frac{2\pi i}{n_\theta},\frac{\omega_2}{\omega_1}\frac{j}{n_\rho})\}_{i,j=1}^{i=n_\theta,j=n_\rho}${, cf. \eqref{eqn:K}}. Then for any $r,n_r\in\mathbb{N}$ such that $\frac{n_\theta}{n_r},\frac{n_\rho}{n_r}\in \mathbb{N}$, there exists an $n_r \times n_r$ low-rank (of rank $r$) block $K_r$ such that 
$|K-K_r|_\infty\leq(\frac{\pi}{2} \omega_2)^r(r!)^{-1}n_r^{-2r}$.
\end{lemma}

\begin{proof}
This lemma can be viewed as a corollary of Theorem \ref{thm:low_rank}.
After dividing the domain $[0,2\pi]\times[0,\frac{\omega_2}{\omega_1}]$ into $n_r^2$ subdomains of same size, i.e., $\{{\bf D}_{pq}=[\frac{2\pi (p-1)}{n_r},\frac{2\pi p}{n_r}]\times[\frac{\omega_2}{\omega_1}\frac{q-1}{n_r},\frac{\omega_2}{\omega_1}\frac{q}{n_r}]\}_{p,q=1}^{n_r}$ with centers $\{(\tau_p,\varrho_q)=(\frac{\pi (2p-1)}{n_r},\frac{\omega_2}{2\omega_1}\frac{2q-1}{n_r})\}_{p,q=1}^{n_r}$,
we claim that there exists a rank-$r$ approximation to the discretization of $k$ over each subdomain. 
In fact, we can define a piecewise continuous function $k_r$ by \eqref{eqn:low_rank_func} with $(t_0,\rho_0)=(\tau_p,\varrho_q)$ over the subdomain ${\bf D}_{pq}$ for any $p,q=1,\cdots n_r$, i.e., 
\begin{align*}
k_{r}(t,\rho)=\sum_{l=0}^{r-1}\frac{(-i\omega_1)^l}{l !}(\cos t-\cos \tau_p)^l e^{-i\omega_1\varrho_q(\cos t-\cos \tau_p)}(\rho-\varrho_q)^l e^{-i\omega_1\rho \cos \tau_p}.
\end{align*}
Then, by \eqref{eqn:low_rank_c} with $c_{r,\epsilon}:=\sup_{p,q}\sup_{(t,\rho)\in{\bf D}_{pq}}\omega_1|t-\tau_p||\rho-\varrho_q|\leq \omega_1 \frac{\pi}{n_r}\frac{\omega_2}{2\omega_1}\frac{1}{n_r}= \frac{\pi}{2} \omega_2 n_r^{-2}$, we have
\begin{align*}
|k(t,\rho)-k_{r}(t,\rho)|\leq\frac{c_{r,\epsilon}^r}{r!}\leq\frac{(\frac{\pi}{2} \omega_2)^r}{r!}n_r^{-2r},
\quad \forall  (t,\rho)\in[0,2\pi]\times[0,\frac{\omega_2}{\omega_1}].
\end{align*}
Finally, let $K_r$ be the discretization of $k_r(t,\rho)$ on the gird $\{(t_i,\rho_j)=(\frac{2\pi i}{n_\theta},\frac{\omega_2}{\omega_1}\frac{j}{n_\rho})\}_{i,j=1}^{i=n_\theta,j=n_\rho}$, there holds
\begin{align*}
|K-K_r|_\infty=\sup_{i,j}|k(t_i,\rho_j)-k_r(t_i,\rho_j)|\leq\frac{(\frac{\pi}{2} \omega_2)^r}{r!}n_r^{-2r}.
\end{align*}
Further, the complementary low-rank feature of $k_r$ gives the complementary low-rank property of $K_r$. In particular, if we denote the discretization of $k_r(t,\rho)$ over the subdomain ${\bf D}_{pq}$ by $K_{r,pq}$, then any submatrix of $K_{r,pq}$ has rank of at most $r$.
\end{proof}

Now, we construct a compressed $1$-level BFNN (corresponding to the middle level butterfly factorization), denoted by $\phi^1$, where the exact discretization matrix $K$ in $\phi^0$ is replaced by the $n_r\times n_r$ low-rank block matrix $K_r$ described in Lemma \ref{lem:K_r}.

\begin{theorem}\label{thm:BF_1_app}
Let $F^*$ be the adjoint operator of the linearized forward operator $F$ defined in \eqref{eqn:prob_linear_sys}. Then for any $r,n_r\in\mathbb{N}$ satisfying $\frac{n_\theta}{n_r},\frac{n_\rho}{n_r}\in \mathbb{N}$ and $n_r^{2r}=c_r n_\theta$ for some constant $c_r>0$,
there exists a neural network $\phi^1$ on the input $\Lambda$, which incorporates the complementary low-rank approximation into $\phi^0$ proposed in Theorem \ref{thm:BF_0_app} by replacing $K$ {defined in \eqref{eqn:K}} with $K_r$ defined in Lemma \ref{lem:K_r}, such that
\begin{align*}
|\Phi-\phi^1(\Lambda)|_{\infty}\leq c_{{\bf bf}_1} n_\theta^{-1}
\end{align*}
with the constant $c_{{\bf bf}_1}=8\pi^2(c_{\omega_1}^2+c_{\omega_2}^2)\big( (4\omega_2 +3)\pi+ \frac{(\frac{\pi}{2} \omega_2)^r}{r!c_r}\big)(\|\eta(y)\|_{L^1(\mathcal{B}_R)}+\|\gamma(y)\|_{L^1(\mathcal{B}_R)})$.
\end{theorem}
\begin{proof}
We follow the analysis in the proof of Theorem \ref{thm:BF_0_app}. By replacing the operator $\phi^0$ with $\phi^1$ in \eqref{eqn:bf_bd0}, there holds
\begin{align}\label{eqn:bf_bd1}
|\Phi-\phi^1(\Lambda)|_{\infty}\leq\max_{m,n,\tilde{j}}\sum_{\tilde{i}=1}^2 {\rm I}_{\tilde{i}\tilde{j}}^1
\end{align}
where
\begin{align*}
{\rm I}_{\tilde{i}1}^1=& c_{\omega_{\tilde{i}}}\big|-\int_0^{2\pi} \overline{k(\theta_{\hat{x}},\rho_n)}\int_0^{2\pi} k(\theta_{\zeta},\rho_n)\cos(\theta_{\hat{x}}-\theta_{\zeta})  \Lambda_{\theta_m}^{\omega_{\tilde{i}}}(\theta_{\hat{x}},\theta_{\zeta}){\rm d}\theta_{\zeta}{\rm d}\theta_{\hat{x}}\\
&\quad\;\;-\frac{4\pi^2}{n_\theta^2}\sum_{i=1}^{n_\theta} \overline{K_r[i,n]}\sum_{j=1}^{n_\theta} K_r[j,n] C[i,j]\Lambda_{\theta_m}^{\omega_{\tilde{i}}}[i,j]\big|\\
\leq &{\rm I}_{\tilde{i}1}^0+c_{\omega_{\tilde{i}}}\frac{4\pi^2}{n_\theta^2}\big|\sum_{i=1}^{n_\theta} \overline{K[i,n]}\sum_{j=1}^{n_\theta} K[j,n] C[i,j]\Lambda_{\theta_m}^{\omega_{\tilde{i}}}[i,j]-\sum_{i=1}^{n_\theta} \overline{K_r[i,n]}\sum_{j=1}^{n_\theta} K_r[j,n] C[i,j]\Lambda_{\theta_m}^{\omega_{\tilde{i}}}[i,j]\big|,\\
{\rm I}_{\tilde{i}2}^1=& c_{\omega_{\tilde{i}}}\big|\int_0^{2\pi} \overline{k(\theta_{\hat{x}},\rho_n)}\int_0^{2\pi} k(\theta_{\zeta},\rho_n) \Lambda_{\theta_m}^{\omega_{\tilde{i}}}(\theta_{\hat{x}},\theta_{\zeta}){\rm d}\theta_{\zeta}{\rm d}\theta_{\hat{x}}-\frac{4\pi^2}{n_\theta^2}\sum_{i=1}^{n_\theta} \overline{K_r[i,n]}\sum_{j=1}^{n_\theta} K_r[j,n]\Lambda_{\theta_m}^{\omega_{\tilde{i}}}[i,j]\big|\\
\leq &{\rm I}_{\tilde{i}2}^0+c_{\omega_{\tilde{i}}}\frac{4\pi^2}{n_\theta^2}\big|\sum_{i=1}^{n_\theta} \overline{K[i,n]}\sum_{j=1}^{n_\theta} K[j,n] \Lambda_{\theta_m}^{\omega_{\tilde{i}}}[i,j]-\sum_{i=1}^{n_\theta} \overline{K_r[i,n]}\sum_{j=1}^{n_\theta} K_r[j,n] \Lambda_{\theta_m}^{\omega_{\tilde{i}}}[i,j]\big|
\end{align*}
with ${\rm I}_{\tilde{i}1}^0$ and ${\rm I}_{\tilde{i}2}^0$ being defined in the proof of Theorem \ref{thm:BF_0_app}.
We decompose the first error for the low-rank approximation stated above into
\begin{align*}
{\rm I_{e}}:=
&\big|\sum_{i=1}^{n_\theta} \overline{K[i,n]}\sum_{j=1}^{n_\theta} K[j,n] C[i,j]\Lambda_{\theta_m}^{\omega_{\tilde{i}}}[i,j]-\sum_{i=1}^{n_\theta} \overline{K_r[i,n]}\sum_{j=1}^{n_\theta} K_r[j,n]C[i,j] \Lambda_{\theta_m}^{\omega_{\tilde{i}}}[i,j]\big|\\
%\leq&\big|\sum_{i=1}^{n_\theta} \overline{K[i,n]}\sum_{j=1}^{n_\theta} K[j,n] C[i,j]\Lambda_{\theta_m}^{\omega_{\tilde{i}}}[i,j]-\sum_{i=1}^{n_\theta} \overline{K_r[i,n]}\sum_{j=1}^{n_\theta} K[j,n] C[i,j]\Lambda_{\theta_m}^{\omega_{\tilde{i}}}[i,j]\big|\\
%&+\big|\sum_{i=1}^{n_\theta} \overline{K_r[i,n]}\sum_{j=1}^{n_\theta} K[j,n] C[i,j]\Lambda_{\theta_m}^{\omega_{\tilde{i}}}[i,j]-\sum_{i=1}^{n_\theta} \overline{K_r[i,n]}\sum_{j=1}^{n_\theta} K_r[j,n] C[i,j]\Lambda_{\theta_m}^{\omega_{\tilde{i}}}[i,j]\big|\\
\leq&\big|\sum_{i=1}^{n_\theta} (\overline{K[i,n]}-\overline{K_r[i,n]})\sum_{j=1}^{n_\theta} K[j,n] C[i,j]\Lambda_{\theta_m}^{\omega_{\tilde{i}}}[i,j]\big|+\big|\sum_{i=1}^{n_\theta} \overline{K_r[i,n]}\sum_{j=1}^{n_\theta} (K[j,n]-K_r[j,n]) C[i,j]\Lambda_{\theta_m}^{\omega_{\tilde{i}}}[i,j]\big|\\
\leq&\sum_{i=1}^{n_\theta} \big|\overline{K[i,n]}-\overline{K_r[i,n]}\big|\sum_{j=1}^{n_\theta} \|k\|_{L^{\infty}}\|\Lambda^{\omega_{\tilde{i}}}\|_{L^{\infty}}+\sum_{i=1}^{n_\theta} \|\bar{k}\|_{L^{\infty}}\sum_{j=1}^{n_\theta} \big|K[j,n]-K_r[j,n]\big|\|\Lambda^{\omega_{\tilde{i}}}\|_{L^{\infty}}.
\end{align*}
Then, Lemma \ref{lem:est1} with $R=1$ and Lemma \ref{lem:K_r} give
\begin{align*}
{\rm I_{e}}
\leq& 2\sum_{i=1}^{n_\theta} \sum_{j=1}^{n_\theta} \frac{(\frac{\pi}{2} \omega_2)^r}{r!}n_r^{-2r}\|\Lambda^{\omega_{\tilde{i}}}\|_{L^{\infty}}
\leq 2 c_{\omega_{\tilde{i}}}\frac{(\frac{\pi}{2} \omega_2)^r}{r!}(\|\eta(y)\|_{L^1(\mathcal{B}_R)}+\|\gamma(y)\|_{L^1(\mathcal{B}_R)})n_\theta^2 n_r^{-2r},
\end{align*}
which implies
\begin{align*}
{\rm I}_{\tilde{i}1}^1\leq &{\rm I}_{\tilde{i}1}^0+c_{\omega_{\tilde{i}}}\frac{4\pi^2}{n_\theta^2}{\rm I_e}\leq {\rm I}_{\tilde{i}1}^0+8\pi^2 c_{\omega_{\tilde{i}}}^2 \frac{(\frac{\pi}{2} \omega_2)^r}{r!}(\|\eta(y)\|_{L^1(\mathcal{B}_R)}+\|\gamma(y)\|_{L^1(\mathcal{B}_R)}) n_r^{-2r}.
\end{align*}
Similarly, we can bound the second term ${\rm I}_{\tilde{i}2}^1$ by
\begin{align*}
{\rm I}_{\tilde{i}2}^1\leq&{\rm I}_{\tilde{i}2}^0+8\pi^2 c_{\omega_{\tilde{i}}}^2 \frac{(\frac{\pi}{2} \omega_2)^r}{r!}(\|\eta(y)\|_{L^1(\mathcal{B}_R)}+\|\gamma(y)\|_{L^1(\mathcal{B}_R)}) n_r^{-2r}.
\end{align*}
By the estimates above and in the proof of Theorem \ref{thm:BF_0_app}, we derive from \eqref{eqn:bf_bd1} that
\begin{align*}
|\Phi-\phi^1(\Lambda)|_{\infty}\leq&\max_{m,n,\tilde{j}}\sum_{\tilde{i}=1}^2 ({\rm I}_{\tilde{i}\tilde{j}}^0+8\pi^2 c_{\omega_{\tilde{i}}}^2 \frac{(\frac{\pi}{2} \omega_2)^r}{r!}(\|\eta(y)\|_{L^1(\mathcal{B}_R)}+\|\gamma(y)\|_{L^1(\mathcal{B}_R)}) n_r^{-2r})\\
\leq& c_{{\bf bf}_0} n_\theta^{-1}+8\pi^2 (c_{\omega_1}^2+c_{\omega_2}^2) \frac{(\frac{\pi}{2} \omega_2)^r}{r!}(\|\eta(y)\|_{L^1(\mathcal{B}_R)}+\|\gamma(y)\|_{L^1(\mathcal{B}_R)}) n_r^{-2r}.
\end{align*}
Finally, setting $n_r^{2r}=c_r n_\theta$ for some constant $c_r>0$ indicates
\begin{align*}
|\Phi-\phi^1(\Lambda)|_{\infty}
\leq& c_{{\bf bf}_1} n_\theta^{-1}.
\end{align*}
with the constant $c_{{\bf bf}_1}=8\pi^2(c_{\omega_1}^2+c_{\omega_2}^2)\big( (4\omega_2 +3)\pi+ \frac{(\frac{\pi}{2} \omega_2)^r}{r!c_r}\big)(\|\eta(y)\|_{L^1(\mathcal{B}_R)}+\|\gamma(y)\|_{L^1(\mathcal{B}_R)})$.
This completes the proof.
\end{proof}

\begin{remark}
In Theorem \ref{thm:BF_1_app}, we did not use the exact butterfly structure of $K_r$ to analyze the approximation capability of the compressed $1$-level BFNN $\phi^1$. In fact, in middle level butterfly factorization, the low-rank block matrix $K_r\in\mathbb{C}^{n_\theta\times n_\rho}$ is decomposed into $K_r=U M V$ with $M\in \mathbb{C}^{n_r^2 r\times n_r^2 r}$ being a weighted permutation matrix generated by the singular values of the low-rank block matrices (see Section \ref{sec:hypothesis} for details) and $U\in\mathbb{C}^{n_\theta \times n_r^2 r}$ and $V\in\mathbb{C}^{n_r^2 r \times n_\rho}$ being diagonal block matrices generated by the left and right singular vectors, respectively, see \cite{LiYangMartin:2015} for the structure of the decomposition. This decomposition reduces the computational complexity but theoretically does not affect the approximation. In particular, there are $n_r^2r$, $n_\theta n_r r$ and $n_\rho n_r r$ nonzero entries in $M$, $U$ and $V$ respectively and half amount of free parameters due to the symmetry of $K_r$.
Further, with the condition $n_r^{2r}=c_r n_\theta$ in Theorem \ref{thm:BF_1_app}, the number of {nonzero parameters} in the decomposition $K_r=U M V$ is given by
\begin{align*}
|M|+|U|+|V|=(n_\theta+n_\rho+n_r)n_r r=(n_\theta+n_\rho+n_r)\ln (c_r n_\theta) \frac{n_r}{2\ln n_r}.
\end{align*}
\end{remark}

\subsection{Convolutional neural networks (CNNs) for $(F^*F+\alpha I)^{-1}$}\label{sec:app_CNN}
We shall construct a CNN with scaled residual layers, denoted by $\psi$, to approximate the convolution-type inverse operator $(F^*F+\alpha I)^{-1}: \left(L^2(\mathcal{B}_R)\right)^2\to \left(L^2(\mathcal{B}_R)\right)^2$ ,
after converting the polar coordinates to the Cartesian coordinates by a operator $\psi_0$ which transforms the output of the BFNN $\phi^0$ or $\phi^1$ in $\mathbb{R}^{2 n_\theta\times \frac{\omega_1}{\omega_2}n_\rho}$ (on the domain $[0,2\pi]\times[0,1]$) into the input of the CNN $\psi$ in $\mathbb{R}^{2 n_c\times n_c}$, denoted by $O_I$ (on the domain $[-1,1]^2$). 
For the sake of simplifying the analysis, we define $n_c=n_\theta$.
We discretize the domain $[-1,1]^2$ of the input, intermediate output and final output uniformly with $\{x_{ij}=(\frac{2i-1}{n_c}-1,\frac{2j-1}{n_c}-1)\}_{i,j=1}^{n_c}$ (symmetric about the origin) and assume that the input is extended to be a piecewise constant function (corresponding to the proposed grid) over the domain $[-1,1]^2$ when employing the exact (continuous) convolution kernel.

By noting the structure of $(F^*F+\alpha I)^{-1}$ in Proposition \ref{prop:F*F_formula}, i.e., 
\begin{align}\label{eqn:F*F}
(F^*F+\alpha I)^{-1}(O_{I})=&(g_2 I_2)*\sum_{i=1}^2\begin{bmatrix}
g_{22}^{\omega_i}+\frac\alpha2\delta &-g_{12}^{\omega_i} \\
-g_{21}^{\omega_i} &g_{11}^{\omega_i}+\frac\alpha2\delta
\end{bmatrix}*O_{I}
:=(g_2I_2)*g_1*O_{I},
\end{align}
with $g_{21}^\omega=g_{12}^\omega$ and $g_2=\mathcal{F}^{-1}\big(\det(\mathcal{F}({\bf g}))^{-1}\big)=\mathcal{F}^{-1}\big(\mathcal{F}(g_{de})^{-1}\big)$ being the kernel for the deconvolution of the convolution with kernel $g_{de}=\sum_{i,i'=1}^2(g_{11}^{\omega_i}+\frac\alpha2\delta)*(g_{22}^{\omega_{i'}}+\frac\alpha2\delta) -g_{12}^{\omega_i}*g_{12}^{\omega_{i'}}$, 
this operator $(F^*F+\alpha I)^{-1}$ can be approximated by 
\begin{enumerate}
\item[(i)] a single layer CNN, with $3\times 2$ filters $\{g_{11}^{\omega_i},\; g_{22}^{\omega_i},\;-g_{12}^{\omega_i}\}_{i=1}^2$, combined with explicit matrix transformations and a scaled residual layer for simulating $g_1*$, denoted by $\psi_1$,
\item[(ii)] followed by a CNN with a scaled residual layer for simulating $g_2 I_2*$, denoted by $\psi_2$. 
\end{enumerate}
Further, we denote the discretizations of $[g_1*O_I](y)$ and $[(g_2I_2)*g_1*O_I](y)$ on the domain $[-1,1]^2$ by $\Psi_1$ and $\Psi_2$ respectively.

\begin{theorem}\label{thm:CNN_app}
Let $(F^*F+\alpha I)^{-1}$ be the inverse pseudo-differential
operator of the linearized forward operator $F$ defined in \eqref{eqn:prob_linear_sys}. Then there exists a convolutional neural network with scaled residual layers $\psi=\psi_2\circ \psi_1$ on the input $O_I%=\begin{bmatrix}
%O_{I1}\\
%O_{I2}
%\end{bmatrix}
\in \mathbb{R}^{2 n_c\times n_c}$ given by
\begin{align}
\psi_1(O_I)=&\sum_{i=1}^2 (P_u^TG_{22}^{\omega_i}*P_u +P_l^TG_{11}^{\omega_i}*P_l-P_u^TG_{12}^{\omega_i}*P_l - P_l^TG_{12}^{\omega_i}*P_u)O_I+\alpha O_I:=\tilde{O}_I,\label{eqn:Psi_1}\\
\psi_2(\tilde{O}_I)=&\alpha^{-2}\tilde{O}_I -(P_u^T G_{2,\alpha}*P_u +P_l^T G_{2,\alpha}*P_l )\tilde{O}_I\label{eqn:Psi_2}
\end{align}
with explicit matrix transformations $P_u$ and $P_l$ defined in \eqref{eqn:P_ul} and convolution kernels $G_{jj'}^{\omega_i}, G_{2,\alpha}\in \mathbb{R}^{2n_c\times 2n_c}$ for any $j,j',i=1,2$, such that
\begin{align*}
|\Psi_2-\psi_2\circ \psi_1(O_I)|_{\infty}\leq c_{{\bf cnn}}|O_{I}|_{\infty}\alpha^{-2}n_c^{-1}
\end{align*}
with the constant $c_{{\bf cnn}}=2^7\sqrt{2} \pi\omega_{2}( c_{\omega_{1}}^2+ c_{\omega_{2}}^2)(\pi+2^3\omega_2^2)$.

\end{theorem}
\begin{proof}
\noindent$\textbf{(i)}$
We start with $\psi_1$ defined in \eqref{eqn:Psi_1} for approximating $g_1$ on the input $O_I=(O_{I1},O_{I2})^T$, i.e.,
\begin{align*}
\psi_1(O_I)
=&\sum_{i=1}^2 P_u^TG_{22}^{\omega_i}*P_uO_{I}+P_l^TG_{11}^{\omega_i}*P_lO_{I}-(P_u^TG_{12}^{\omega_i}*P_lO_{I}+ P_l^TG_{12}^{\omega_i}*P_uO_{I})+\alpha O_I\\
=&\sum_{i=1}^2\begin{bmatrix}
G_{22}^{\omega_i}*O_{I1} -G_{12}^{\omega_i}*O_{I2}\\
-G_{12}^{\omega_i}*O_{I1}+G_{11}^{\omega_i}*O_{I2}
\end{bmatrix}+\alpha O_I,
\end{align*}
where the filter $G_{jj'}^{\omega_i}\in\mathbb{R}^{2n_c\times 2n_c}$ is the uniform discretization of $g_{jj'}^{\omega_i}$ over $[-2,2]^2$ (for any $j,j',i=1,2$) and $P_u$, $P_l$ are the explicit matrix transformations defined in \eqref{eqn:P_ul}. 
By the definition of $g_1$ in \eqref{eqn:F*F}, we can bound the element-wise difference between the discretization of $g_1*O_I$ (i.e., $\Psi_1$) and the output $\psi_1(O_I)$ by
\begin{align*}
|\Psi_1-\psi_1(O_I)|_{\infty}
%=&\max_{m,n}\bigg|\sum_{i=1}^2\begin{bmatrix}
%(g_{22}^{\omega_i} *O_{I1} -g_{12}^{\omega_i}*O_{I2} \\
%-g_{12}^{\omega_i}*O_{I1}+(g_{11}^{\omega_i} *O_{I2}
%\end{bmatrix}(x_{mn})-\sum_{i=1}^2\begin{bmatrix}
%G_{22}^{\omega_i}*O_{I1} -%G_{12}^{\omega_i}*O_{I2} \\
%-G_{12}^{\omega_i}*O_{I1}+G_{11}^{\omega_i}*O_{I2}
%\end{bmatrix}[m,n]\bigg|_{\infty}\\
=&\max_{m,n}\bigg|\sum_{i=1}^2\bigg(\begin{bmatrix}
g_{22}^{\omega_i} *O_{I1} -g_{12}^{\omega_i}*O_{I2} \\
-g_{12}^{\omega_i}*O_{I1}+g_{11}^{\omega_i} *O_{I2}
\end{bmatrix}(x_{mn})-\begin{bmatrix}
G_{22}^{\omega_i}*O_{I1} -G_{12}^{\omega_i}*O_{I2} \\
-G_{12}^{\omega_i}*O_{I1}+G_{11}^{\omega_i}*O_{I2}
\end{bmatrix}[m,n]\bigg)\bigg|_{\infty}\\
%\leq &\max_{m,n}\sum_{i=1}^2\max\Big(\big|\big( [g_{22}^{\omega_i} *O_{I1}](x_{mn}) -(G_{22}^{\omega_i}*O_{I1})[m,n]\big)-\big([g_{12}^{\omega_i}*O_{I2}](x_{mn})-(G_{12}^{\omega_i}*O_{I2})[m,n]\big)\big|,\\
%&\qquad\qquad\;\;\big|\big( [g_{11}^{\omega_i} *O_{I2}](x_{mn})-(G_{11}^{\omega_i}*O_{I2})[m,n]\big)-\big([g_{12}^{\omega_i}*O_{I1}](x_{mn})-(G_{12}^{\omega_i}*O_{I1})[m,n]\big)\big|\Big)\\
\leq&\max_{m,n}\sum_{\tilde{i}=1}^2\max({\rm I}^2_{\tilde{i}1},{\rm I}^2_{\tilde{i}2})=\max_{m,n,\tilde{j}}\sum_{\tilde{i}=1}^2 {\rm I}^2_{\tilde{i}\tilde{j}},
\end{align*}
where
\begin{align*}
{\rm I}^2_{\tilde{i}1}=&\big|\big( [g_{22}^{\omega_{\tilde{i}}} *O_{I1}](x_{mn}) -(G_{22}^{\omega_{\tilde{i}}}*O_{I1})[m,n]\big)-\big([g_{12}^{\omega_{\tilde{i}}}*O_{I2}](x_{mn})-(G_{12}^{\omega_{\tilde{i}}}*O_{I2})[m,n]\big)\big|,\\
{\rm I}^2_{\tilde{i}2}=&\big|\big( [g_{11}^{\omega_{\tilde{i}}} *O_{I2}](x_{mn})-(G_{11}^{\omega_{\tilde{i}}}*O_{I2})[m,n]\big)-\big([g_{12}^{\omega_{\tilde{i}}}*O_{I1}](x_{mn})-(G_{12}^{\omega_{\tilde{i}}}*O_{I1})[m,n]\big)
\big|.
\end{align*}
We follow the techniques used in the proof of Theorem \ref{thm:BF_0_app}. 
For the first term ${\rm I}^2_{\tilde{i}1}$,
we have
\begin{align*}
{\rm I}^2_{\tilde{i}1}
\leq&\big|\int_{[-1,1]^2}g_{22}^{\omega_{\tilde{i}}}(x_{mn}-y)O_{I1}(y){\rm d}y -\frac{4}{n_c^2}\sum_{i,j=1}^{n_c}G_{22}^{\omega_{\tilde{i}}}[m-i,n-j]O_{I1}[i,j]\big|\\
&+\big|\int_{[-1,1]^2}g_{12}^{\omega_{\tilde{i}}}(x_{mn}-y)O_{I2}(y){\rm d}y -\frac{4}{n_c^2}\sum_{i,j=1}^{n_c}G_{12}^{\omega_{\tilde{i}}}[m-i,n-j]O_{I2}[i,j]\big|\\
\leq&\sum_{i,j=1}^{n_c}\big|\int_{\frac{2(i-1)}{n_c}-1}^{\frac{2i}{n_c}-1}\int_{\frac{2(j-1)}{n_c}-1}^{\frac{2j}{n_c}-1}g_{22}^{\omega_{\tilde{i}}}(x_{mn}-y)O_{I1}[i,j]{\rm d}y -\frac{4}{n_c^2}G_{22}^{\omega_{\tilde{i}}}[m-i,n-j]O_{I1}[i,j]\big|\\
&+\sum_{i,j=1}^{n_c}\big|\int_{\frac{2(i-1)}{n_c}-1}^{\frac{2i}{n_c}-1}\int_{\frac{2(j-1)}{n_c}-1}^{\frac{2j}{n_c}-1}
g_{12}^{\omega_{\tilde{i}}}(x_{mn}-y)O_{I2}[i,j]{\rm d}y -\frac{4}{n_c^2}G_{12}^{\omega_{\tilde{i}}}[m-i,n-j]O_{I2}[i,j]\big|\\
\leq&|O_{I1}|_{\infty}\sum_{i,j=1}^{n_c}\int_{\frac{2(i-1)}{n_c}-1}^{\frac{2i}{n_c}-1}\int_{\frac{2(j-1)}{n_c}-1}^{\frac{2j}{n_c}-1}\big|(g_{22}^{\omega_{\tilde{i}}}(x_{mn}-y) -g_{22}^{\omega_{\tilde{i}}}(x_{mn}-x_{ij}))\big|{\rm d}y\\
&+|O_{I2}|_{\infty}\sum_{i,j=1}^{n_c}\int_{\frac{2(i-1)}{n_c}-1}^{\frac{2i}{n_c}-1}\int_{\frac{2(j-1)}{n_c}-1}^{\frac{2j}{n_c}-1}\big|g_{12}^{\omega_{\tilde{i}}}(x_{mn}-y)-g_{12}^{\omega_{\tilde{i}}}(x_{mn}-x_{ij})\big|{\rm d}y \\
\leq&\frac{4}{n_c^2}\sum_{i,j=1}^{n_c}w_{g_{22}^{\omega_{\tilde{i}}}}(\tfrac{2}{n_c})|O_{I1}|_{\infty}+\frac{4}{n_c^2}\sum_{i,j=1}^{n_c}w_{g_{12}^{\omega_{\tilde{i}}}}(\tfrac{2}{n_c}) |O_{I2}|_{\infty}=4\big(w_{g_{22}^{\omega_{\tilde{i}}}}(\tfrac{2}{n_c})+w_{g_{12}^{\omega_{\tilde{i}}}}(\tfrac{2}{n_c})\big) |O_{I}|_{\infty}.
\end{align*}
Further, Lemma \ref{lem:filter} indicates that
\begin{align*}
{\rm I}^2_{\tilde{i}1}
\leq  4\big(2^4\sqrt{2} \pi^2\omega_{\tilde{i}} c_{\omega_{\tilde{i}}}^2 n_c^{-1}
+2^4\sqrt{2} \pi^2\omega_{\tilde{i}} c_{\omega_{\tilde{i}}}^2 n_c^{-1}\big)|O_{I}|_{\infty}
\leq 2^7\sqrt{2} \pi^2\omega_{\tilde{i}} c_{\omega_{\tilde{i}}}^2
|O_{I}|_{\infty}n_c^{-1}.
\end{align*}
Similarly, for the second term ${\rm I}^2_{\tilde{i}2}$, there holds
\begin{align*}
{\rm I}^2_{\tilde{i}2}
\leq 2^7\sqrt{2} \pi^2\omega_{\tilde{i}} c_{\omega_{\tilde{i}}}^2
|O_{I}|_{\infty}n_c^{-1}.
\end{align*}
Thus, there holds
\begin{align}\label{eqn:app_g_1}
&|\Psi_1-\psi_1(O_I)|_{\infty}
\leq\max_{m,n,\tilde{j}}\sum_{\tilde{i}=1}^2 {\rm I}^2_{\tilde{i}\tilde{j}}\leq 2^7\sqrt{2} \pi^2(\omega_{1} c_{\omega_{1}}^2+\omega_{2} c_{\omega_{2}}^2)
|O_{I}|_{\infty}n_c^{-1}.
\end{align}
%It is worth noting that we can represent any single convolution operation as a matrix-vector multiplication by flattening the input (and output) matrices $O_{I1},O_{I2}\in\mathbb{R}^{n_c\times n_c}$ along the column into vectors $o_{I1},o_{I2}\in\mathbb{R}^{n_c^2}$ and generating an $n_c^2\times n_c^2$ doubly block Toeplitz matrix (each block is itself an $n_c\times n_c$ Toeplitz matrix) from the filter where each filter coefficient appears once per column and row (see the circulant case in \cite[Section 3.2]{GoukFrankaPfahringer:2021}). 
%We denote the doubly block Toeplitz matrix generated from $G_{jj'}^{\omega_i}$ by $G_{o,jj'}^{\omega_i}$, then
%$|G_{o,jj'}^{\omega_i}|_\infty=|G_{jj'}^{\omega_i}|_\infty$.
%Furthermore, in view of the radially symmetry of $g_{jj'}^{\omega_i}$ (Lemma \ref{lem:filter}), the entire doubly block Toeplitz matrix $G_{o,jj'}^{\omega_i}$ and all blocks are symmetric.

\vspace{5mm}
\noindent$\textbf{(ii)}$
Now, we discuss the second part of the neural network, i.e., $\psi_2$ which simulates the operator $g_2I_2$. 
In light of Lemma \ref{lem:det}, we decompose $g_2$ as $g_2=\alpha^{-2}\delta-g_{2,\alpha}$ where
\begin{align}\label{eqn:g_2}
w_{g_{2,\alpha}}(s)\leq&\frac{4\sqrt{2}\omega_2^3}{3\pi\alpha^{2}}s \quad \mbox{and}\quad \|g_{2,\alpha}^{\omega}\|_{L^{\infty}} \leq \frac{\omega_2^2}{\pi\alpha^{2}}.
\end{align}  
Similar to the analysis for $\psi_1$, we define $G_{2,\alpha}\in \mathbb{R}^{2n_c\times 2n_c}$ as the uniform discretization of $g_{2,\alpha}$ over $[-2,2]^2$ and $\psi_2$ as the operator given by \eqref{eqn:Psi_2}. %$\psi_2(\cdot)=\alpha^{-2}\cdot -(P_u^T G_{2,\alpha}*P_u \cdot+P_l^T G_{2,\alpha}*P_l \cdot)$. 
Then, we bound the element-wise difference between the discretization of $(g_2 I)*g_1*O_I$ (i.e., $\Psi_2$) and the output $\psi_2\circ \psi_1(O_I)$ by 
\begin{align}
&|\Psi_2-\psi_2\circ \psi_1(O_I)|_{\infty}\nonumber\\
=&\max_{m,n}\big|[(g_2 I)*(g_1*O_I-\psi_1(O_I))](x_{mn})+[(g_2 I)*\psi_1(O_I)](x_{mn})-\psi_2\circ \psi_1(O_I)[m,n]\big|_{\infty}\nonumber\\
=&\max_{m,n}\big|\alpha^{-2}\big(\Psi_1-\psi_1(O_I)\big)[m,n]\big|_{\infty}+\big|(g_{2,\alpha} I)*\big(g_1*O_I-\psi_1(O_I)\big)(x_{mn})\big|_{\infty}\nonumber\\
&\qquad\qquad\qquad\qquad\qquad\qquad\qquad\quad\;\;\;\;\;+\big|[(g_{2} I)*\psi_1(O_I)](x_{mn})-\psi_2\circ \psi_1(O_I)[m,n]\big|_{\infty}\nonumber\\
\leq& \alpha^{-2}|\Psi_1-\psi_1(O_I)|_{\infty}+\big(\int_{[-1,1]^2}\big|g_{2,\alpha}(x_{mn}-y)|{\rm d}y\big) |\Psi_1-\psi_1(O_I)|_{\infty}\nonumber\\
&\qquad\qquad\qquad\qquad\qquad\qquad\qquad\quad\;\;\;\;\;+\max_{m,n}\bigg|\begin{bmatrix}
g_{2,\alpha}*\tilde{O}_{I1} \\
g_{2,\alpha}*\tilde{O}_{I2}
\end{bmatrix}(x_{mn})
-\begin{bmatrix}
G_{2,\alpha}*\tilde{O}_{I1}\\
G_{2,\alpha}*\tilde{O}_{I2}
\end{bmatrix}[m,n]\bigg|_{\infty}\nonumber\\
\leq& (\alpha^{-2}+4\|g_{2,\alpha}\|_{L^\infty})|\Psi_1-\psi_1(O_I)|_{\infty}+\max_{m,n}\max({\rm I}^3_{1},{\rm I}^3_{2})\label{eqn:app_g_20}
\end{align}
where
\begin{align*}
{\rm I}^3_{1}=&\big|[g_{2,\alpha}*\tilde{O}_{I1}](x_{mn}) -\big(G_{2,\alpha}*\tilde{O}_{I1}\big)[m,n]\big| \quad\mbox{and}\quad
{\rm I}^3_{2}=\big|[g_{2,\alpha}*\tilde{O}_{I2}](x_{mn}) -\big(G_{2,\alpha}*\tilde{O}_{I2}\big)[m,n]\big|,
\end{align*}
with the notations $\tilde{O}_{I1}=\sum_{\tilde{i}=1}^2(G_{22}^{\omega_{\tilde{i}}}*O_{I1} -G_{12}^{\omega_{\tilde{i}}}*O_{I2})$ and $\tilde{O}_{I2}=\sum_{\tilde{i}=1}^2(-G_{12}^{\omega_{\tilde{i}}}*O_{I1}+G_{11}^{\omega_{\tilde{i}}}*O_{I2})$.
Similar to the analysis in part (i), we obtain that
\begin{align*}
{\rm I}^3_{1}%=&\big|[g_{2,\alpha}*\tilde{O}_{I1}](x_{mn}) -\big(G_{2,\alpha}*\tilde{O}_{I1}\big)[m,n]\big|\\
\leq&\big|\int_{[-1,1]^2}g_{2,\alpha}(x_{mn}-y)\tilde{O}_{I1}(y){\rm d}y -\frac{4}{n_c^2}\sum_{i,j=1}^{n_c}G_{2,\alpha}[m-i,n-j]\tilde{O}_{I1}[i,j]\big|\\
\leq&\sum_{i,j=1}^{n_c}\big|\int_{\frac{2(i-1)}{n_c}-1}^{\frac{2i}{n_c}-1}\int_{\frac{2(j-1)}{n_c}-1}^{\frac{2j}{n_c}-1}g_{2,\alpha}(x_{mn}-y)\tilde{O}_{I1}[i,j]{\rm d}y -\frac{4}{n_c^2}G_{2,\alpha}[m-i,n-j]\tilde{O}_{I1}[i,j]\big|\\
\leq&|\tilde{O}_{I1}|_{\infty}\sum_{i,j=1}^{n_c}\int_{\frac{2(i-1)}{n_c}-1}^{\frac{2i}{n_c}-1}\int_{\frac{2(j-1)}{n_c}-1}^{\frac{2j}{n_c}-1}\big|g_{2,\alpha}(x_{mn}-y) -g_{2,\alpha}(x_{mn}-x_{ij})\big|{\rm d}y\\
\leq&\frac{4}{n_c^2}\sum_{i,j=1}^{n_c}w_{g_{2,\alpha}}(\tfrac{2}{n_c})|\tilde{O}_{I1}|_{\infty}=4w_{g_{2,\alpha}}(\tfrac{2}{n_c})|\tilde{O}_{I1}|_{\infty}
\end{align*}
and ${\rm I}^3_{2}\leq 4w_{g_{2,\alpha}}(\tfrac{2}{n_c})|\tilde{O}_{I2}|_{\infty}$.
With the estimates in \eqref{eqn:app_g_1} and \eqref{eqn:g_2}, we derive from \eqref{eqn:app_g_20} that
\begin{align*}
&|\Psi_2-\psi_2\circ \psi_1(O_I)|_{\infty}
\leq (\alpha^{-2}+4\|g_{2,\alpha}\|_{L^\infty})|\Psi_1-\psi_1(O_I)|_{\infty}+4w_{g_{2}}(\tfrac{2}{n_c})\max(|\tilde{O}_{I1}|_{\infty},|\tilde{O}_{I2}|_{\infty})\\
\leq &2^7\sqrt{2} \pi^2(\omega_{1} c_{\omega_{1}}^2+\omega_{2} c_{\omega_{2}}^2)(1+4\pi^{-1}\omega_2^2)
|O_{I}|_{\infty}\alpha^{-2}n_c^{-1}+2^5\sqrt{2}\omega_2^3(3\pi)^{-1}\max(|\tilde{O}_{I1}|_{\infty},|\tilde{O}_{I2}|_{\infty})\alpha^{-2}n_c^{-1}.
\end{align*}
Finally, the estimate $\|g_{j{j'}}^{\omega_{\tilde{i}}}\|_{L^{\infty}} \leq 4\pi^2 c_\omega^2$ (for any $j,j',i=1,2$) in Lemma \ref{lem:filter} yields 
\begin{align*}
|\tilde{O}_{I1}|_{\infty}=&\sum_{\tilde{i}=1}^2\Big(\big|\frac{4}{n_c^2}\sum_{i,j=1}^{n_c}G_{22}^{\omega_{\tilde{i}}}[m-i,n-j]O_{I1}[i,j]\big|+\big|\frac{4}{n_c^2}\sum_{i,j=1}^{n_c}G_{12}^{\omega_{\tilde{i}}}[m-i,n-j]O_{I2}[i,j]\big|\Big)\\
\leq&\sum_{\tilde{i}=1}^2\Big(\frac{4}{n_c^2}\sum_{i,j=1}^{n_c}\|g_{22}^{\omega_{\tilde{i}}}\|_{L^{\infty}}+\frac{4}{n_c^2}\sum_{i,j=1}^{n_c}\|g_{12}^{\omega_{\tilde{i}}}\|_{L^{\infty}}\Big)|O_{I}|_{\infty}
\leq2^5\pi^2 (c_{\omega_1}^2+c_{\omega_2}^2)|O_{I}|_{\infty},
\end{align*}
and similarly $|\tilde{O}_{I2}|_{\infty}\leq 2^5\pi^2 (c_{\omega_1}^2+c_{\omega_2}^2)|O_{I}|_{\infty}$, which imply that
\begin{align*}
|\Psi_2-\psi_2\circ \psi_1(O_I)|_{\infty}
\leq &2^7\sqrt{2} \pi^2(\omega_{1} c_{\omega_{1}}^2+\omega_{2} c_{\omega_{2}}^2)(1+4\pi^{-1}\omega_2^2)
|O_{I}|_{\infty}\alpha^{-2}n_c^{-1}\\
&+2^{9}\sqrt{2}\pi\omega_2^3 (c_{\omega_1}^2+c_{\omega_2}^2)|O_{I}|_{\infty}\alpha^{-2}n_c^{-1}\\
\leq &2^7\sqrt{2} \pi\omega_{2}( c_{\omega_{1}}^2+ c_{\omega_{2}}^2)(\pi+2^3\omega_2^2)|O_{I}|_{\infty}\alpha^{-2}n_c^{-1}.
\end{align*}
This completes the proof.
\end{proof}

%%%%%%%%%%%%%%%%%%%%%%%%%%%%%%%%%%%%%%%%%%%%%%%%%%%%%%%%%%%%%%%%%%

\section{Generalization of the combined neural networks}\label{sec:gene}
In this section, we shall {give a generalization} analysis of the combined neural networks. For the purpose of analysis, we replace the element-wise Hadamard multiplication in the BFNNs $\phi^0$ and $\phi^1$ in Theorems \ref{thm:BF_0_app} and \ref{thm:BF_1_app} by matrix multiplication introduced in Remark \ref{rem:F*_cos} such that $C\odot \Lambda=-\frac12 (D^* \Lambda D+D \Lambda D^*)$, where $D\in \mathbb{C}^{n_\theta\times n_\theta}$ is a diagonal matrix with ${\rm diag}(D)$ being the discretization of the kernel $k_0(t):=e^{-it}$. 
We assume that the wave frequencies $\omega_1$ and $\omega_2$ are known and the output of the BFNN shares the same size of the input of the CNN, i.e., $n_c=n_\theta=\frac{\omega_1}{\omega_2}n_\rho$, and the linear coordinates converting operator $\psi_0$ (convert the polar coordinates to the Cartesian coordinates) is $2$-Lipschitz continuous such that $\|\psi_0(\cdot)\|_F\leq 2\|\cdot\|_F$.
Furthermore, we consider the CNN $\psi:=\psi_2\circ\psi_1$ as a single operator whose operator norm $\|\psi\|$ is bounded by $2\|(F^*F+\alpha I)^{-1}\|\leq 2\alpha^{-1}$ such that $\|\psi(\cdot)\|_F\leq 2\alpha^{-1}\|\cdot\|_F$.

\subsection{Training set, hypothesis space and loss function}\label{sec:hypothesis}
First, we collect a set of training data $\mathcal{S}:=\{\Lambda_s\in \mathbb{C}^{2n_\theta\times n_\theta}, (\gamma_s,\eta_s)^T\in \mathbb{R}^{2n_c\times n_c}\}_{s=1}^{N}$ with scattering data $\Lambda_s$ being the inputs and the perturbations $(\gamma_s,\eta_s)^T$ being the targets defined in \eqref{eqn:prob_linear_sys} such that $\|\Lambda_s\|_F\leq B_{in}$, $\|(\gamma_s,\eta_s)^T\|_F\leq B_{tgt}$ and $(\Lambda_s,(\gamma_s,\eta_s)^T)\thicksim \mu=\mu_1\times \mu_2$. 
We define the hypothesis spaces, i.e., the parameterized set of all proposed neural networks without (uncompressed $0$-level BFNN) and with (compressed $1$-level BFNN) the butterfly structure as
\begin{align*}
&\mathcal{H}^0:=\{h=\psi\circ\psi_0\circ\phi^0:\|\psi\|\leq 2\alpha^{-1}, \|\psi_0\|\leq2, \phi^0\in\mathcal{H}^0_{BFNN}\}\\
\mbox{and} \quad&\mathcal{H}^1:=\{h=\psi\circ\psi_0\circ\phi^1:\|\psi\|\leq 2\alpha^{-1}, \|\psi_0\|\leq2, \phi^1\in\mathcal{H}^1_{BFNN}\}
\end{align*}
respectively, where 
\begin{align*}
\mathcal{H}^0_{BFNN}=\bigg\{&\phi^0(\cdot)
=\sum_{i=1}^2 (P_u^T (\phi_1^0 R_1^{\omega_i}(\cdot),\cdots,\phi_1^0 R_{n_\theta}^{\omega_i}(\cdot) )^T P_{\omega_i}+ P_l^T (\phi_2^0 R_1^{\omega_i}(\cdot),\cdots,\phi_2^0 R_{n_\theta}^{\omega_i}(\cdot))^T P_{\omega_i}):\\ 
&\phi_1^0 R_m^{\omega}(\Lambda^\omega)=\phi_1^0\Lambda_{\theta_m}^\omega=-\frac{2\pi^2 c_{\omega}}{n_\theta^2}{\rm diag}\big(K^* ( D^*\Lambda_{\theta_m}^{\omega}D+D\Lambda_{\theta_m}^{\omega}D^*) K\big),\\
&\phi_2^0 R_m^{\omega}(\Lambda^\omega)=\phi_2^0\Lambda_{\theta_m}^\omega=\frac{4\pi^2 c_{\omega}}{n_\theta^2}{\rm diag}\big(K^* \Lambda_{\theta_m}^{\omega} K\big),\\
&K\in\mathbb{C}^{n_\theta\times n_\rho},\; D\in\mathbb{C}^{n_\theta\times n_\theta},\;\|K\|_2\leq B_K, \; \|D\|_2\leq B_D \bigg\}\\
\mathcal{H}^1_{BFNN}=\bigg\{&\phi^1(\cdot)
=\sum_{i=1}^2 (P_u^T (\phi_1^1 R_1^{\omega_i}(\cdot),\cdots,\phi_1^1 R_{n_\theta}^{\omega_i}(\cdot) )^T P_{\omega_i}+ P_l^T (\phi_2^1 R_1^{\omega_i}(\cdot),\cdots,\phi_2^1 R_{n_\theta}^{\omega_i}(\cdot))^T P_{\omega_i}):\\ 
&\phi_1^1 R_m^{\omega}(\Lambda^\omega)=\phi_1^1\Lambda_{\theta_m}^\omega=-\frac{2\pi^2 c_{\omega}}{n_\theta^2}{\rm diag}\big(V^*M^*U^* ( D^*\Lambda_{\theta_m}^{\omega}D+D\Lambda_{\theta_m}^{\omega}D^*) UMV\big),\\
&\phi_2^1 R_m^{\omega}(\Lambda^\omega)=\phi_2^1\Lambda_{\theta_m}^\omega=\frac{4\pi^2 c_{\omega}}{n_\theta^2}{\rm diag}\big(V^*M^*U^* \Lambda_{\theta_m}^{\omega} UMV\big),\\
& U=diag(U_1,\cdots,U_{n_r}) \mbox{ with } U_i\in \mathbb{C}^{\frac{n_\theta}{n_r}\times n_r r},  \;V=diag(V_1,\cdots,V_{n_r}) \mbox{ with } V_i\in \mathbb{C}^{n_r r\times \frac{n_\rho}{n_r}},\\
&M=(M_{ij})_{i,j=1}^{n_r} \mbox{ where } M_{ij}=(M_{mn}^{ij})_{m,n=1}^{n_r} \mbox{ with } M_{ji}^{ij}\in \mathbb{C}^{r\times r}, \mbox{ otherwise } M_{mn}^{ij}=O_{r\times r},\\
&D\in\mathbb{C}^{n_\theta\times n_\theta},\;\|U\|_2\leq B_U, \;\|V\|_2\leq B_V,\;\|M\|_2\leq B_M,\;\|D\|_2\leq B_D \bigg\},
\end{align*}
We assume the outputs of the neural networks are bounded by $\|h(\Lambda_s)\|_F\leq B_{out}$.
Given a loss function $\ell(h,\Lambda,(\gamma,\eta)^T):=\|h(\Lambda)-(\gamma,\eta)^T\|_F$, the empirical risk (i.e., the numerical approximation error in $\|\cdot\|_F$) is the average reconstruction error on the training set $\mathcal{S}$ given by $\hat{\mathcal{L}}(h)=\frac1N \sum_{s=1}^N \ell(h,\Lambda_s,(\gamma_s,\eta_s)^T)$ while the true loss (i.e., the generalization error) is given by $\mathcal{L}(h)=\E_{(\Lambda,(\gamma,\eta)^T)\thicksim \mu} [\ell(h,\Lambda,(\gamma,\eta)^T)]$.

\subsection{Generalization error of $\mathcal{H}^0$}\label{sec:gene_0}
We shall start with the hypothesis class $\mathcal{H}^0$. 
The Rademacher complexity is often used to estimate the generalization error $\mathcal{L}(h)$ and the following result collected from \cite[Theorem 26.5]{Shalev-ShwartzBen-David2014} provides a way to estimate the generalization error in terms of the empirical Rademacher complexity.
\begin{lemma}\label{lem:Rcomp_gene}
Let $\mathcal{S}$ be a training set of $N$ pairs of data drawn i.i.d. from the distribution $\mu$, $\mathcal{H}$ be a hypothesis class and $\ell$ be a real-valued loss function such that $|\ell(h,z)|\leq c$ for any $z\in\mathcal{S}$ and $h\in\mathcal{H}$. 
Then for any $\delta\in(0,1)$ and $h\in\mathcal{H}$, with probability of at least $1-\delta$, there holds
\begin{align*}
\mathcal{L}(h)\leq\hat{\mathcal{L}}(h)+2\mathcal{R}_\mathcal{S}(\ell\circ \mathcal{H}\circ \mathcal{S})+4c\sqrt{2\ln(4\delta^{-1})}N^{-\frac12},
\end{align*}
where the empirical Rademacher complexity $\mathcal{R}_\mathcal{S}(\ell\circ \mathcal{H}\circ \mathcal{S})=\E_\epsilon \sup_{h\in\mathcal{H}}\frac{1}{N}\sum_{s=1}^N \epsilon_s \ell(h,z_s)$ with independent Rademacher sequence $\epsilon=(\epsilon_1,\cdots,\epsilon_N)$ taking the value $\pm 1$ with equal probability.
\end{lemma}
\noindent
By bounding the loss function $\ell$, we obtain the next result directly from the above lemma.
\begin{cor}\label{cor:Rcomp_gene}
Let the training set $\mathcal{S}$, the hypothesis class $\mathcal{H}^0$ and the loss function $\ell$ be defined in Section \ref{sec:hypothesis}.
Then for any $\delta\in(0,1)$ and $h\in\mathcal{H}^0$, with probability of at least $1-\delta$, there holds
\begin{align*}
\mathcal{L}(h)\leq\hat{\mathcal{L}}(h)+ 2\mathcal{R}_\mathcal{S}(\ell\circ \mathcal{H}^0\circ \mathcal{S})+4(B_{out}+B_{tgt})\sqrt{2\ln(4\delta^{-1})}N^{-\frac12},
\end{align*}
where the empirical Rademacher complexity $\mathcal{R}_\mathcal{S}(\ell\circ \mathcal{H}^0\circ \mathcal{S})=\E_\epsilon \sup_{h\in\mathcal{H}^0}\frac{1}{N}\sum_{s=1}^N \epsilon_s \|h(\Lambda_s)-(\gamma_s,\eta_s)^T\|_F$.
\end{cor}
\begin{proof}
It follows directly from Lemma \ref{lem:Rcomp_gene} with the estimate 
\begin{align*}
|\ell(h,\Lambda,(\gamma,\eta)^T)|=\|h(\Lambda)-(\gamma,\eta)^T\|_F\leq \|h(\Lambda)\|_F+\|(\gamma,\eta)^T\|_F\leq B_{out}+B_{tgt}:=c.
\end{align*}
\end{proof}

Next, in order to bound the above empirical Rademacher complexity, we state the matrix contraction inequality, i.e., a variant of the vector contraction inequality proved in \cite[Corollary 1]{Maurer2016}, obtained by flattening the matrices into vectors in the next lemma.

\begin{lemma}\label{lem:Rcomp_mci}
Let $\mathcal{S}=\{z_s\}_{s=1}^N\in \mathcal{X}^N$ be a training set of $N$ pairs of data drawn i.i.d. from the distribution $\mu$, $\mathcal{H}$ be a hypothesis class with elements $h: \mathcal{X}\to\mathbb{R}^{n\times n}$ and $f:\mathbb{R}^{n\times n}\to \mathbb{R}$ be a real-valued $L_f$-Lipschitz continuous function. Then
\begin{align*}
\E_\epsilon \sup_{h\in\mathcal{H}}\sum_{s=1}^N \epsilon_s f\circ h(z_s)\leq\sqrt{2}L_f \E_\epsilon \sup_{h\in\mathcal{H}}\sum_{s=1}^N\sum_{i,j=1}^{n}\epsilon_{sij}h_{ij}(z_s)
\end{align*}
where $\epsilon$ is either an independent Rademacher sequence $\epsilon=(\epsilon_{s})_{s=1}^N$ or an independent triply indexed Rademacher sequence $\epsilon=(\epsilon_{sij})_{1\leq s\leq N; 1\leq i,j\leq n}$ and $h_{ij}(z_s)$ is the $(i,j)$-entry of $h(z_s)$. 
\end{lemma}
\noindent
With above lemma, we derive an intermediate result which shows a relationship between the empirical Rademacher complexities of $\ell\circ h$ and $\phi^0$. 

\begin{cor}\label{cor:Rcomp_Rcomp}
Let the training set $\mathcal{S}$, the hypothesis class $\mathcal{H}^0$ and the loss function $\ell$ be defined in Section \ref{sec:hypothesis}. Then there holds
\begin{align*}
\mathcal{R}_\mathcal{S}(\ell\circ \mathcal{H}^0\circ \mathcal{S})\leq 2^4\alpha^{-1}n_c^2 \mathcal{R}_\mathcal{S}(\mathcal{H}^0_{BFNN}\circ \mathcal{S})
\end{align*}
where the empirical Rademacher complexity \begin{align*}
\mathcal{R}_\mathcal{S}(\mathcal{H}^0_{BFNN}\circ \mathcal{S})=\E_\epsilon \sup_{\phi^0\in\mathcal{H}^0_{BFNN}}\frac{1}{N}\sum_{s=1}^N\sum_{i=1}^{2n_c}\sum_{j=1}^{n_c}\epsilon_{sij}\phi^0_{ij}(\Lambda_s).
\end{align*}
\end{cor}
\begin{proof}
By Lemma \ref{lem:Rcomp_mci} with a $1$-Lipschitz continuous function $f:\mathbb{R}^{2n_c\times n_c}\to \mathbb{R}$ given by $f(\cdot):=\|\cdot\|_F$, there holds
\begin{align*}
\mathcal{R}_\mathcal{S}(\ell\circ \mathcal{H}^0\circ \mathcal{S})=&\E_\epsilon \sup_{h\in\mathcal{H}^0}\frac{1}{N}\sum_{s=1}^N \epsilon_s \|h(\Lambda_s)-(\gamma_s,\eta_s)^T\|_F\\
\leq&\sqrt{2} \E_\epsilon \sup_{h\in\mathcal{H}^0}\frac{1}{N}\sum_{s=1}^N\sum_{i=1}^{2n_c}\sum_{j=1}^{n_c}\epsilon_{sij}(h_{ij}(\Lambda_s)-(\gamma_s,\eta_s)^T_{ij})\\
=&\sqrt{2} \E_\epsilon \sup_{h\in\mathcal{H}^0}\frac{1}{N}\sum_{s=1}^N\sum_{i=1}^{2n_c}\sum_{j=1}^{n_c}\epsilon_{sij}h_{ij}(\Lambda_s):=\sqrt{2} \mathcal{R}_\mathcal{S}( \mathcal{H}^0\circ \mathcal{S}).
\end{align*}
By the definition of the hypothesis class $\mathcal{H}^0$, we rewrite $\mathcal{R}_\mathcal{S}( \mathcal{H}^0\circ \mathcal{S})$ as
\begin{align*}
\mathcal{R}_\mathcal{S}(\mathcal{H}^0\circ \mathcal{S})=\E_\epsilon \sup_{\|\psi\|\leq 2\alpha^{-1},\|\psi_0\|\leq 2}\;\sup_{\phi^0\in\mathcal{H}^0_{BFNN}}\frac{1}{N}\sum_{s=1}^N\sum_{i=1}^{2n_c}\sum_{j=1}^{n_c}\epsilon_{sij}[\psi\circ\psi_0]_{ij}(\phi^0(\Lambda_s))
\end{align*}
where $[\psi\circ\psi_0]_{ij}: \cdot \mapsto (\psi\circ\psi_0(\cdot))[i,j]\in \mathbb{R}$ such that $$\|[\psi\circ\psi_0]_{ij}(\cdot)\|_F=|[\psi\circ\psi_0]_{ij}(\cdot)|\leq \|\psi\circ\psi_0(\cdot)\|_F\leq 4\alpha^{-1}\|\cdot\|_F$$ is a $(4\alpha^{-1})$-Lipschitz continuous function for any $i,j$. Thus, by Lemma \ref{lem:Rcomp_mci} with $f=[\psi\circ\psi_0]_{ij}$, we have
\begin{align*}
\mathcal{R}_\mathcal{S}(\mathcal{H}^0\circ \mathcal{S})\leq&\sum_{i=1}^{2n_c}\sum_{j=1}^{n_c}\E_\epsilon \sup_{\|[\psi\circ\psi_0]_{ij}\|\leq 4\alpha^{-1}}\;\sup_{\phi^0\in\mathcal{H}^0_{BFNN}}\frac{1}{N}\sum_{s=1}^N\epsilon_{s}[\psi\circ\psi_0]_{ij}(\phi^0(\Lambda_s))\\
\leq& 2^3\sqrt{2}\alpha^{-1}n_c^2\E_\epsilon \sup_{\phi^0\in\mathcal{H}^0_{BFNN}}\frac{1}{N}\sum_{s=1}^N\sum_{i=1}^{2n_c}\sum_{j=1}^{n_c}\epsilon_{sij}\phi^0_{ij}(\Lambda_s):= 2^3\sqrt{2}\alpha^{-1}n_c^2 \mathcal{R}_\mathcal{S}(\mathcal{H}^0_{BFNN}\circ \mathcal{S})
\end{align*}
which completes the proof.
\end{proof}

It is worth noting that the Rademacher process is centered stochastic and has sub-Gaussian increments \cite[Definition 8.22]{FoucartRauhut2013}, where the Dudley's inequality \cite[Theorem 8.23]{FoucartRauhut2013} holds, and thus we can bound the empirical Rademacher complexity $\mathcal{R}_\mathcal{S}(\mathcal{H}^0_{BFNN}\circ \mathcal{S})$ in terms of covering numbers associated with the pseudo-metric $\|\phi^0-\tilde{\phi^0}\|:=\sqrt{\sum_{s=1}^{N}\|\phi^0(\Lambda_s)-\tilde{\phi^0}(\Lambda_s)\|_F^2}$ \cite[(8.44)]{FoucartRauhut2013}.
To this end, we assume that the intermediate output $\|\phi^0(\Lambda_s)\|_F\leq B_{out}$ for any possible hypothesis $\phi^0\in\mathcal{H}^0_{BFNN}$ and define the intermediate output tensor set which consists of the outputs of any $\phi_0\in\mathcal{H}^0_{BFNN}$ applied to the training set $\mathcal{S}$ by $$\mathcal{M}_{\mathcal{H}^0_{BFNN}}:=\big\{\big(\phi^0(\Lambda_s)\big)_{s=1}^N\in \mathbb{R}^{2n_c\times n_c \times N}: \phi^0\in \mathcal{H}^0_{BFNN}\big\}.$$ 
The covering number of $\mathcal{M}_{\mathcal{H}^0_{BFNN}}$, denoted by $\mathcal{N}(\mathcal{M}_{\mathcal{H}^0_{BFNN}},\|\cdot\|,\epsilon)$, is defined as the smallest size of the subset $\mathcal{M} \subset \mathcal{M}_{\mathcal{H}^0_{BFNN}}$ such that $\min_{z\in\mathcal{M}}\|z-z'\|\leq\epsilon$ for all $z'\in\mathcal{M}_{\mathcal{H}^0_{BFNN}}$.
We state the result obtained directly from the Dudley's inequality in the following lemma.
\begin{lemma}\label{lem:dudley}
Let $\mathcal{R}_\mathcal{S}(\mathcal{H}^0_{BFNN}\circ \mathcal{S})$ be the empirical Rademacher complexity defined in Corollary \ref{cor:Rcomp_Rcomp}. 
Then, we have
\begin{align*}
\mathcal{R}_\mathcal{S}(\mathcal{H}^0_{BFNN}\circ \mathcal{S})\leq 12N^{-1}\int_{0}^{\sqrt{N}B_{out}}\sqrt{\ln \big(\mathcal{N}(\mathcal{M}_{\mathcal{H}^0_{BFNN}},\|\cdot\|,\epsilon)\big)}{\rm d}\epsilon.
\end{align*}
\end{lemma}
\begin{proof}
The Dudley's inequality gives the estimate
\begin{align*}
\mathcal{R}_\mathcal{S}(\mathcal{H}^0_{BFNN}\circ \mathcal{S})=&\frac1N\E_\epsilon \sup_{\phi^0\in\mathcal{H}^0_{BFNN}}\sum_{s=1}^N\sum_{i=1}^{2n_c}\sum_{j=1}^{n_c}\epsilon_{sij}\phi^0_{ij}(\Lambda_s)\\
\leq& \frac{12}{N}\int_{0}^{\frac12 \Delta(\mathcal{M}_{\mathcal{H}^0_{BFNN}})}\sqrt{\ln \big(\mathcal{N}(\mathcal{M}_{\mathcal{H}^0_{BFNN}},\|\cdot\|,\epsilon)\big)}{\rm d}\epsilon,
\end{align*}
where the diameter of $\mathcal{M}_{\mathcal{H}^0_{BFNN}}$ is given by 
\begin{align*}
\Delta(\mathcal{M}_{\mathcal{H}^0_{BFNN}}):= \sup_{\phi^0,\tilde{\phi^0}\in \mathcal{H}^0_{BFNN}}\|\phi^0-\tilde{\phi^0}\|=\sup_{\phi^0,\tilde{\phi^0}\in \mathcal{H}^0_{BFNN}}\sqrt{\sum_{s=1}^{N}\|\phi^0(\Lambda_s)-\tilde{\phi^0}(\Lambda_s)\|_F^2}.
\end{align*}
Then the inequalities
\begin{align*}
\sum_{s=1}^{N}\|\phi^0(\Lambda_s)-\tilde{\phi^0}(\Lambda_s)\|_F^2\leq 2\sum_{s=1}^{N}(\|\phi^0(\Lambda_s)\|_F^2+\|\tilde{\phi^0}(\Lambda_s)\|_F^2)\leq 4NB_{out}^2, \quad \forall \phi^0, \tilde{\phi^0}\in \mathcal{H}^0_{BFNN}
\end{align*}
{complete} the proof.
\end{proof}

The remaining part of the analysis is to bound the covering number $\mathcal{N}(\mathcal{M}_{\mathcal{H}^0_{BFNN}},\|\cdot\|,\epsilon)$ by giving the Lipschitz estimates of the neural networks in $\mathcal{H}^0_{BFNN}$ with respect to the parameters 
\begin{equation}\label{eqn:para_set}
\Theta_0:=\{K\in\mathcal{K}\subset B_K \mathcal{B}_K,\; D\in\mathcal{D}\subset B_D \mathcal{B}_D\}:=\mathcal{K}\times \mathcal{D} 
\end{equation} 
with the unit balls given by $\mathcal{B}_K=\{K\in\mathbb{C}^{n_\theta \times n_\rho} \backsimeq \mathbb{R}^{2n_\theta n_\rho}: \|K\|_2\leq 1\}$ and $\mathcal{B}_D=\{D\in\mathbb{C}^{n_\theta \times n_\theta}\backsimeq \mathbb{R}^{2n_\theta^2}: \|D\|_2\leq 1\}$. 
Based on the Lipschitz estimates in Lemma \ref{lem:cover_lip_0} that
\begin{align*}
\|\phi^0-\tilde{\phi^0}\|\leq (c_{K}\|K-\tilde{K}\|_2+c_{D}\|D-\tilde{D}\|_2)n_c^{-\frac32}\sqrt{N}
\end{align*}
with the constants $c_{K}=2^3\pi^2 (c_{\omega_1}+c_{\omega_2})(B_D^4+1)^\frac12 B_K B_{in}$ and $c_{D}=2^3\pi^2 (c_{\omega_1}+c_{\omega_2})B_DB_K^2 B_{in}$,
we can bound the covering number of $\mathcal{M}_{\mathcal{H}^0_{BFNN}}$ by the covering number of the parameter product set $\Theta_0=\mathcal{K}\times \mathcal{D}$ with the metric $\|(K,D)\|_{\Theta_0}:=c_{K}\|K\|_2+c_{D}\|D\|_2$ as
\begin{equation}\label{eqn:cover_cover_convert}
\mathcal{N}(\mathcal{M}_{\mathcal{H}^0_{BFNN}},\|\cdot\|,\epsilon)\leq\mathcal{N}(\Theta_0,\|\cdot\|_{\Theta_0},n_c^{\frac32}N^{-\frac12}\epsilon)=\mathcal{N}(\mathcal{K}\times \mathcal{D},\|\cdot\|_{\Theta_0}, n_c^{\frac32}N^{-\frac12}\epsilon).
\end{equation}
To estimate $\mathcal{N}(\mathcal{K}\times \mathcal{D},\|\cdot\|_{\Theta_0},n_c^{\frac32}N^{-\frac12}\epsilon)$, we state Lemma \ref{lem:cover_product} (i.e., \cite[Lemma 5.3]{SchnoorBehboodiRauhut2023}), which provides the relationship between the covering number of a product set and its component sets, followed by Lemma \ref{lem:cover_ball} (e.g., \cite[Lemma 5.2]{SchnoorBehboodiRauhut2023}) which can be applied to derive the covering number of each component.
\begin{lemma}\label{lem:cover_product}
Let $(S_1, d_1), \cdots, (S_p, d_p)$ be $p$ metric spaces, $c_1, \cdots,c_p$ be $p$ positive constants and $S=(S_1\times\cdots \times S_p,d)$ be the product space equipped with the metric $d(z,z')=\sum_{i=1}^{p}c_id_i(z_i,z_i')$ for $z=(z_1,\cdots,z_p),z'=(z_1',\cdots,z_p')\in S$. Then for any $\epsilon>0$, the covering number of $S$ can be bounded by
\begin{align*}
\mathcal{N}(S,d,\epsilon)\leq\Pi_{i=1}^p \mathcal{N}(S_i,d_i,\frac{\epsilon}{c_i p}).
\end{align*}
\end{lemma}
\begin{lemma}\label{lem:cover_ball}
Let $S\subset\{z\in\mathbb{R}^{n}:d(z,0)\leq 1\}$ be a subset contained in the unit ball with an associate metric $d$. Then for any $\epsilon>0$, there holds
\begin{align*}
\mathcal{N}(S,d,\epsilon)\leq\Big(1+\frac{2}{\epsilon}\Big)^n.
\end{align*}
\end{lemma}

\noindent
The above two lemmas yield the next result on the upper bound of $\mathcal{N}(\Theta_0,\|\cdot\|_{\Theta_0},n_c^{\frac32}N^{-\frac12}\epsilon)$.
\begin{cor}\label{cor:cover_integral}
Let $\Theta_0=\mathcal{K}\times \mathcal{D}$ be the parameter set, equipped with the metric $\|\cdot\|_{\Theta_0}$, defined in \eqref{eqn:para_set}. Then for any $\epsilon>0$, there holds
\begin{align*}
\mathcal{N}(\Theta_0,\|\cdot\|_{\Theta_0}, n_c^{\frac32}N^{-\frac12}\epsilon)
%\leq \Big(1+\frac{4B_K c_{K}}{\epsilon}n_c^{-\frac32}\sqrt{N}\Big)^{2\frac{\omega_2}{\omega_1}n_c^2}
%\Big(1+\frac{4B_D c_{D}}{\epsilon}n_c^{-\frac32}\sqrt{N}\Big)^{2n_c^2}
\leq\Big(1+\frac{c_0}{\epsilon}n_c^{-\frac32}\sqrt{N}\Big)^{2(\frac{\omega_2}{\omega_1}+1)n_c^2},
\end{align*}
where $c_0=2^5\pi^2 (c_{\omega_1}+c_{\omega_2})(B_D^4+1)^\frac12 B_K^2 B_{in}$.
\end{cor}
\begin{proof}
By setting $S_1=\mathcal{K}$, $S_2=\mathcal{D}$, $c_1=c_{K}$, $c_2=c_{D}$ and $d_1=d_2=\|\cdot\|_2$ in Lemma \ref{lem:cover_product}, 
we obtain
\begin{align*}
\mathcal{N}(\Theta_0,\|\cdot\|_{\Theta_0},n_c^{\frac32}N^{-\frac12}\epsilon)%=&\mathcal{N}(\mathcal{K}\times \mathcal{D},\|\cdot\|_{\Theta_0},n_c^{\frac32}N^{-\frac12}\epsilon)\\
\leq& \mathcal{N}(\mathcal{K},\|\cdot\|_2,\frac{\epsilon}{2c_{K}}n_c^{\frac32}N^{-\frac12})\mathcal{N}(\mathcal{D},\|\cdot\|_2,\frac{\epsilon}{2c_{D}}n_c^{\frac32}N^{-\frac12})\\
\leq & \mathcal{N}(B_K \mathcal{B}_K,\|\cdot\|_2,\frac{\epsilon}{2c_{K}}n_c^{\frac32}N^{-\frac12})\mathcal{N}(B_D \mathcal{B}_D,\|\cdot\|_2,\frac{\epsilon}{2c_{D}}n_c^{\frac32}N^{-\frac12})\\
\leq & \mathcal{N}(\mathcal{B}_K,\|\cdot\|_2,\frac{\epsilon}{2B_K c_{K}}n_c^{\frac32}N^{-\frac12})\mathcal{N}( \mathcal{B}_D,\|\cdot\|_2,\frac{\epsilon}{2B_Dc_{D}}n_c^{\frac32}N^{-\frac12}).
\end{align*}
Then, with Lemma \ref{lem:cover_ball} and the definitions of $\mathcal{B}_K$ and $\mathcal{B}_D$, there hold
\begin{align*}
\mathcal{N}(\mathcal{B}_K,\|\cdot\|_2,\frac{\epsilon}{2B_K c_{K}}n_c^{\frac32}N^{-\frac12})\leq &\Big(1+\frac{4B_K c_{K}}{\epsilon}n_c^{-\frac32}\sqrt{N}\Big)^{2n_\theta n_\rho},\\
\mathcal{N}(\mathcal{B}_D,\|\cdot\|_2,\frac{\epsilon}{2B_D c_{D}}n_c^{\frac32}N^{-\frac12})\leq &\Big(1+\frac{4B_D c_{D}}{\epsilon}n_c^{-\frac32}\sqrt{N}\Big)^{2n_\theta^2}.
\end{align*}
Finally, the assumption $n_c=n_\theta=\frac{\omega_1}{\omega_2}n_\rho$ and the notation $$c_0=4\max(B_Kc_K,B_Dc_D)=2^5\pi^2 (c_{\omega_1}+c_{\omega_2})(B_D^4+1)^\frac12 B_K^2 B_{in}$$ give the assertion in the corollary.
\end{proof}

\noindent
To estimate the integrals in Corollary \ref{cor:cover_integral}, we collect the following lemma from \cite[Lemma 5.4]{SchnoorBehboodiRauhut2023}.
\begin{lemma}\label{lem:integral}
For any $\beta_1,\beta_2>0$, there holds
\begin{align*}
\int_0^{\beta_1}\sqrt{\ln \Big(1+\frac{\beta_2}{\epsilon}\Big)}{\rm d}\epsilon\leq \beta_1 \Phi\Big(\frac{\beta_2}{\beta_1}\Big),
\end{align*}
where $\Phi(t):=\sqrt{\ln (1+t)+t\ln(1+\frac{1}{t})}\leq \sqrt{\ln(1+t)+1}$ for any $t>0$ and $\lim_{t\to 0^+}\Phi(t)=0$.
\end{lemma}

Now, based on all above corollaries, Lemma \ref{lem:dudley} and Lemma \ref{lem:integral}, we prove  the first main result on the generalization analysis.
\begin{theorem}[Generalization of uncompressed neural networks]\label{thm:gen_0}
Let the training set $\mathcal{S}$, the hypothesis class $\mathcal{H}^0$ and the loss function $\ell$ be defined in Section \ref{sec:hypothesis}.
Then for any $\delta\in(0,1)$ and $h\in\mathcal{H}^0$, with probability of at least $1-\delta$, there holds
\begin{align*}
\mathcal{L}(h)\leq \hat{\mathcal{L}}(h)+2^2\Big(\sqrt{2\ln(4\delta^{-1})}(B_{out}+B_{tgt})+2^6 3\sqrt{\omega_2/\omega_1}B_{out}\alpha^{-1} (1+c_0 B_{out}^{-1} n_c^{-\frac32})^{\frac12}n_c^3 \Big)N^{-\frac12},
\end{align*}
with the constant $c_0=2^5\pi^2 (c_{\omega_1}+c_{\omega_2})(B_D^4+1)^\frac12 B_K^2 B_{in}$.
\end{theorem}

\begin{proof}
It follows directly from Corollaries \ref{cor:Rcomp_gene}, \ref{cor:Rcomp_Rcomp}, Lemma \ref{lem:dudley}, equation \eqref{eqn:cover_cover_convert} and Corollary \ref{cor:cover_integral} that
\begin{align*}
\mathcal{L}(h)\leq& \hat{\mathcal{L}}(h)+4(B_{out}+B_{tgt})\sqrt{2\ln(4\delta^{-1})}N^{-\frac12}+{\rm I}_\mathcal{N}^0,
\end{align*}
where 
\begin{align*}
{\rm I}_\mathcal{N}^0:=&\frac{2^7 3n_c^2}{\alpha N} \int_{0}^{\sqrt{N}B_{out}}\sqrt{\ln \Big(\big(1+\frac{c_0}{\epsilon}n_c^{-\frac32}\sqrt{N}\big)^{2(\frac{\omega_2}{\omega_1}+1)n_c^2}\Big)}{\rm d}\epsilon\\
%=&\frac{2^7 3n_c^2}{\alpha N} \int_{0}^{\sqrt{N}B_{out}}\sqrt{{2(\frac{\omega_2}{\omega_1}+1)n_c^2}\ln \big(1+\frac{c_0}{\epsilon}n_c^{-\frac32}\sqrt{N}\big)}{\rm d}\epsilon\\
=&\frac{2^7 3\sqrt{2(\omega_2/\omega_1+1)}n_c^3}{\alpha N} \int_{0}^{\sqrt{N}B_{out}}\sqrt{\ln \big(1+\frac{c_0}{\epsilon}n_c^{-\frac32}\sqrt{N}\big)}{\rm d}\epsilon.
\end{align*}
Further, with Lemma \ref{lem:integral} and the assumption $\omega_2>\omega_1$, we bound ${\rm I}_\mathcal{N}^0$ by
\begin{align*}
{\rm I}_\mathcal{N}^0\leq&\frac{2^8 3\sqrt{\omega_2/\omega_1}n_c^3}{\alpha N} \sqrt{N}B_{out}\Phi\Big(\frac{c_0 n_c^{-\frac32}\sqrt{N}}{\sqrt{N}B_{out}}\Big)
=\frac{2^8 3\sqrt{\omega_2/\omega_1}B_{out} n_c^3}{ \alpha\sqrt{N}} \Phi\Big(c_0 B_{out}^{-1} n_c^{-\frac32}\Big)\\
\leq&\frac{2^8 3\sqrt{\omega_2/\omega_1}B_{out} n_c^3}{ \alpha\sqrt{N}} \sqrt{1+\ln(1+c_0 B_{out}^{-1} n_c^{-\frac32})}
\leq 2^8 3\sqrt{\omega_2/\omega_1}B_{out}\alpha^{-1} (1+c_0 B_{out}^{-1} n_c^{-\frac32})^{\frac12}n_c^3 N^{-\frac12}.
\end{align*}
This completes the proof.
\end{proof}

\subsection{Generalization error of $\mathcal{H}^1$}\label{sec:gene_1}
We extent the result on the generalization of the uncompressed hypothesis space $\mathcal{H}^0$ in Section \ref{sec:gene_0} to the compressed hypothesis space $\mathcal{H}^1$ where the parameter $K$ is further approximated by $K_r=UMV$ and prove the second main result on the generalization analysis.

\begin{theorem}[Generalization of compressed neural networks]\label{thm:gen_1}
Let the training set $\mathcal{S}$, the hypothesis class $\mathcal{H}^1$ and the loss function $\ell$ be defined in Section \ref{sec:hypothesis} with $\frac{\omega_2}{\omega_1} n_r r\leq n_c$.
Then for any $\delta\in(0,1)$ and $h\in\mathcal{H}^0$, with probability of at least $1-\delta$, there holds
\begin{align*}
\mathcal{L}(h)\leq \hat{\mathcal{L}}(h)+2^2\Big(\sqrt{2\ln(4\delta^{-1})}(B_{out}+B_{tgt})+2^6 3 \sqrt{2}B_{out} \alpha^{-1}\big(1+c_1 B_{out}^{-1} n_c^{-\frac32}\big)^\frac12 n_c^3\Big)N^{-\frac12},
\end{align*}
with the constant $c_1:=2^6\pi^2 (c_{\omega_1}+c_{\omega_2})(B_D^4+1)^\frac12 B_U^2 B_M^2 B_V^2 B_{in}$.   
\end{theorem}

\begin{proof}
By Corollaries \ref{cor:Rcomp_gene}, \ref{cor:Rcomp_Rcomp} and Lemma \ref{lem:dudley} for the hypothesis space $\mathcal{H}^1$,
we have
\begin{align*}
\mathcal{L}(h)\leq& \hat{\mathcal{L}}(h)+4(B_{out}+B_{tgt})\sqrt{2\ln(4\delta^{-1})}N^{-\frac12}+{\rm I}_\mathcal{N}^1,
\end{align*}
where 
\begin{align*}
{\rm I}_\mathcal{N}^1\leq 2^5\alpha^{-1}n_c^2 \mathcal{R}_\mathcal{S}(\mathcal{H}^1_{BFNN}\circ \mathcal{S})\leq 2^7 3\alpha^{-1}n_c^2 N^{-1}\int_{0}^{\sqrt{N}B_{out}}\sqrt{\ln \big(\mathcal{N}(\mathcal{M}_{\mathcal{H}^1_{BFNN}},\|\cdot\|,\epsilon)\big)}{\rm d}\epsilon.
\end{align*}
Based on the Lipschitz estimates in Lemma \ref{lem:cover_lip_1} that
\begin{align*}
\|\phi^1-\tilde{\phi^1}\|\leq (c_{U} \|U-\tilde{U}\|_2+c_{M}\|M-\tilde{M}\|_2+c_{V}\|V-\tilde{V}\|_2+c_{D_r}\|D-\tilde{D}\|_2)n_c^{-\frac32}\sqrt{N},
\end{align*}
where the constants $c_{U}=c_{K_r}B_M B_V$,  $c_{M}=c_{K_r}B_U B_V$, $c_{V}=c_{K_r}B_U B_M$ and $c_{D_r}=2^3\pi^2 (c_{\omega_1}+c_{\omega_2})B_D B_U^2 B_M^2 B_V^2 B_{in}$ with $c_{K_r}=2^3\pi^2 (c_{\omega_1}+c_{\omega_2})(B_D^4+1)^\frac12 B_U B_M B_V B_{in}$,
we bound the covering number of $\mathcal{M}_{\mathcal{H}^1_{BFNN}}$ by the covering number of the parameter product set 
\begin{align*}
\Theta_1=\mathcal{U}\times\mathcal{M}\times \mathcal{V}\times\mathcal{D},\quad \mbox{where}\quad \mathcal{U}\subset B_U \mathcal{B}_U, \mathcal{M}\subset B_M \mathcal{B}_M, \mathcal{V}\subset B_V \mathcal{B}_V \; \mbox{and} \;\mathcal{D}\subset B_D \mathcal{B}_D,
\end{align*}
with the unit balls 
\begin{align*}
\mathcal{B}_U=&\{U=diag(U_1,\cdots,U_{n_r})\in\mathbb{C}^{\frac{n_\theta}{n_r}\times n_r r}\times \cdots \times \mathbb{C}^{\frac{n_\theta}{n_r}\times n_r r}\backsimeq \mathbb{R}^{2n_\theta n_r r}: \|U\|_2\leq 1\},\\
\mathcal{B}_V=&\{V=diag(V_1,\cdots,V_{n_r})\in\mathbb{C}^{n_r r\times \frac{n_\rho}{n_r}}\times \cdots \times \mathbb{C}^{n_r r\times \frac{n_\rho}{n_r}}\backsimeq \mathbb{R}^{2n_\rho n_r r}: \|V\|_2\leq 1\},\\
\mathcal{B}_M=&\{M=(M_{ij}^{ji})_{i,j=1}^{n_r} \in\mathbb{C}^{r\times r}\times \cdots \times \mathbb{C}^{r\times r}\backsimeq \mathbb{R}^{2n_r^2 r^2}: \|M\|_2\leq 1\},\\
\mathcal{B}_D=&\{D\in\mathbb{C}^{n_\theta \times n_\theta}\backsimeq \mathbb{R}^{2n_\theta^2}: \|D\|_2\leq 1\}
\end{align*}
and the metric $\|(U,M,V,D)\|_{\Theta_1}:=c_{U}\|U\|_2+c_{M}\|M\|_2+c_{V}\|V\|_2+c_{D_r}\|D\|_2$, as
\begin{equation}\label{eqn:cover_cover_convert_1}
\mathcal{N}(\mathcal{M}_{\mathcal{H}^1_{BFNN}},\|\cdot\|,\epsilon)\leq\mathcal{N}(\Theta_1,\|\cdot\|_{\Theta_1},n_c^{\frac32}N^{-\frac12}\epsilon)=\mathcal{N}(\mathcal{U}\times\mathcal{M}\times \mathcal{V}\times\mathcal{D},\|\cdot\|_{\Theta_1}, n_c^{\frac32}N^{-\frac12}\epsilon).
\end{equation}
{By following} the same strategy used in the proof of Corollary \ref{cor:cover_integral}, we derive
\begin{align*}
\mathcal{N}(\Theta_1,\|\cdot\|_{\Theta_1},n_c^{\frac32}N^{-\frac12}\epsilon)
\leq& \mathcal{N}(\mathcal{U},\|\cdot\|_2,\frac{\epsilon}{4c_{U}}n_c^{\frac32}N^{-\frac12})\mathcal{N}(\mathcal{M},\|\cdot\|_2,\frac{\epsilon}{4c_{M}}n_c^{\frac32}N^{-\frac12})\\
&\;\;\mathcal{N}(\mathcal{V},\|\cdot\|_2,\frac{\epsilon}{4c_{V}}n_c^{\frac32}N^{-\frac12})\mathcal{N}(\mathcal{D},\|\cdot\|_2,\frac{\epsilon}{4c_{D_r}}n_c^{\frac32}N^{-\frac12})\\
\leq & \mathcal{N}(\mathcal{B}_U,\|\cdot\|_2,\frac{\epsilon}{4B_U c_{U}}n_c^{\frac32}N^{-\frac12})\mathcal{N}(\mathcal{B}_M,\|\cdot\|_2,\frac{\epsilon}{4B_M c_{M}}n_c^{\frac32}N^{-\frac12})\\
&\;\;\mathcal{N}(\mathcal{B}_V,\|\cdot\|_2,\frac{\epsilon}{4B_V c_{V}}n_c^{\frac32}N^{-\frac12})\mathcal{N}( \mathcal{B}_D,\|\cdot\|_2,\frac{\epsilon}{4B_Dc_{D_r}}n_c^{\frac32}N^{-\frac12})\\
\leq& \Big(1+\frac{8 B_U c_{U}}{\epsilon}n_c^{-\frac32}\sqrt{N}\Big)^{2n_\theta n_r r} \Big(1+\frac{8 B_M c_{M}}{\epsilon}n_c^{-\frac32}\sqrt{N}\Big)^{2n_r^2 r^2}\\
&\;\;\Big(1+\frac{8 B_V c_{V}}{\epsilon}n_c^{-\frac32}\sqrt{N}\Big)^{2n_\rho n_r r}\Big(1+\frac{8 B_D c_{D_r}}{\epsilon}n_c^{-\frac32}\sqrt{N}\Big)^{2n_\theta^2}.
\end{align*}
With the assumption $n_c=n_\theta=\frac{\omega_1}{\omega_2}n_\rho$ and the notation
\begin{align*}
c_{1}:=&8\max(B_Uc_{U},B_Mc_{M},B_Vc_{V},B_Dc_{D_r})=8\max(B_UB_MB_Vc_{K_r},B_D c_{D_r})\\
=&8B_UB_MB_Vc_{K_r}=2^6\pi^2 (c_{\omega_1}+c_{\omega_2})(B_D^4+1)^\frac12 B_U^2 B_M^2 B_V^2 B_{in},
\end{align*}
we simplify the above estimates of $\mathcal{N}(\Theta_1,\|\cdot\|_{\Theta_1},n_c^{\frac32}N^{-\frac12}\epsilon)$ and derive
\begin{align*}
\mathcal{N}(\Theta_1,\|\cdot\|_{\Theta_1},n_c^{\frac32}N^{-\frac12}\epsilon)
\leq&\Big(1+\frac{c_{1}}{\epsilon}n_c^{-\frac32}\sqrt{N}\Big)^{2(n_c n_r r+n_r^2 r^2+\frac{\omega_2}{\omega_1}n_c n_r r+n_c^2)}.
\end{align*}
Finally, with Lemma \ref{lem:integral} and the assumptions $\omega_2>\omega_1$ and $\frac{\omega_2}{\omega_1} n_r r\leq n_c$, we can bound ${\rm I}_\mathcal{N}^1$ by
\begin{align*}
{\rm I}_\mathcal{N}^1\leq &\frac{2^7 3n_c^2}{\alpha N}\int_{0}^{\sqrt{N}B_{out}}\sqrt{\ln \Big(\big(1+\frac{c_{1}}{\epsilon}n_c^{-\frac32}\sqrt{N}\big)^{2(n_c n_r r+n_r^2 r^2+\frac{\omega_2}{\omega_1}n_c n_r r+n_c^2)}\Big)}{\rm d}\epsilon\\
%\leq &\frac{2^7 3n_c^2}{\alpha N}\int_{0}^{\sqrt{N}B_{out}}\sqrt{2(n_c n_r r+n_r^2 r^2+\frac{\omega_2}{\omega_1}n_c n_r r+n_c^2)\ln \big(1+\frac{c_{1}}{\epsilon}n_c^{-\frac32}\sqrt{N}\big)}{\rm d}\epsilon\\
\leq &\frac{2^7 3 \sqrt{2} n_c^2}{\alpha N}\sqrt{n_r^2 r^2+2(\omega_2/\omega_1) n_c n_r r+n_c^2}\int_{0}^{\sqrt{N}B_{out}}\sqrt{\ln \big(1+\frac{c_{1}}{\epsilon}n_c^{-\frac32}\sqrt{N}\big)}{\rm d}\epsilon\\
%\leq &\frac{2^8 3 \sqrt{2} n_c^3}{\alpha N}\sqrt{N}B_{out}\Phi\Big(\frac{c_1 n_c^{-\frac32}\sqrt{N}}{\sqrt{N}B_{out}}\Big)\\
\leq&\frac{2^8 3 \sqrt{2} n_c^3}{\alpha N}\sqrt{N}B_{out}\Phi\big(c_1 B_{out}^{-1} n_c^{-\frac32}\big)
\leq \frac{2^8 3 \sqrt{2} n_c^3}{\alpha N}\sqrt{N}B_{out}\sqrt{1+\ln(1+c_1 B_{out}^{-1} n_c^{-\frac32})}\\
\leq& 2^8 3 \sqrt{2}B_{out} \alpha^{-1}\big(1+c_1 B_{out}^{-1} n_c^{-\frac32}\big)^\frac12 n_c^3 N^{-\frac12}.
\end{align*}
This completes the proof.
\end{proof}
%%%%%%%%%%%%%%%%%%%%%%%%%%%%%

\section{Numerical experiments}\label{num}

In this section, we present preliminary numerical results to showcase the feasibility of the proposed combined neural networks. 

\subsection{Training setup and neural networks}\label{sec:train}
In the numerical experiments below, we consider a unit disc $\Omega$.
The target and input training datasets consist of perturbations of the inhomogeneity in $\Omega$ and corresponding scattering data at several different wave frequencies. 
The perturbations of interest, i.e., $\gamma=a-1$ and $\eta=n-1$, are generated by
\begin{enumerate}
\item[(i)] \textbf{[Smooth perturbations]} The weighted sums of %two-dimensional 
Gaussian functions combined with a mask, given by 
\begin{equation}\label{eqn:dataset_s}
\gamma=\chi_\Omega \sum_{j=1}^{J} c_{a,j}e^{-\frac{(x - x_{a,j})^2 + (y - y_{a,j})^2}{2\sigma_{a,j}^2}} \quad \mbox{and} \quad \eta=\chi_\Omega\sum_{j=1}^{J} c_{n,j}e^{-\frac{(x - x_{n,j})^2 + (y - y_{n,j})^2}{2\sigma_{n,j}^2}},
\end{equation}
where $\chi_\Omega$ is the indicator function of $\Omega$, the peaks $$(x_{a,j},y_{a,j}) = (r_{a,j}\cos\theta_{a,j},r_{a,j}\sin\theta_{a,j})\quad \mbox{and}\quad(x_{n,j},y_{n,j}) = (r_{n,j}\cos\theta_{n,j},r_{n,j}\sin\theta_{n,j})$$ locate inside the circle of radius $0.5$ centered at the origin, the standard deviations $\sigma_{a,j} = R_{a,j}(8\ln2)^{-\frac12}$ and $\sigma_{n,j} = R_{n,j}(8\ln2)^{-\frac12}$,
with $c_{a,j},c_{n,j}\in [0,J^{-1})$, $r_{a,j},r_{n,j}\in[0,0.5)$, $\theta_{a,j},\theta_{n,j}\in[0,2\pi)$, $R_{a,j}\in \big[0,1-\max(|x_{a,j}|,|y_{a,j}|)\big)$ and $R_{n,j}\in \big[0,1-\max(|x_{n,j}|,|y_{n,j}|)\big)$ randomly sampled from uniform distributions. The bounds for $r_{a,j},r_{n,j}$ and $R_{a,j},R_{n,j}$ guarantee the significant part of the Gaussian functions are inside the unit circle and the bounds for $c_{a,j},c_{n,j}$ ensure that the values of $a$ and $n$ are not too far away from $1$;

We also consider a sub-training dataset where $R_{a,j}\in  \big(1-\max(|x_{a,j}|,|y_{a,j}|)\big)[0.3,1)$ and $R_{n,j}\in \big(1-\max(|x_{n,j}|,|y_{n,j}|)\big)[0.3,1)$ to ensure that the significant part of the Gaussian functions are not too concentrated, thus avoiding small inhomogeneities;
\item[(ii)]\textbf{[Nonsmooth \& discontinuous perturbations]} The weighted sums of absolute value of trigonometric functions combined with a mask, given by
\begin{equation}\label{eqn:dataset_nons}
\begin{array}{rl}
\gamma=&\chi_\Omega\sum_{i,j=1}^{J} \frac{c_{a,ij}}{i^3+j^3}|\cos\big(i(x+x_{a,i})\pi\big)\cos\big(j(y+y_{a,j})\pi\big)|\vspace{0.4cm}\\
\mbox{and}\quad \eta=&\chi_\Omega\sum_{i,j=1}^{J} \frac{c_{n,ij}}{i^3+j^3}|\cos\big(i(x+x_{n,i})\pi\big)\cos\big(j(y+y_{n,j})\pi\big)|,
\end{array}
\end{equation}
where the amplitudes $c_{a,ij}$ and $c_{n,ij}$ follow the standard normal distribution $N(0,1)$, and the relative phase shifts $x_{a,i}$, $x_{n,i}$, $y_{a,j}$ and $y_{a,j}$ are randomly selected from $[-1,1)$.
\end{enumerate}

We simulate $N= 1000$ samples of the perturbations $\gamma,\eta$ and corresponding scattering data $\Lambda^\omega$ of wave frequency $\omega\in\{2.5,5,10,20\}$ generated by solving the nonlinear problem \eqref{eqn:prob1}, where the radiation boundary conditions are implemented using the perfectly matched layer with order 4, using the Finite Element Method (implemented in NGsolve, available from https://github.com/NGSolve) with polynomials of order $4$ on a mesh of sizes $\frac{\pi}{4\sqrt{2}\omega}$ and $\frac{\pi}{4\omega}$ in $\Omega$ and in the air respectively.
The input scattering data $\Lambda^\omega$ are measured with $n_\theta\in\{32,64,128\}$ receivers and sources with the same equiangular directions, while the target perturbations $\gamma$ and $\eta$ are discretized in the computational domain $[-1,1)^2$ with a equispaced mesh of $n_c\times n_c$ points (i.e., $n_c\times n_c$ pixels) where $n_c\in\{64,128,256\}$ such that $n_c\geq n_\theta$.

We focus on the uncompressed neural networks $\psi\circ \psi_0\circ\phi^0$ (without butterfly structure) to observe the dependencies of the behavior of proposed neural networks on the sample sizes, wave frequencies and resolutions of the visualization for inputs and targets. 
We recall the construction of $\phi^0$ in Theorem \ref{thm:BF_0_app}:
\begin{align*}
\phi^0(\Lambda)
=\sum_{i=1}^2\begin{bmatrix}
(\phi_1^0 \Lambda_{\theta_1}^{\omega_i},\cdots,\phi_1^0 \Lambda_{\theta_{n_\theta}}^{\omega_i} )^T P_{\omega_i}\\
(\phi_2^0 \Lambda_{\theta_1}^{\omega_i},\cdots,\phi_2^0 \Lambda_{\theta_{n_\theta}}^{\omega_i})^T P_{\omega_i}
\end{bmatrix}
%=\sum_{i=1}^2\begin{bmatrix}
%\phi_1^0 \Lambda_{\theta_1}^{\omega_i}\\
%\vdots\\
%\phi_1^0 \Lambda_{\theta_{n_\theta}}^{\omega_i}\\
%\phi_2^0 \Lambda_{\theta_1}^{\omega_i}\\
%\vdots\\
%\phi_2^0 \Lambda_{\theta_{n_\theta}}^{\omega_i}
%\end{bmatrix}P_{\omega_i}
=\sum_{i=1}^2(
\phi_1^0 \Lambda_{\theta_1}^{\omega_i},
\cdots,
\phi_1^0 \Lambda_{\theta_{n_\theta}}^{\omega_i},
\phi_2^0 \Lambda_{\theta_1}^{\omega_i},
\cdots,
\phi_2^0 \Lambda_{\theta_{n_\theta}}^{\omega_i})^T P_{\omega_i}
\end{align*}
where
\begin{align*}
\phi_1^0\Lambda_{\theta_m}^\omega=\frac{4\pi^2 c_{\omega}}{n_\theta^2}b^T\big(\overline{K} \odot((C\odot\Lambda_{\theta_m}^{\omega}) K)\big)\quad \mbox{and}\quad
\phi_2^0\Lambda_{\theta_m}^\omega=\frac{4\pi^2 c_{\omega}}{n_\theta^2}b^T\big(\overline{K} \odot(\Lambda_{\theta_m}^{\omega} K)\big)
\end{align*}
with shifted inputs
\begin{align*}
\Lambda_{\theta_m}^{\omega_1}
=P_{1,\theta_m}P_u\Lambda P_{2,\theta_m}\quad \mbox{and}\quad  \Lambda_{\theta_m}^{\omega_2}=P_{1,\theta_m}P_l\Lambda P_{2,\theta_m}.
\end{align*}
To implement this process numerically and inspired by
\cite{ZhangZepeda-NunezLi:2024}, we decompose the complex-valued weight $K$ into two real-valued weights $K_{\rm cos}$ (the discretization of $\cos(\omega_1\rho\cos t)$) and $K_{\rm sin}$ (the discretization of $\sin(\omega_1\rho\cos t)$) such that $K=K_{\rm cos}-iK_{\rm sin}$. With the inputs $\Lambda^\omega_{\theta_m}=\Lambda^\omega_{\theta_m,{\rm Re}}+i\Lambda^\omega_{\theta_m,{\rm Im}}$ and Remark \ref{rem:real_out}, we rewrite $\phi_1^0$ and $\phi_2^0$ as 
\begin{align*}
\phi_1^0\Lambda_{\theta_m}^\omega=&\frac{4\pi^2 c_{\omega}}{n_\theta^2}b^T\big(K_{\rm cos} \odot((C\odot\Lambda^\omega_{\theta_m,{\rm Re}}) K_{\rm cos})+K_{\rm sin} \odot((C\odot\Lambda^\omega_{\theta_m,{\rm Re}}) K_{\rm sin})\\
&\qquad\quad\;\;+K_{\rm cos} \odot((C\odot\Lambda^\omega_{\theta_m,{\rm Im}}) K_{\rm sin})-K_{\rm sin} \odot((C\odot\Lambda^\omega_{\theta_m,{\rm Im}}) K_{\rm cos})\big),\\
\phi_2^0\Lambda_{\theta_m}^\omega=&\frac{4\pi^2 c_{\omega}}{n_\theta^2}b^T\big(K_{\rm cos} \odot(\Lambda^\omega_{\theta_m,{\rm Re}} K_{\rm cos})+K_{\rm sin} \odot(\Lambda^\omega_{\theta_m,{\rm Re}} K_{\rm sin})\\
&\qquad\quad\;\;+K_{\rm cos} \odot(\Lambda^\omega_{\theta_m,{\rm Im}} K_{\rm sin})-K_{\rm sin} \odot(\Lambda^\omega_{\theta_m,{\rm Im}} K_{\rm cos})\big).
\end{align*}
Then, based on the analysis in Section \ref{sec:app_CNN} (when $\alpha\to 0^+$) which suggests $\psi$ as a single layer CNN with $6$ filters followed by a deconvolutional neural network, we construct a $2$-dimensional deep CNN with $6$ filters of the same size for implementation. 
The proposed neural network is shown in the Fig. \ref{fig:BFNNCNN}.

\begin{figure}[hbt!]
\centering
\begin{tikzpicture}[x=1.5cm, y=1.5cm, >=stealth]

% Input node
\node [draw, rectangle, fill=white] (input) at (0,0) {$\Lambda=\begin{bmatrix}
\Lambda^{\omega_1}\\
\Lambda^{\omega_2}    
\end{bmatrix}$};

% Layer 1 nodes
\foreach \m [count=\y] in {1,2}
\node [draw, rectangle, fill=white] (layer1-1) at (1.6,2 -1) {$\{\Lambda^{\omega_1}_{\theta_m}\}_{m=1}^{n_\theta}$};
\node [draw, rectangle, fill=white] (layer1-2) at (1.6,1 -2) {$\{\Lambda^{\omega_2}_{\theta_m}\}_{m=1}^{n_\theta}$};

% Layer 2 nodes
\foreach \m [count=\y] in {1,2,3,4}
\node [draw, rectangle, fill=white] (layer2-1) at (3.3,2.5-1) {$
O_{1,\omega_1}^{0'}$};
\node [draw, rectangle, fill=white] (layer2-2) at (3.3,2.5-2) {$
O_{2,\omega_1}^{0'}$};
\node [draw, rectangle, fill=white] (layer2-3) at (3.3,2.5-3) {$
O_{1,\omega_2}^{0'}$};
\node [draw, rectangle, fill=white] (layer2-4) at (3.3,2.5-4) {$
O_{2,\omega_2}^{0'}$};

% Layer 3 nodes
\foreach \m [count=\y] in {1,2}
\node [draw, rectangle, fill=white] (layer3-1) at (5,2.5-1) {$O_{1,\omega_1}^0:=O_{1,\omega_1}^{0'} P_{\omega_1}$};
\node [draw, rectangle, fill=white] (layer3-2) at (5,2.5-2) {$O_{2,\omega_1}^0:=O_{2,\omega_1}^{0'} P_{\omega_1}$};
\node [draw, rectangle, fill=white] (layer3-3) at (5,2.5-3) {$O_{1,\omega_2}^0:=O_{1,\omega_2}^{0'} P_{\omega_2}$};
\node [draw, rectangle, fill=white] (layer3-4) at (5,2.5-4) {$O_{2,\omega_2}^0:=O_{2,\omega_2}^{0'} P_{\omega_2}$};

% Layer 4 nodes
\foreach \m [count=\y] in {1,2}
\node [draw, rectangle, fill=white] (layer4-1) at (7.3,1.6-1) {$\sum_{i=1}^2 O_{1,\omega_i}^0$};
\node [draw, rectangle, fill=white] (layer4-2) at (7.3,0.4-1) {$\sum_{i=1}^2 O_{2,\omega_i}^0$};

% Layer 5 nodes
\foreach \m [count=\y] in {1}
  \node [draw, rectangle, fill=white] (layer5-\m) at (8.8,0) {$\phi^0(\Lambda)$};

% Layer 6 nodes
\foreach \m [count=\y] in {1}
  \node [draw, rectangle, fill=white] (layer6-\m) at (8.8,-1.5) {$\psi_0\circ\phi^0(\Lambda)$};

% Layer 7 nodes
\foreach \m [count=\y] in {1}
  \node [draw, rectangle, fill=white] (layer7-\m) at (8.8,-2.6) {$h(\Lambda)=\psi\circ\psi_0\circ\phi^0(\Lambda)$};

% Connect layers
\foreach \i in {1,2}
  \draw [->] (input) -- (layer1-\i);

\foreach \i in {1,2}
  \foreach \j in {1}
    \draw [->] (layer1-\j) -- (layer2-\i);

\foreach \i in {3,4}
  \foreach \j in {2}
    \draw [->] (layer1-\j) -- (layer2-\i);

\foreach \i in {1}
  \foreach \j in {1}
    \draw [->] (layer2-\j) -- (layer3-\i);
    
\foreach \i in {2}
  \foreach \j in {2}
    \draw [->] (layer2-\j) -- (layer3-\i);

\foreach \i in {3}
  \foreach \j in {3}
    \draw [->] (layer2-\j) -- (layer3-\i);

\foreach \i in {4}
  \foreach \j in {4}
    \draw [->] (layer2-\j) -- (layer3-\i);

\foreach \i in {1}
  \foreach \j in {1}
    \draw [->] (layer3-\j) -- (layer4-\i);

\foreach \i in {1}
  \foreach \j in {3}
    \draw [->] (layer3-\j) -- (layer4-\i);

\foreach \i in {2}
  \foreach \j in {2}
    \draw [->] (layer3-\j) -- (layer4-\i);

\foreach \i in {2}
  \foreach \j in {4}
    \draw [->] (layer3-\j) -- (layer4-\i);

\foreach \i in {1,2}
  \draw [->] (layer4-\i) -- (layer5-1);

\foreach \i in {1}
  \foreach \j in {1}
    \draw [->] (layer5-\j) -- (layer6-\i);

\foreach \i in {1}
  \foreach \j in {1}
    \draw [->] (layer5-\j) -- (layer6-\i);

\foreach \i in {1}
  \foreach \j in {1}
    \draw [->] (layer6-\j) -- (layer7-\i);

% Annotate layers
\node [align=center, above] at (0,0.5) {Input};
\node [align=center, above] at (1.6,1.3) {Shifted inputs};
\node [align=center, above] at (1.6,-1.8) {Shifted inputs};
\node [align=center, above] at (9.7,-1.1) {Coordinate \\transformation};
\node [align=center, above] at (9.2,-2.2) {CNN};
\end{tikzpicture}

\caption{The combined uncompressed neural networks $h=\psi\circ\psi_0\circ \phi^0$ where the intermediate outputs $O_{j,\omega_1}^{0'}:=(\phi_j^0 \Lambda^{\omega_1}_{\theta_1},\cdots,\phi_j^0 \Lambda^{\omega_1}_{\theta_{n_\theta}})^T\in\mathbb{R}^{ n_\theta\times n_\rho}$ and $O_{j,\omega_2}^{0'}:=(\phi_j^0 \Lambda^{\omega_2}_{\theta_1},\cdots,\phi_j^0 \Lambda^{\omega_2}_{\theta_{n_\theta}})^T\in\mathbb{R}^{ n_\theta\times n_\rho}$ for $j=1,2$.
The input $\Lambda\in \mathbb{C}^{2n_\theta\times n_\theta}$ is the discretized scattering data at two distinct wave frequencies $\omega_2>\omega_1$. $0$-level BFNN $\phi^0$ is applied to the shifted inputs and derives the intermediate outputs $\phi^0(\Lambda)\in\mathbb{R}^{2 n_\theta\times \frac{\omega_1}{\omega_2}n_\rho}$. After converting the polar coordinates convert to the Cartesian coordinates by $\psi_0:\mathbb{R}^{2 n_\theta\times \frac{\omega_1}{\omega_2}n_\rho}\to \mathbb{R}^{2 n_c\times n_c}$, the $2$-dimensional deep ($8$-layer) CNN $\psi$ with $6$ filters of the same size is used to derive the output for minimizing the loss function $\ell$.%$\ell(h,\Lambda,(\gamma,\eta)^T):=\|h(\Lambda)-(\gamma,\eta)^T\|_F$ with the target $(\gamma,\eta)^T$.
\label{fig:BFNNCNN}}
\end{figure}

Let $\mathcal{T}_{\rm tr}$ be a training dataset of size $N_{\rm tr}=|\mathcal{T}_{\rm tr}|$ (taken to be $50$, $100$, $300$, and $900$)
and $\mathcal{T}_{\rm te}$ be the testing dataset of size $|\mathcal{T}_{\rm te}|=100$ sharing the same type of perturbations.
We train the neural networks $h$ on the data in the training set $\mathcal{T}_{\rm tr}$ with the loss function 
\begin{align}\label{eqn:loss}
\ell_{\rm tr}(h,\Lambda,(\gamma,\eta)^T):=\sum_{(\Lambda,(\gamma,\eta)^T)\in\mathcal{T}_{\rm tr}}\ell(h,\Lambda,(\gamma,\eta)^T)=\sum_{(\Lambda,(\gamma,\eta)^T)\in\mathcal{T}_{\rm tr}}\|h(\Lambda)-(\gamma,\eta)^T\|_F,
\end{align} 
where $\Lambda$ {denotes the inputs} and $(\gamma,\eta)^T$ are the targets in $\mathcal{T}_{\rm tr}$.
The optimization problem \eqref{eqn:loss} is minimized by the Adam algorithm \cite{KingmaBa:2015} of batch size $50$ or $100$ with a learning rate $l_0=0.005$ for at most $300$ steps.

\subsection{Approximation and generalization capabilities}\label{sec:num_capability}
\begin{figure}
    \centering
    \begin{minipage}{0.15\textwidth}
        \centering
        \textbf{${\color{white}|}\gamma_{\rm exact}{\color{white}|}$}
        \includegraphics[width=\linewidth]{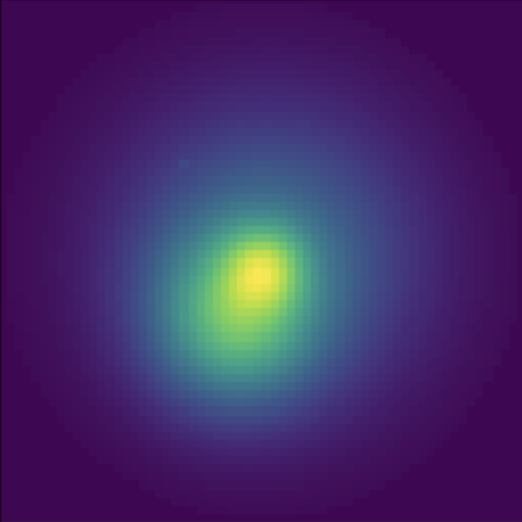}
    \end{minipage}\hfill
    \begin{minipage}{0.15\textwidth}
        \centering
        \textbf{${\color{white}|}\gamma_{\rm NN}{\color{white}|}$}
        \includegraphics[width=\linewidth]{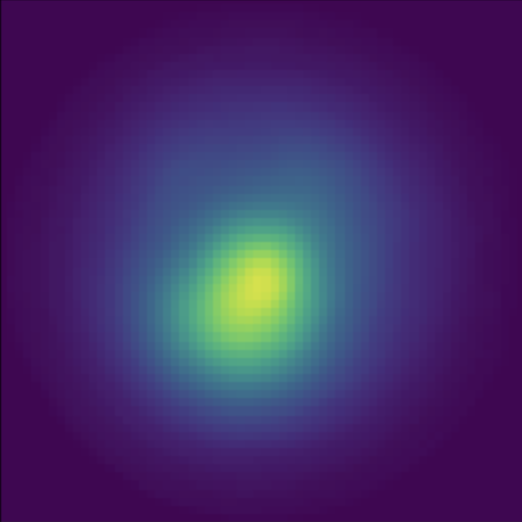}
    \end{minipage}\hfill
    \begin{minipage}{0.15\textwidth}
        \centering
        \textbf{$|\gamma_{\rm exact}-\gamma_{\rm NN}|$}
        \includegraphics[width=\linewidth]{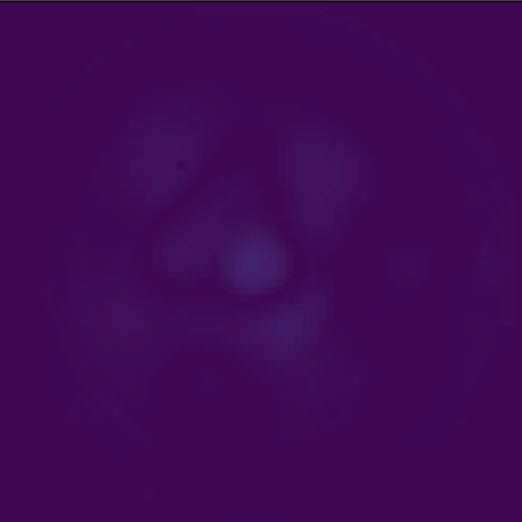}
    \end{minipage}\hfill
    \begin{minipage}{0.15\textwidth}
        \centering
        \textbf{${\color{white}|}\eta_{\rm exact}{\color{white}|}$}
        \includegraphics[width=\linewidth]{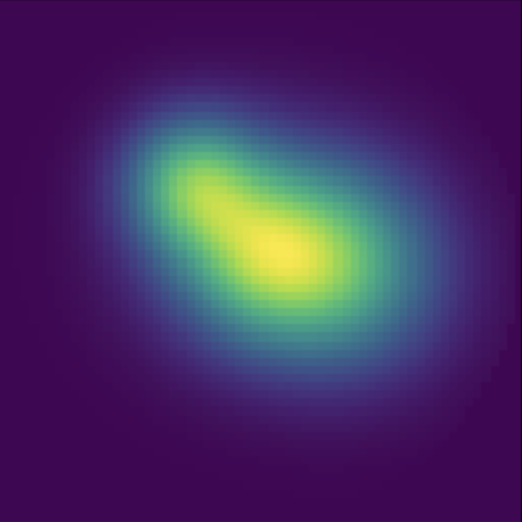}
    \end{minipage}\hfill
    \begin{minipage}{0.15\textwidth}
        \centering
        \textbf{${\color{white}|}\eta_{\rm NN}{\color{white}|}$}
        \includegraphics[width=\linewidth]{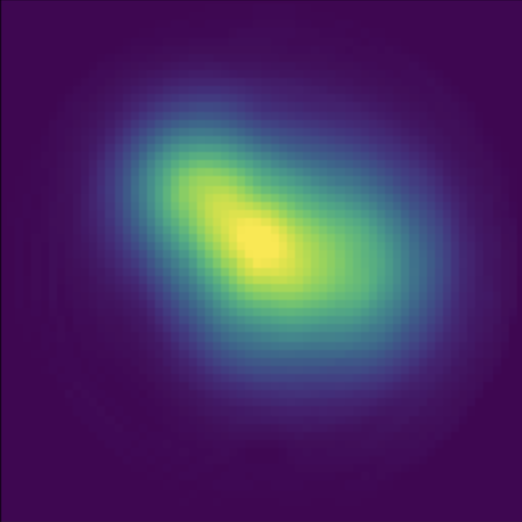}
    \end{minipage}\hfill
    \begin{minipage}{0.15\textwidth}
        \centering
        \textbf{$|\eta_{\rm exact}-\gamma_{\rm NN}|$}
        \includegraphics[width=\linewidth]{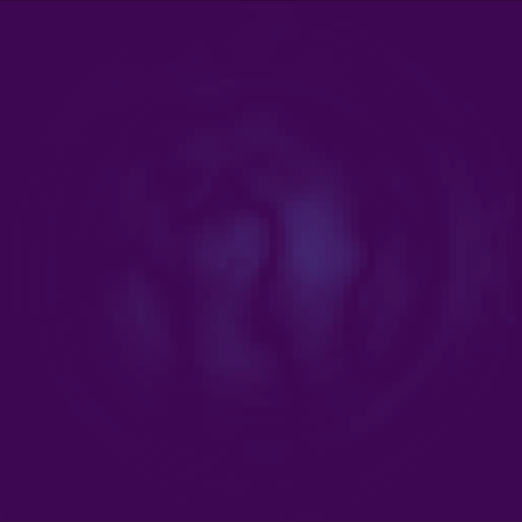}
    \end{minipage}
    
    \medskip

    \begin{minipage}{0.15\textwidth}
        \centering
        %\textbf{$\gamma_{\rm exact}$}
        \includegraphics[width=\linewidth]{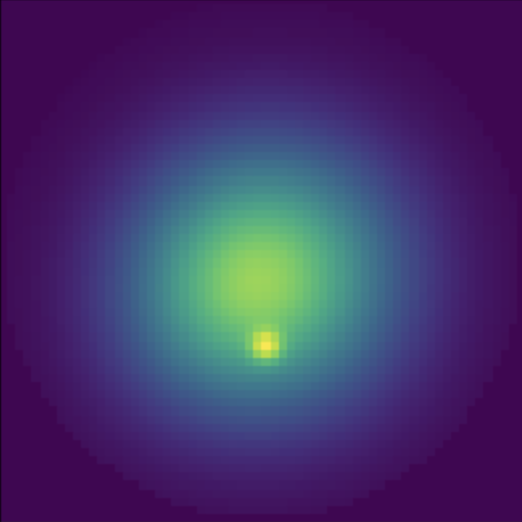}
    \end{minipage}\hfill
    \begin{minipage}{0.15\textwidth}
        \centering
        %\textbf{$\gamma_{\rm NN}$}
        \includegraphics[width=\linewidth]{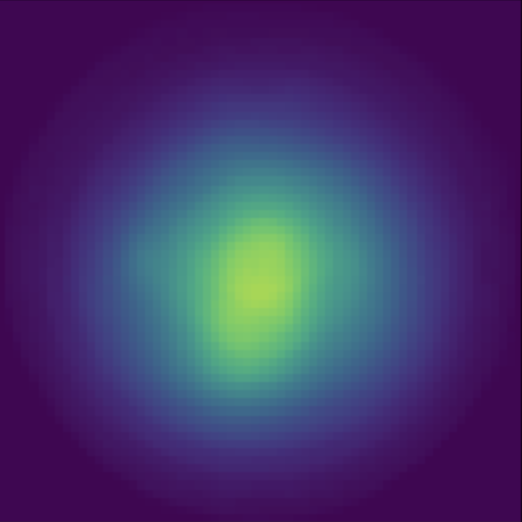}
    \end{minipage}\hfill
    \begin{minipage}{0.15\textwidth}
        \centering
        %\textbf{$|\gamma_{\rm exact}-\gamma_{\rm NN}|$}
        \includegraphics[width=\linewidth]{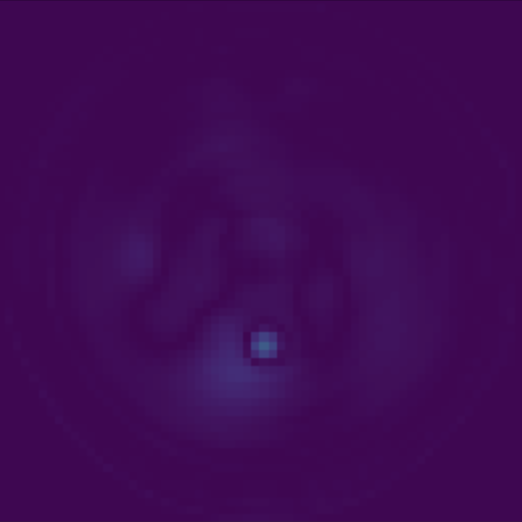}
    \end{minipage}\hfill
    \begin{minipage}{0.15\textwidth}
        \centering
        %\textbf{$\eta_{\rm exact}$}
        \includegraphics[width=\linewidth]{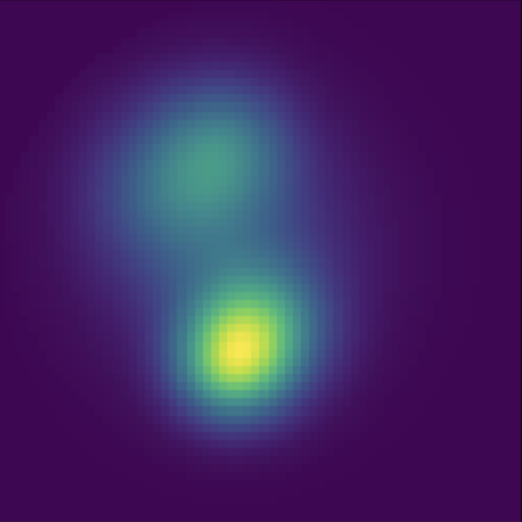}
    \end{minipage}\hfill
    \begin{minipage}{0.15\textwidth}
        \centering
        %\textbf{$\eta_{\rm NN}$}
        \includegraphics[width=\linewidth]{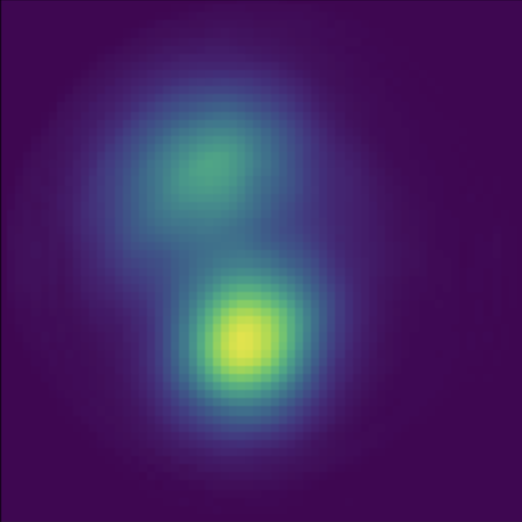}
    \end{minipage}\hfill
    \begin{minipage}{0.15\textwidth}
        \centering
        %\textbf{$|\eta_{\rm exact}-\gamma_{\rm NN}|$}
        \includegraphics[width=\linewidth]{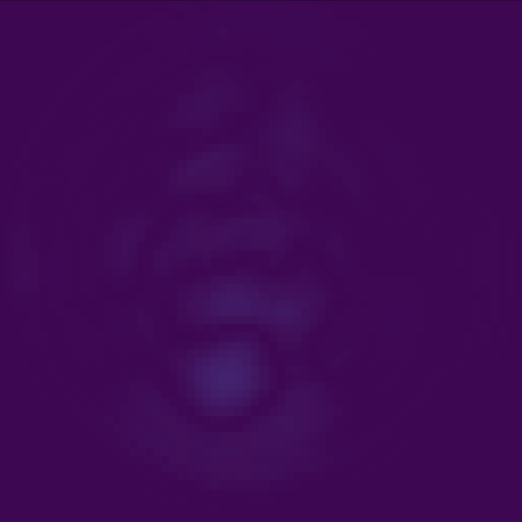}
    \end{minipage}

     \medskip

    \begin{minipage}{0.15\textwidth}
        \centering
        %\textbf{$\gamma_{\rm exact}$}
        \includegraphics[width=\linewidth]{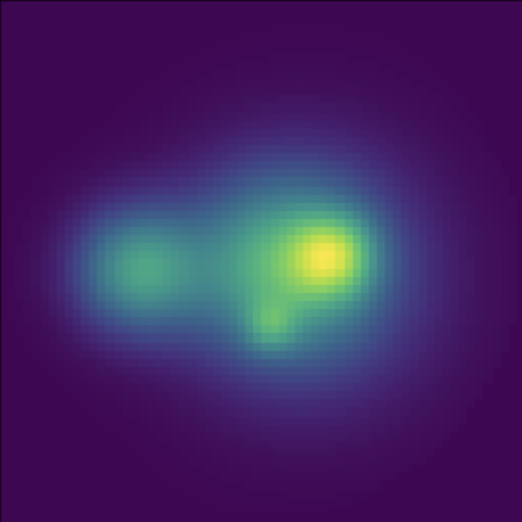}
    \end{minipage}\hfill
    \begin{minipage}{0.15\textwidth}
        \centering
        %\textbf{$\gamma_{\rm NN}$}
        \includegraphics[width=\linewidth]{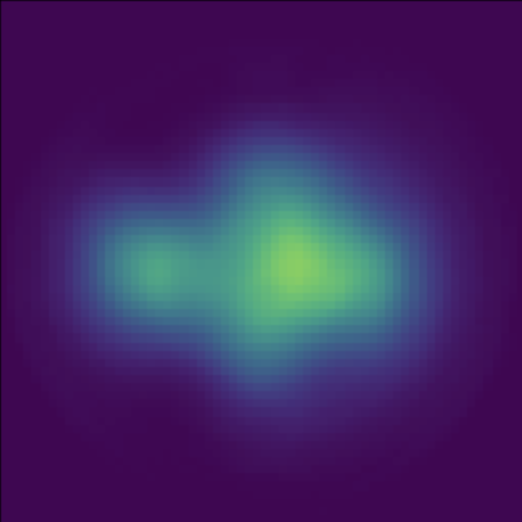}
    \end{minipage}\hfill
    \begin{minipage}{0.15\textwidth}
        \centering
        %\textbf{$|\gamma_{\rm exact}-\gamma_{\rm NN}|$}
        \includegraphics[width=\linewidth]{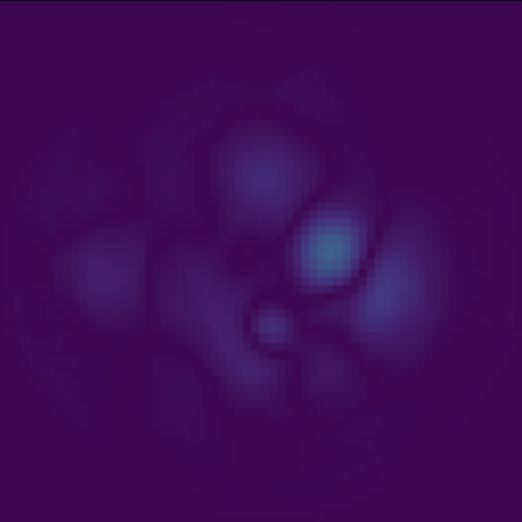}
    \end{minipage}\hfill
    \begin{minipage}{0.15\textwidth}
        \centering
        %\textbf{$\eta_{\rm exact}$}
        \includegraphics[width=\linewidth]{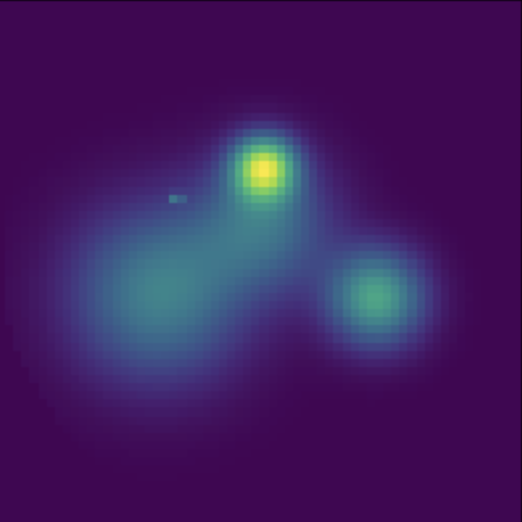}
    \end{minipage}\hfill
    \begin{minipage}{0.15\textwidth}
        \centering
        %\textbf{$\eta_{\rm NN}$}
        \includegraphics[width=\linewidth]{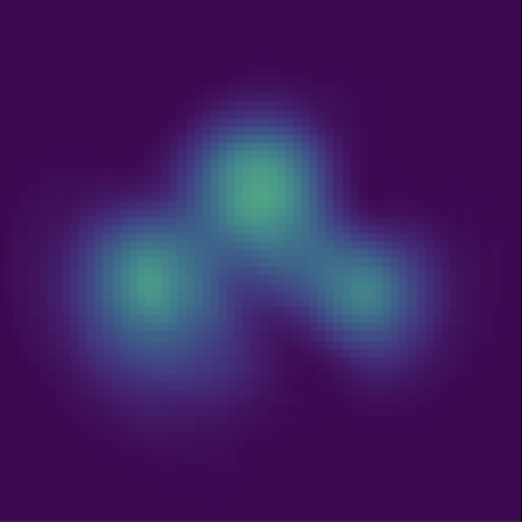}
    \end{minipage}\hfill
    \begin{minipage}{0.15\textwidth}
        \centering
        %\textbf{$|\eta_{\rm exact}-\gamma_{\rm NN}|$}
        \includegraphics[width=\linewidth]{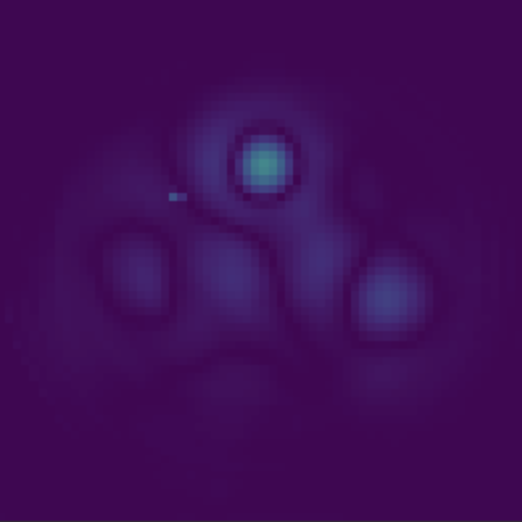}
    \end{minipage}

     \medskip

    \begin{minipage}{0.15\textwidth}
        \centering
        %\textbf{$\gamma_{\rm exact}$}
        \includegraphics[width=\linewidth]{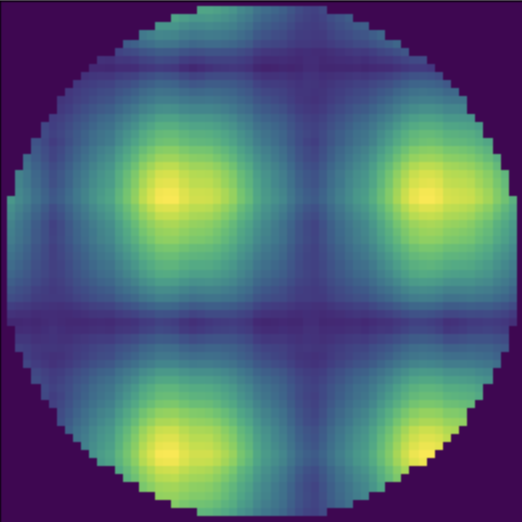}
    \end{minipage}\hfill
    \begin{minipage}{0.15\textwidth}
        \centering
        %\textbf{$\gamma_{\rm NN}$}
        \includegraphics[width=\linewidth]{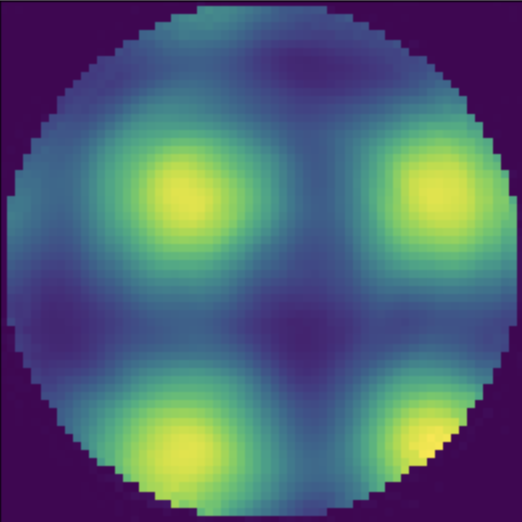}
    \end{minipage}\hfill
    \begin{minipage}{0.15\textwidth}
        \centering
        %\textbf{$|\gamma_{\rm exact}-\gamma_{\rm NN}|$}
        \includegraphics[width=\linewidth]{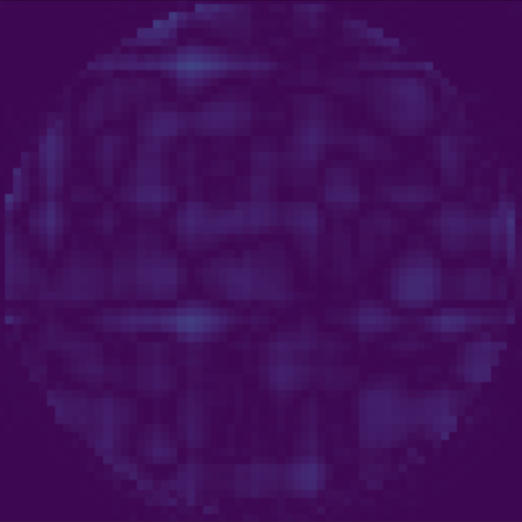}
    \end{minipage}\hfill
    \begin{minipage}{0.15\textwidth}
        \centering
        %\textbf{$\eta_{\rm exact}$}
        \includegraphics[width=\linewidth]{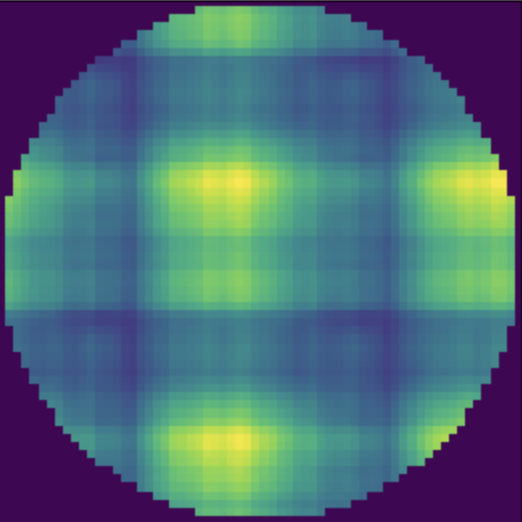}
    \end{minipage}\hfill
    \begin{minipage}{0.15\textwidth}
        \centering
        %\textbf{$\eta_{\rm NN}$}
        \includegraphics[width=\linewidth]{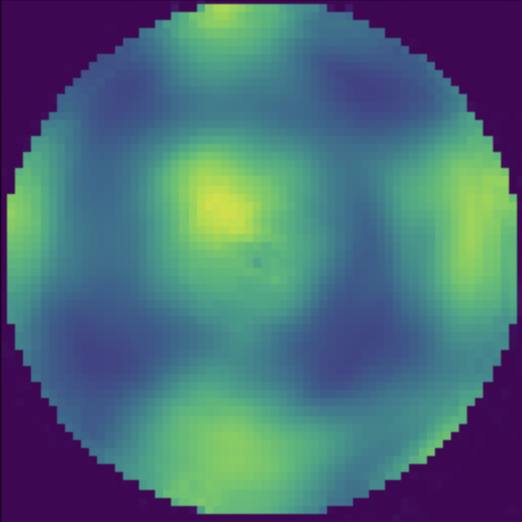}
    \end{minipage}\hfill
    \begin{minipage}{0.15\textwidth}
        \centering
        %\textbf{$|\eta_{\rm exact}-\gamma_{\rm NN}|$}
        \includegraphics[width=\linewidth]{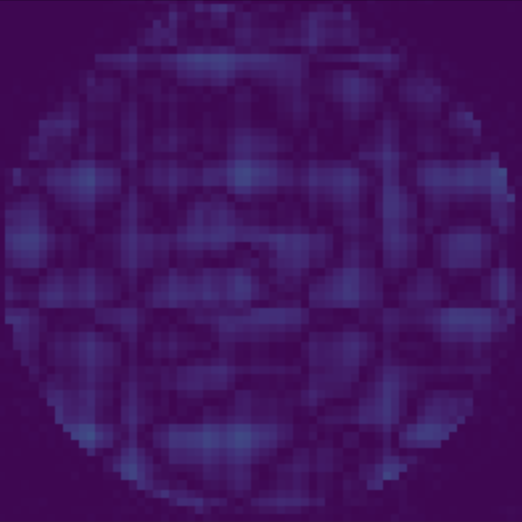}
    \end{minipage}
    
    \medskip

    \begin{minipage}{0.15\textwidth}
        \centering
        %\textbf{$\gamma_{\rm exact}$}
        \includegraphics[width=\linewidth]{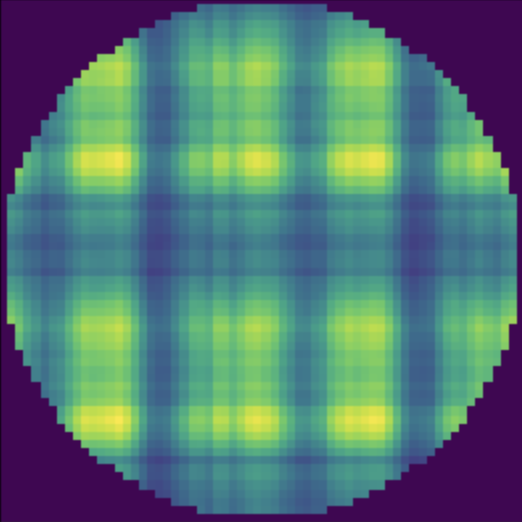}
    \end{minipage}\hfill
    \begin{minipage}{0.15\textwidth}
        \centering
        %\textbf{$\gamma_{\rm NN}$}
        \includegraphics[width=\linewidth]{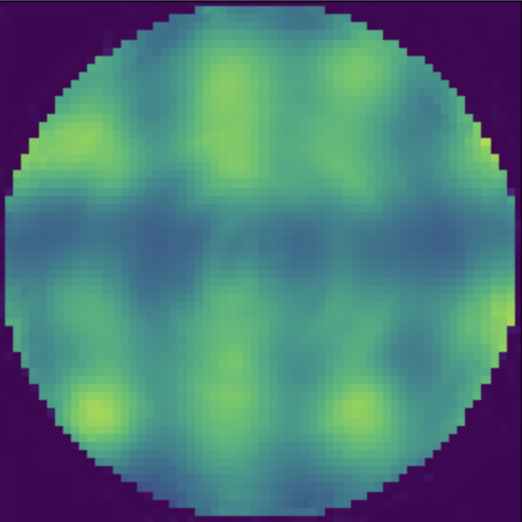}
    \end{minipage}\hfill
    \begin{minipage}{0.15\textwidth}
        \centering
        %\textbf{$|\gamma_{\rm exact}-\gamma_{\rm NN}|$}
        \includegraphics[width=\linewidth]{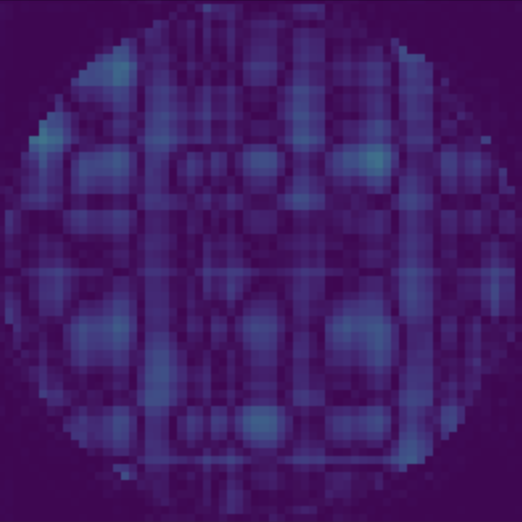}
    \end{minipage}\hfill
    \begin{minipage}{0.15\textwidth}
        \centering
        %\textbf{$\eta_{\rm exact}$}
        \includegraphics[width=\linewidth]{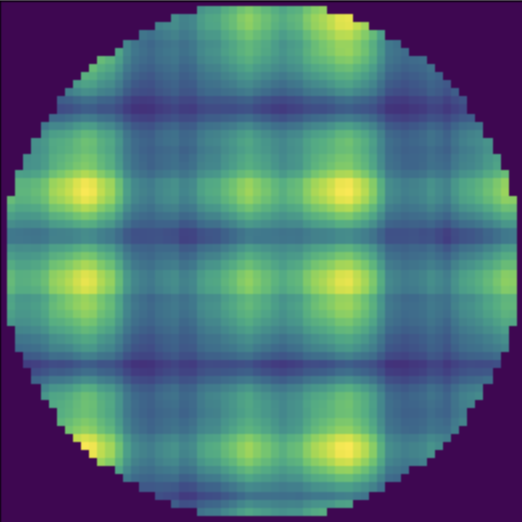}
    \end{minipage}\hfill
    \begin{minipage}{0.15\textwidth}
        \centering
        %\textbf{$\eta_{\rm NN}$}
        \includegraphics[width=\linewidth]{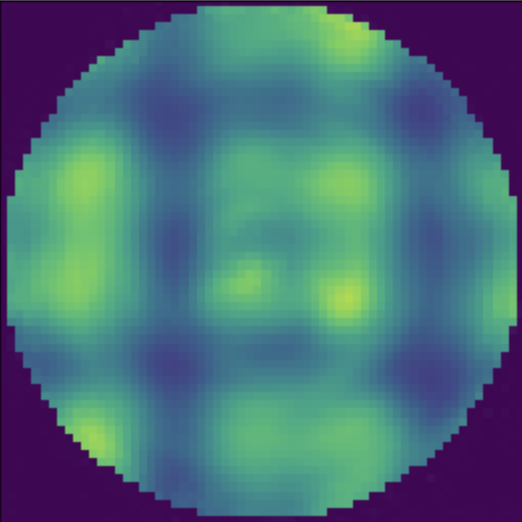}
    \end{minipage}\hfill
    \begin{minipage}{0.15\textwidth}
        \centering
        %\textbf{$|\eta_{\rm exact}-\gamma_{\rm NN}|$}
        \includegraphics[width=\linewidth]{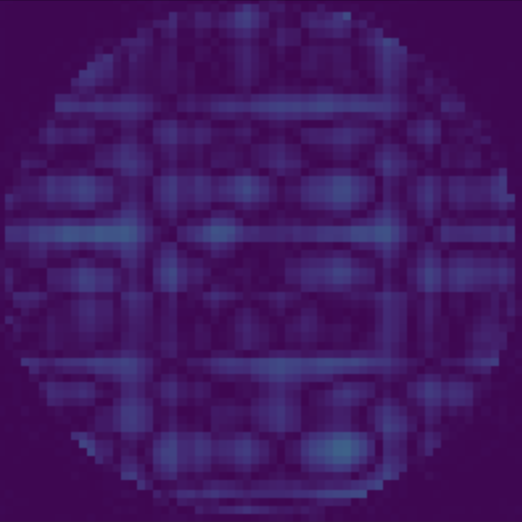}
    \end{minipage}
\caption{The exact and reconstructed perturbations from two different types of datasets. The first three rows are generated by \eqref{eqn:dataset_s} (smooth) and the last two rows are generated by \eqref{eqn:dataset_nons} (nonsmooth and discontinuous).
\label{fig:err}}   
\end{figure}

The numerical results of the proposed neural networks using wave frequencies $\omega_1=2.5$ and $\omega_2=5$ on both smooth (\eqref{eqn:dataset_s} with $J=5$) and nonsmooth (\eqref{eqn:dataset_nons} with $J=20$) perturbations training datasets, introduced in Section \ref{sec:train}, with $N_{\rm tr}=900$ and $n_\theta=n_c=64$ are given in Table \ref{tab:err}, where $e_{g}$ and $e_{a}$ denote respectively the relative generalization
error (on the testing dataset $\mathcal{T}_{\rm te}$) and the relative approximation error
(of $N_{\rm tr}$ training data points) of the neural networks {that is trained} on the training dataset $\mathcal{T}_{\rm tr}$. Several selected samples are presented in Fig. \ref{fig:err}, where $\gamma_{\rm exact},\eta_{\rm exact}$ are the exact targets and $\gamma_{\rm NN},\eta_{\rm NN}$ are the reconstructed perturbations from the trained neural networks.

\begin{table}[htp!]
  \centering
  \begin{threeparttable}
  \caption{Relative errors of the proposed neural networks on the smooth and nonsmooth perturbations training datasets where $N_{\rm tr}=900$, $n_\theta=n_c=64$, $\omega_1=2.5$ and $\omega_2=5$.\label{tab:err}}
    \begin{tabular}{cccc}
    \toprule
    \multicolumn{1}{c}{Training dataset}&
    \multicolumn{1}{c}{\;\;\;\;\;\;Smooth\;\;\;\;\;}&\multicolumn{1}{c}{\;\;\;Nonmooth\;\;\;}&\multicolumn{1}{c}{$50\%$Smooth + $50\%$ Nonmooth}\\
    \midrule
     $e_{g}$ & 12.57\% &  15.80\% & 15.18\%\\
    $e_{a}$ & 11.98\% &  15.39\% & 15.26\%\\
    \bottomrule
    \end{tabular}
    \end{threeparttable}
\end{table}

The numerical result indicates an acceptable approximation and generalization capability of the proposed neural networks on {both} smooth and nonsmooth perturbations datasets as well as the combination of these two types of datasets. 
It can be observed from the results in Fig. \ref{fig:err} that the proposed neural networks can accurately reconstruct the components of suitable frequencies. 
However, their ability to recover relatively concentrated and high-frequency components {is relatively poor}. 
In order to improve the neural networks for achieving higher accuracy of the reconstruction, we investigate the influence of the parameters of the neural networks on the approximation and generalization capabilities in the following sections.

\subsection{Dependence on the wave frequency}

\begin{figure}
    \centering
    $\gamma\;\;$\begin{minipage}{0.15\textwidth}
        \centering
        \textbf{{\color{white}(}Exact{\color{white}(}}
        %\textbf{$\gamma_{\rm exact}$}
        \includegraphics[width=\linewidth]{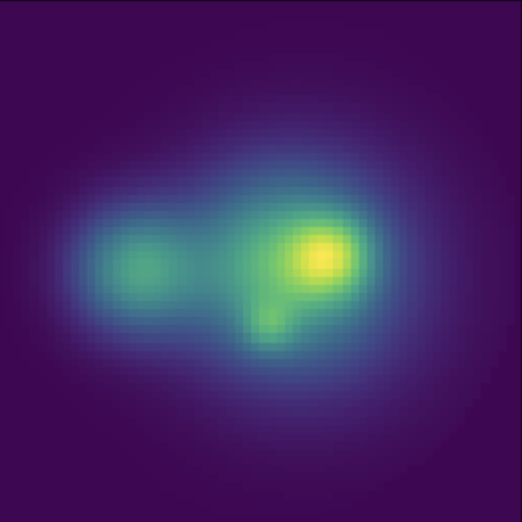}
    \end{minipage}\hfill
    \begin{minipage}{0.15\textwidth}
        \centering
        \textbf{NN(2.5\&5)}
        %\textbf{$\gamma_{\rm NN}(2.5\;\&\;5)$}
        \includegraphics[width=\linewidth]{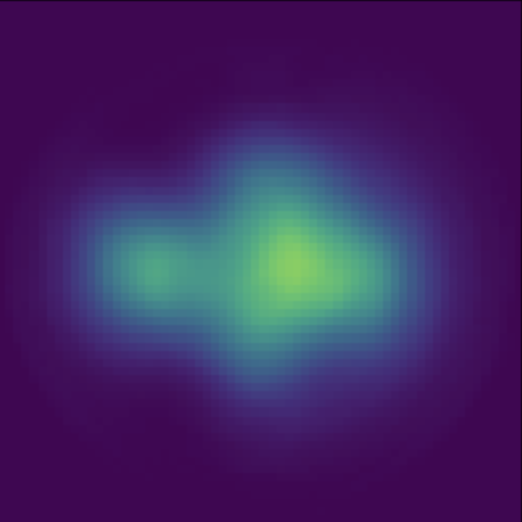}
    \end{minipage}\hfill
    \begin{minipage}{0.15\textwidth}
        \centering
        \textbf{NN(2.5\&10)}
        %\textbf{$\gamma_{\rm NN}(2.5\;\&\;10)$}
        \includegraphics[width=\linewidth]{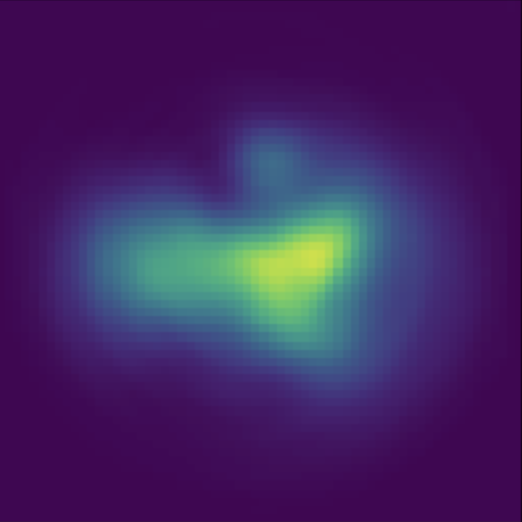}
    \end{minipage}\hfill
    \begin{minipage}{0.15\textwidth}
        \centering
        \textbf{NN(5\&10)}
        %\textbf{$\gamma_{\rm NN}(5\;\&\;10)$}
        \includegraphics[width=\linewidth]{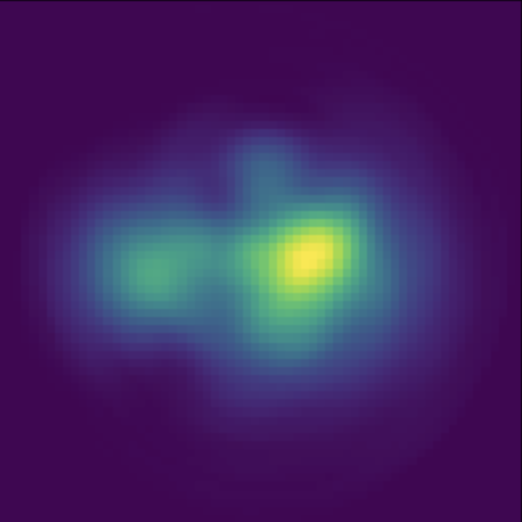}
    \end{minipage}\hfill
    \begin{minipage}{0.15\textwidth}
        \centering
        \textbf{NN(5\&20)}
        %\textbf{$\gamma_{\rm NN}(5\;\&\;10)$}
        \includegraphics[width=\linewidth]{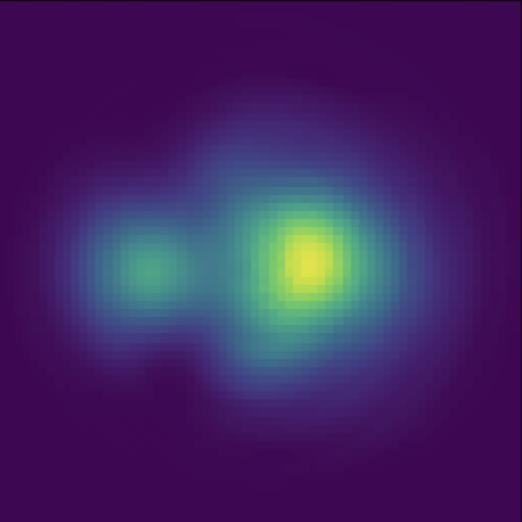}
    \end{minipage}\hfill
    \begin{minipage}{0.15\textwidth}
        \centering
        \textbf{NN(10\&20)}
        %\textbf{$\gamma_{\rm NN}(10\;\&\;20)$}
        \includegraphics[width=\linewidth]{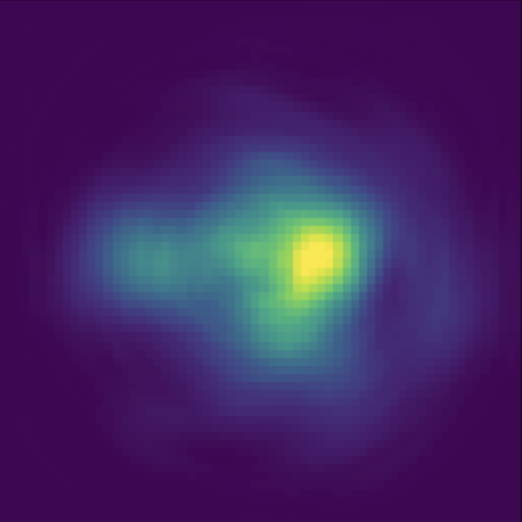}
    \end{minipage}

    $\eta\;\;$\begin{minipage}{0.15\textwidth}
        \centering
       % \textbf{$\eta_{\rm exact}$}
        \includegraphics[width=\linewidth]{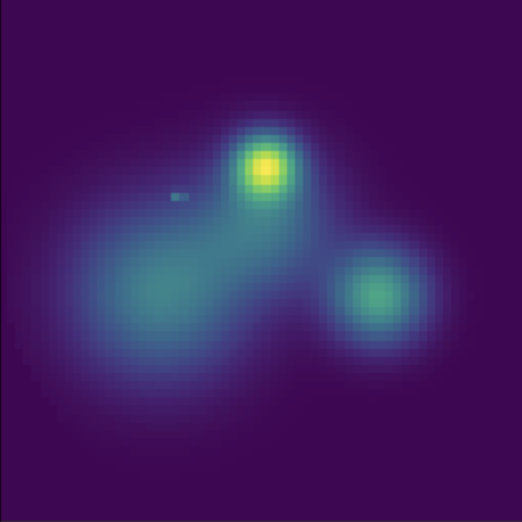}
    \end{minipage}\hfill
    \begin{minipage}{0.15\textwidth}
        \centering
       % \textbf{$\eta_{\rm NN}(2.5\;\&\;5)$}
        \includegraphics[width=\linewidth]{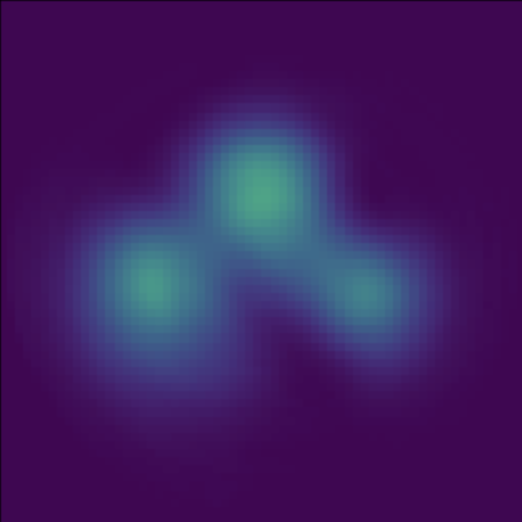}
    \end{minipage}\hfill
    \begin{minipage}{0.15\textwidth}
        \centering
        %\textbf{$\eta_{\rm NN}(2.5\;\&\;10)$}
        \includegraphics[width=\linewidth]{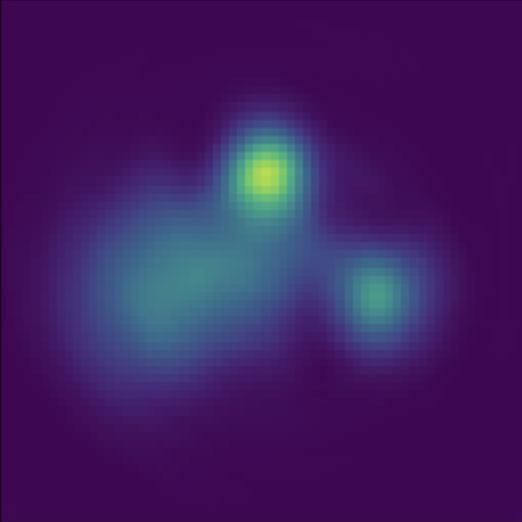}
    \end{minipage}\hfill
    \begin{minipage}{0.15\textwidth}
        \centering
        %\textbf{$\eta_{\rm NN}(5\;\&\;10)$}
        \includegraphics[width=\linewidth]{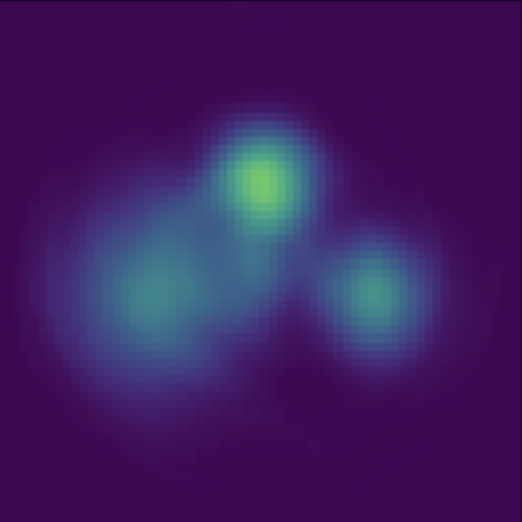}
    \end{minipage}\hfill
    \begin{minipage}{0.15\textwidth}
        \centering
        %\textbf{$\eta_{\rm NN}(5\;\&\;10)$}
        \includegraphics[width=\linewidth]{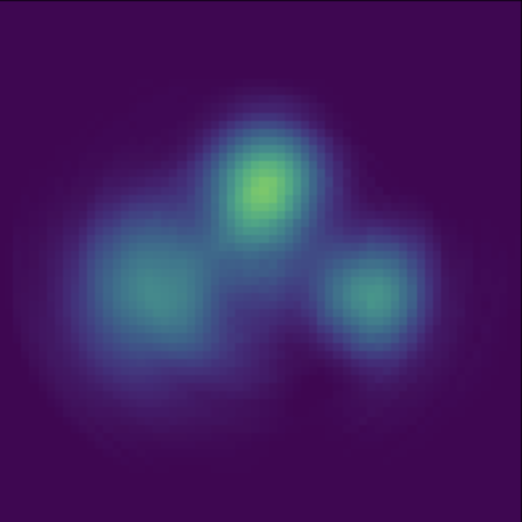}
    \end{minipage}\hfill
    \begin{minipage}{0.15\textwidth}
        \centering
        %\textbf{$\eta_{\rm NN}(10\;\&\;20)$}
        \includegraphics[width=\linewidth]{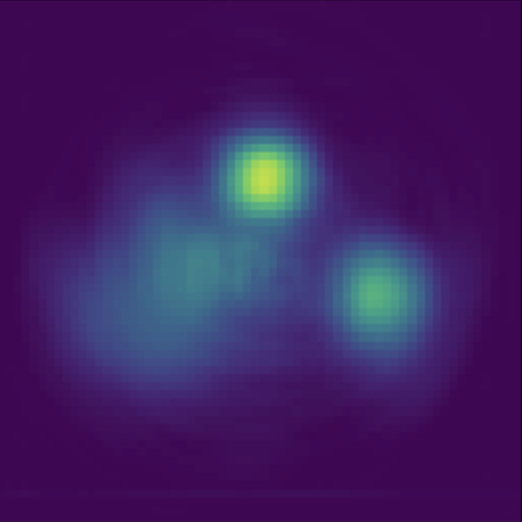}
    \end{minipage}
    
    \bigskip
    
    $\gamma\;\;$\begin{minipage}{0.15\textwidth}
        \centering
       % \textbf{$\eta_{\rm exact}$}
        \includegraphics[width=\linewidth]{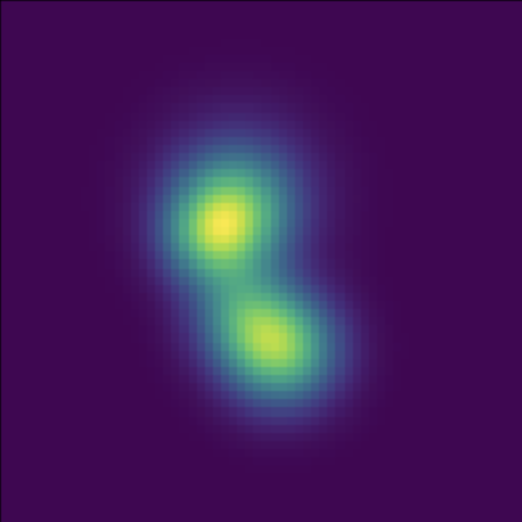}
    \end{minipage}\hfill
    \begin{minipage}{0.15\textwidth}
        \centering
       % \textbf{$\eta_{\rm NN}(2.5\;\&\;5)$}
        \includegraphics[width=\linewidth]{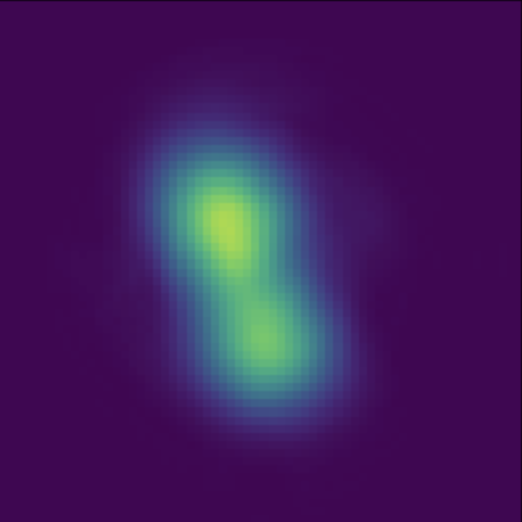}
    \end{minipage}\hfill
    \begin{minipage}{0.15\textwidth}
        \centering
        %\textbf{$\eta_{\rm NN}(2.5\;\&\;10)$}
        \includegraphics[width=\linewidth]{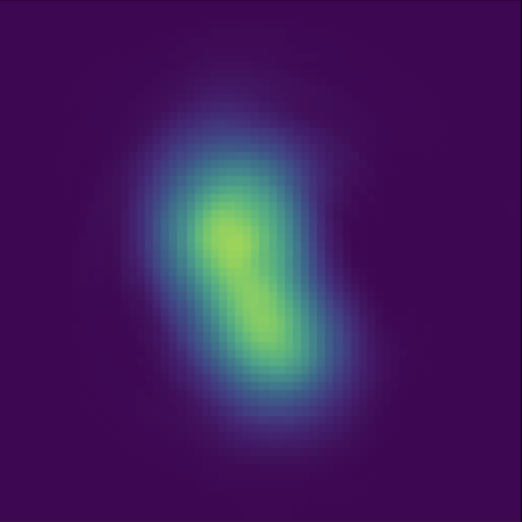}
    \end{minipage}\hfill
    \begin{minipage}{0.15\textwidth}
        \centering
        %\textbf{$\eta_{\rm NN}(5\;\&\;10)$}
        \includegraphics[width=\linewidth]{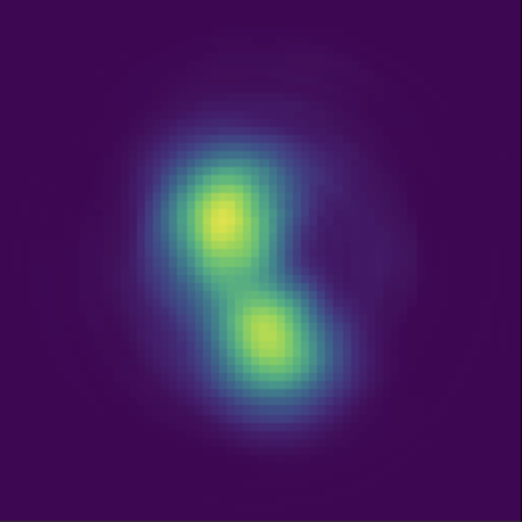}
    \end{minipage}\hfill
    \begin{minipage}{0.15\textwidth}
        \centering
        %\textbf{$\eta_{\rm NN}(5\;\&\;10)$}
        \includegraphics[width=\linewidth]{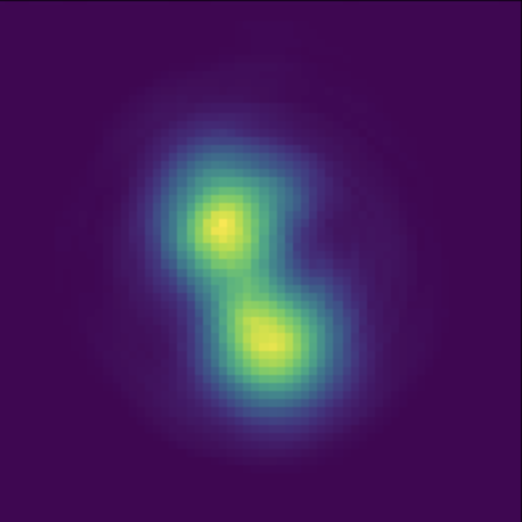}
    \end{minipage}\hfill
    \begin{minipage}{0.15\textwidth}
        \centering
        %\textbf{$\eta_{\rm NN}(10\;\&\;20)$}
        \includegraphics[width=\linewidth]{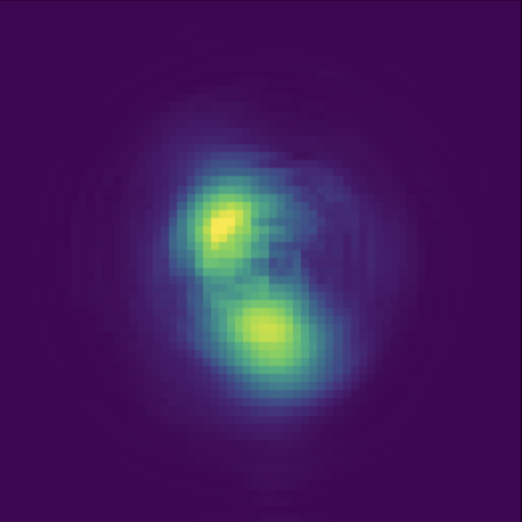}
    \end{minipage}

    $\eta\;\;$\begin{minipage}{0.15\textwidth}
        \centering
       % \textbf{$\eta_{\rm exact}$}
        \includegraphics[width=\linewidth]{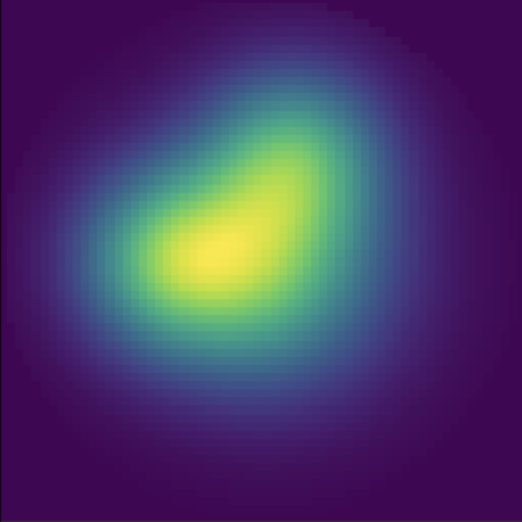}
    \end{minipage}\hfill
    \begin{minipage}{0.15\textwidth}
        \centering
       % \textbf{$\eta_{\rm NN}(2.5\;\&\;5)$}
        \includegraphics[width=\linewidth]{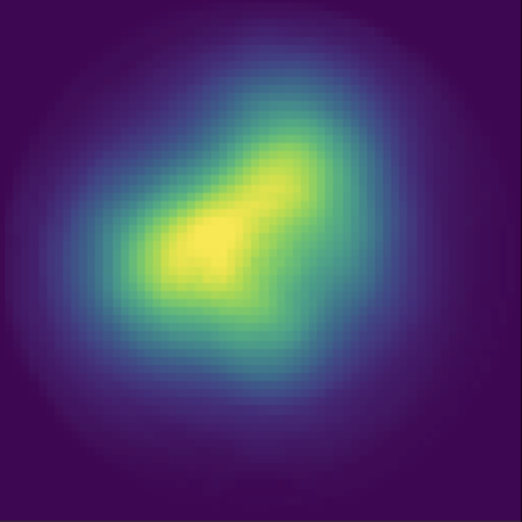}
    \end{minipage}\hfill
    \begin{minipage}{0.15\textwidth}
        \centering
        %\textbf{$\eta_{\rm NN}(2.5\;\&\;10)$}
        \includegraphics[width=\linewidth]{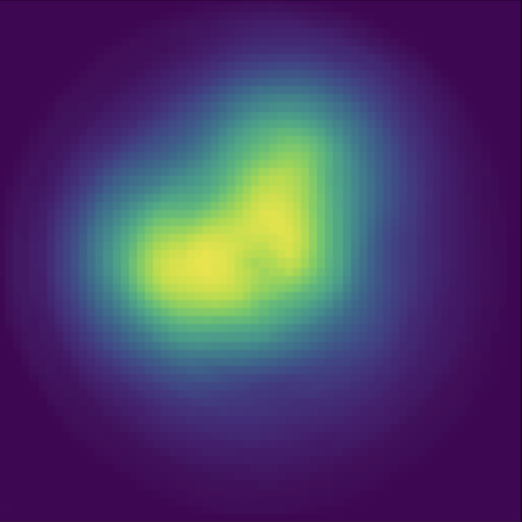}
    \end{minipage}\hfill
    \begin{minipage}{0.15\textwidth}
        \centering
        %\textbf{$\eta_{\rm NN}(5\;\&\;10)$}
        \includegraphics[width=\linewidth]{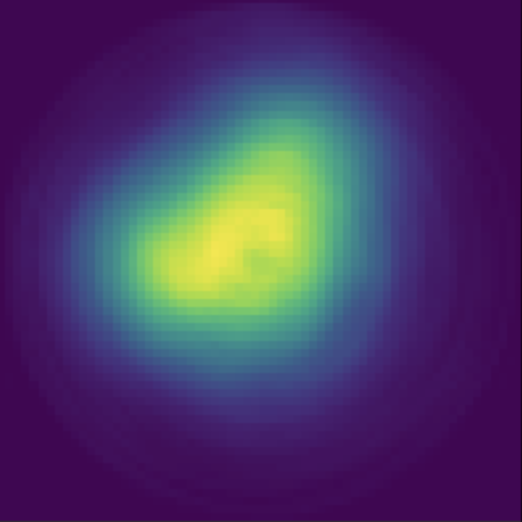}
    \end{minipage}\hfill
    \begin{minipage}{0.15\textwidth}
        \centering
        %\textbf{$\eta_{\rm NN}(5\;\&\;10)$}
        \includegraphics[width=\linewidth]{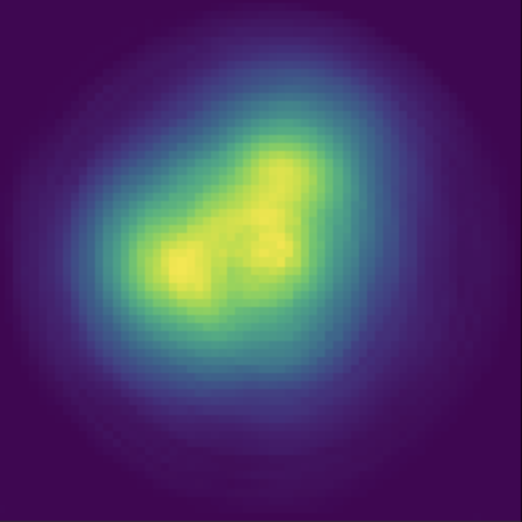}
    \end{minipage}\hfill
    \begin{minipage}{0.15\textwidth}
        \centering
        %\textbf{$\eta_{\rm NN}(10\;\&\;20)$}
        \includegraphics[width=\linewidth]{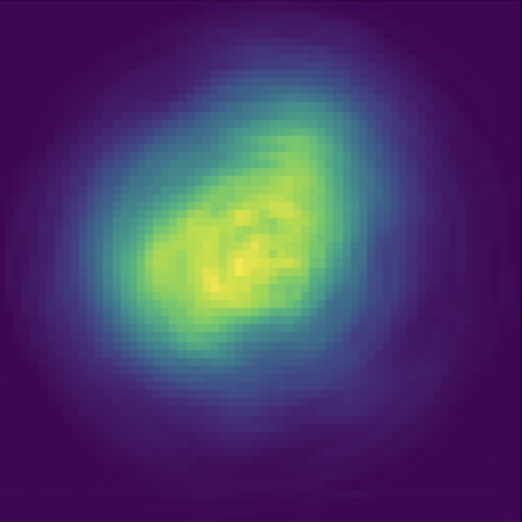}
    \end{minipage}
    
     \bigskip
    
    $\gamma\;\;$\begin{minipage}{0.15\textwidth}
        \centering
       % \textbf{$\eta_{\rm exact}$}
        \includegraphics[width=\linewidth]{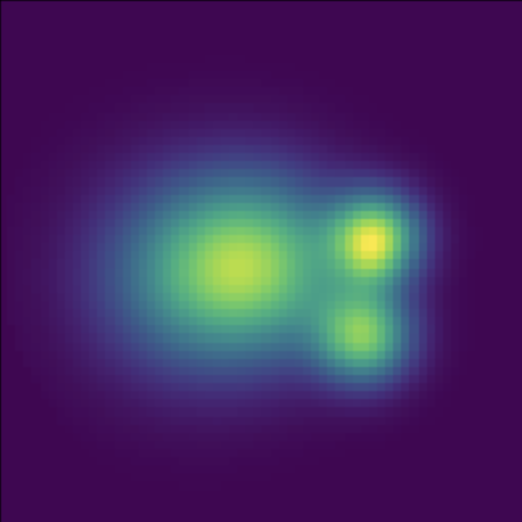}
    \end{minipage}\hfill
    \begin{minipage}{0.15\textwidth}
        \centering
       % \textbf{$\eta_{\rm NN}(2.5\;\&\;5)$}
        \includegraphics[width=\linewidth]{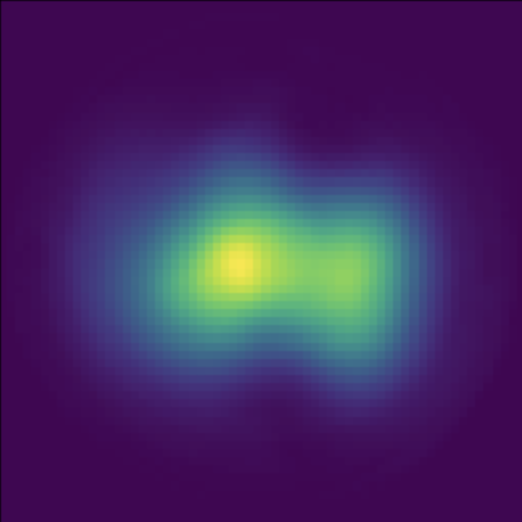}
    \end{minipage}\hfill
    \begin{minipage}{0.15\textwidth}
        \centering
        %\textbf{$\eta_{\rm NN}(2.5\;\&\;10)$}
        \includegraphics[width=\linewidth]{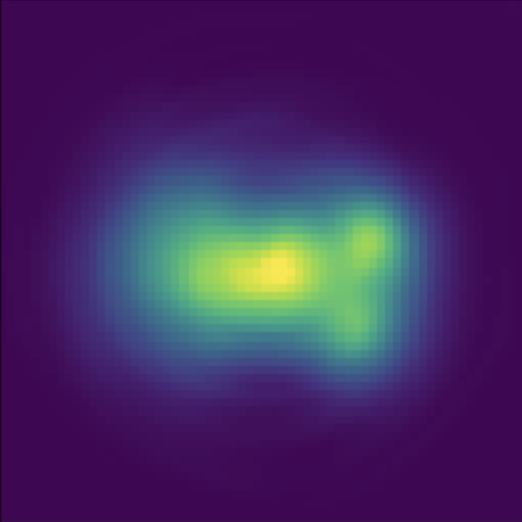}
    \end{minipage}\hfill
    \begin{minipage}{0.15\textwidth}
        \centering
        %\textbf{$\eta_{\rm NN}(5\;\&\;10)$}
        \includegraphics[width=\linewidth]{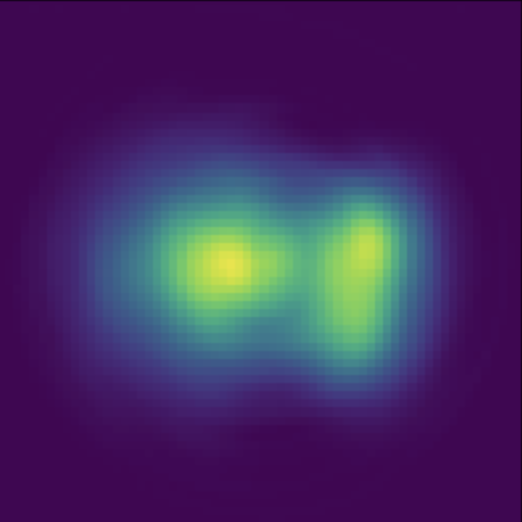}
    \end{minipage}\hfill
    \begin{minipage}{0.15\textwidth}
        \centering
        %\textbf{$\eta_{\rm NN}(5\;\&\;10)$}
        \includegraphics[width=\linewidth]{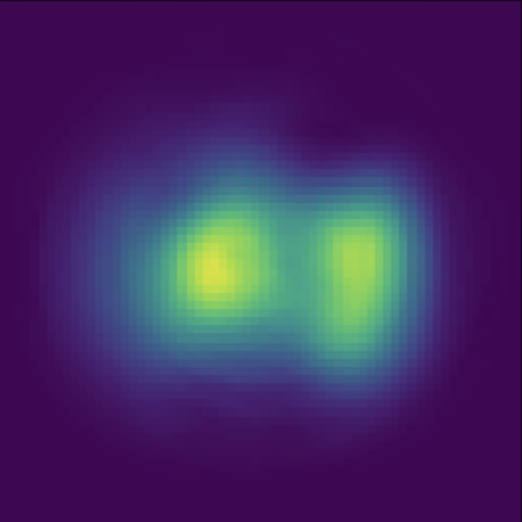}
    \end{minipage}\hfill
    \begin{minipage}{0.15\textwidth}
        \centering
        %\textbf{$\eta_{\rm NN}(10\;\&\;20)$}
        \includegraphics[width=\linewidth]{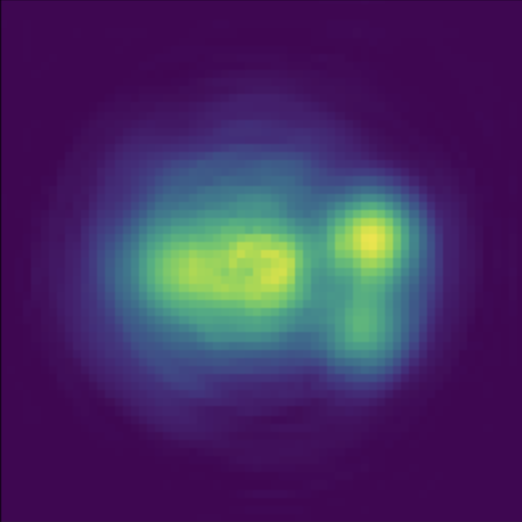}
    \end{minipage}

    $\eta\;\;$\begin{minipage}{0.15\textwidth}
        \centering
       % \textbf{$\eta_{\rm exact}$}
        \includegraphics[width=\linewidth]{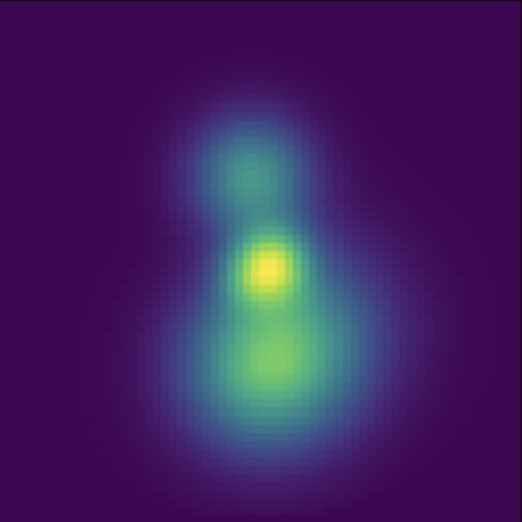}
    \end{minipage}\hfill
    \begin{minipage}{0.15\textwidth}
        \centering
       % \textbf{$\eta_{\rm NN}(2.5\;\&\;5)$}
        \includegraphics[width=\linewidth]{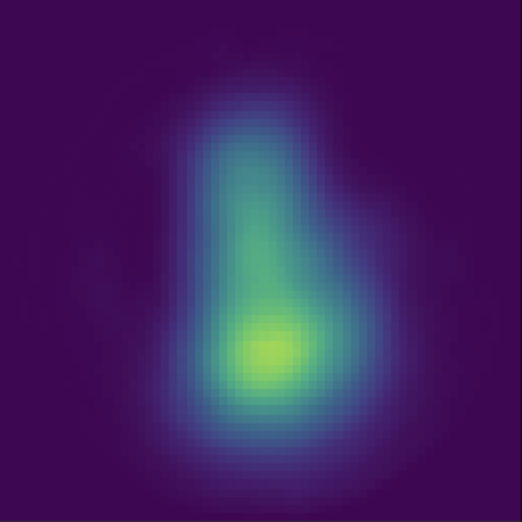}
    \end{minipage}\hfill
    \begin{minipage}{0.15\textwidth}
        \centering
        %\textbf{$\eta_{\rm NN}(2.5\;\&\;10)$}
        \includegraphics[width=\linewidth]{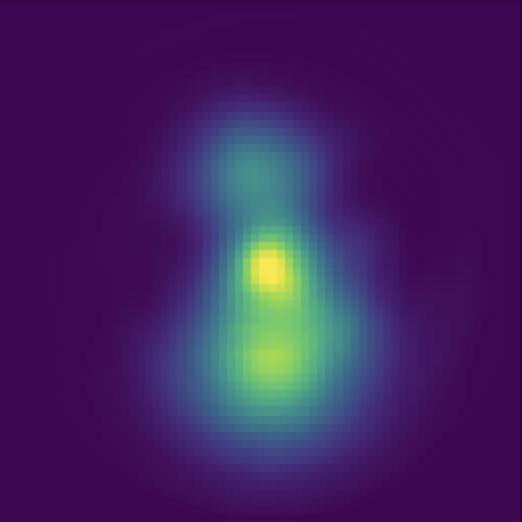}
    \end{minipage}\hfill
    \begin{minipage}{0.15\textwidth}
        \centering
        %\textbf{$\eta_{\rm NN}(5\;\&\;10)$}
        \includegraphics[width=\linewidth]{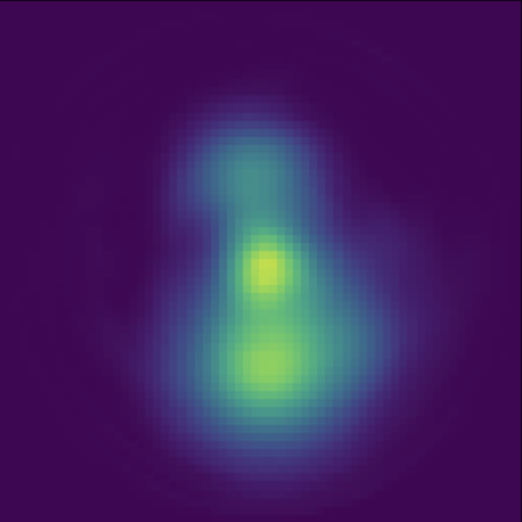}
    \end{minipage}\hfill
    \begin{minipage}{0.15\textwidth}
        \centering
        %\textbf{$\eta_{\rm NN}(5\;\&\;10)$}
        \includegraphics[width=\linewidth]{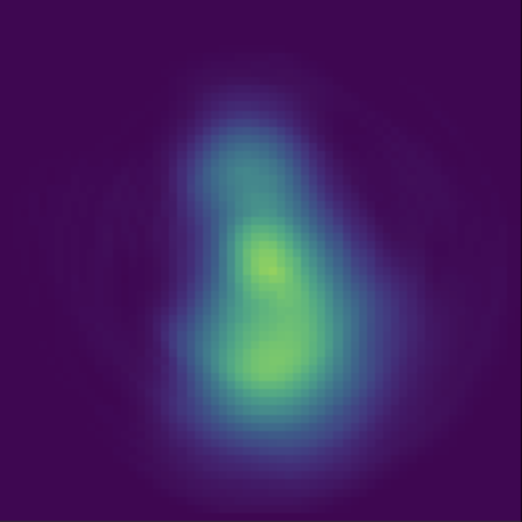}
    \end{minipage}\hfill
    \begin{minipage}{0.15\textwidth}
        \centering
        %\textbf{$\eta_{\rm NN}(10\;\&\;20)$}
        \includegraphics[width=\linewidth]{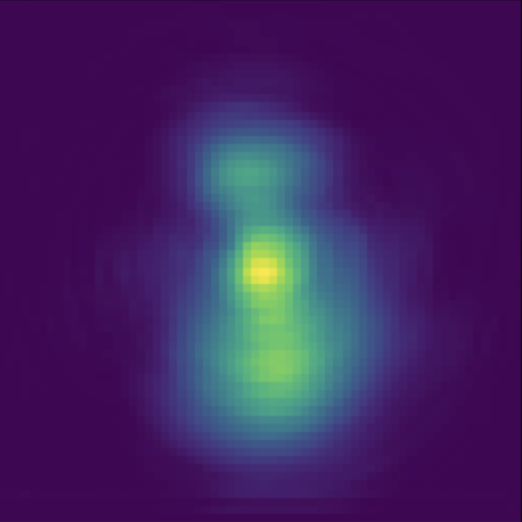}
    \end{minipage}
\caption{The exact and reconstructed perturbations from the proposed neural networks with different wave frequencies $\omega_1,\omega_2$. 
\label{fig:freq}}   
\end{figure}

\begin{figure}
    \centering
    $\gamma\;\;$\begin{minipage}{0.15\textwidth}
        \centering
        \textbf{{\color{white}(}Exact{\color{white}(} {\color{white}aaaaaaa}}
        %\textbf{$\gamma_{\rm exact}$}
        \includegraphics[width=\linewidth]{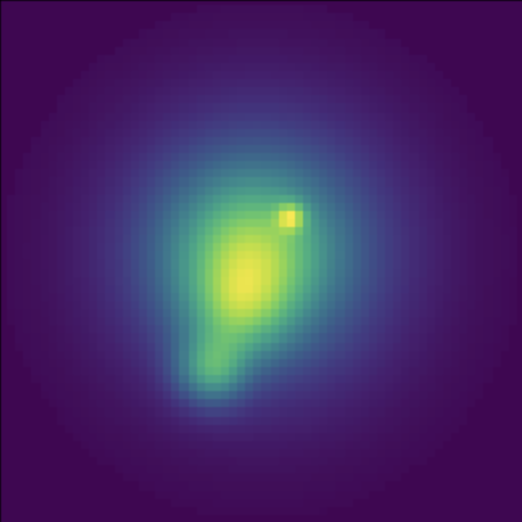}
    \end{minipage}\hfill
    \begin{minipage}{0.15\textwidth}
        \centering
        \textbf{NN(2.5\&5) error}
        %\textbf{$\gamma_{\rm NN}(2.5\;\&\;5)$}
        \includegraphics[width=\linewidth]{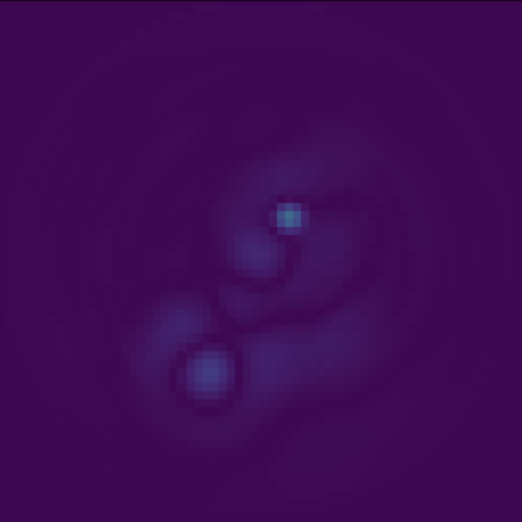}
    \end{minipage}\hfill
    \begin{minipage}{0.15\textwidth}
        \centering
        \textbf{NN(2.5\&10) error}
        %\textbf{$\gamma_{\rm NN}(2.5\;\&\;10)$}
        \includegraphics[width=\linewidth]{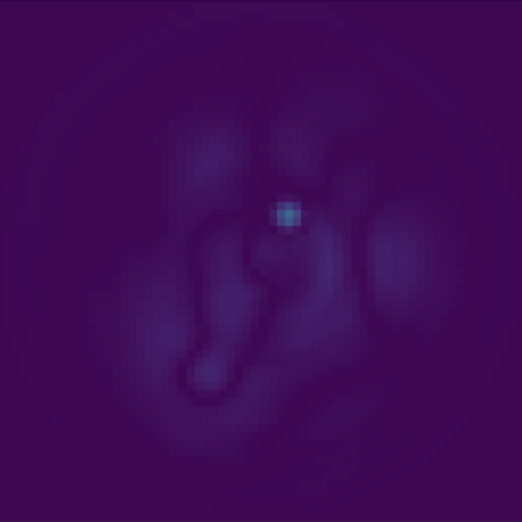}
    \end{minipage}\hfill
    \begin{minipage}{0.15\textwidth}
        \centering
        \textbf{NN(2.5\&20) error}
        %\textbf{$\gamma_{\rm NN}(5\;\&\;10)$}
        \includegraphics[width=\linewidth]{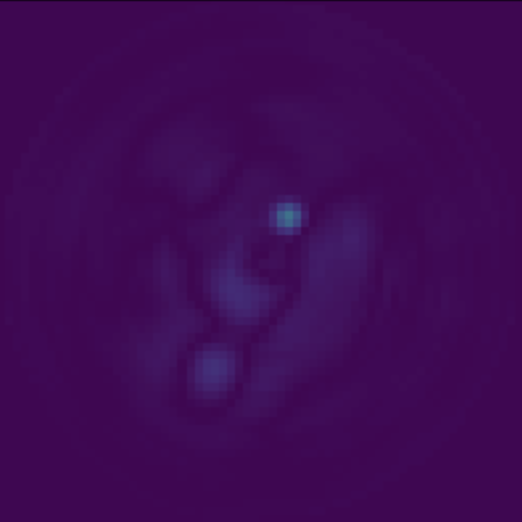}
    \end{minipage}\hfill
    \begin{minipage}{0.15\textwidth}
        \centering
        \textbf{NN(5\&10) error}
        %\textbf{$\gamma_{\rm NN}(5\;\&\;10)$}
        \includegraphics[width=\linewidth]{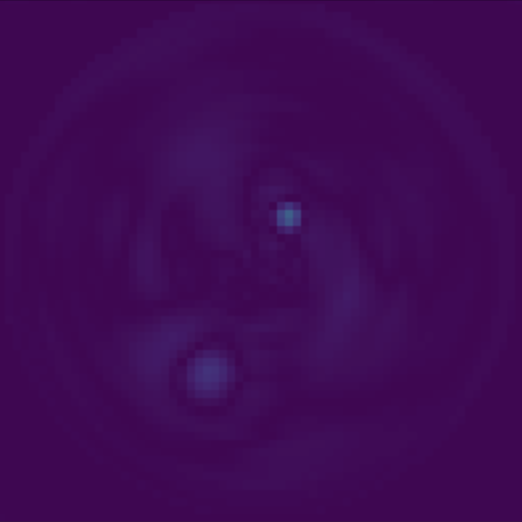}
    \end{minipage}\hfill
    \begin{minipage}{0.15\textwidth}
        \centering
        \textbf{NN(10\&20) error}
        %\textbf{$\gamma_{\rm NN}(10\;\&\;20)$}
        \includegraphics[width=\linewidth]{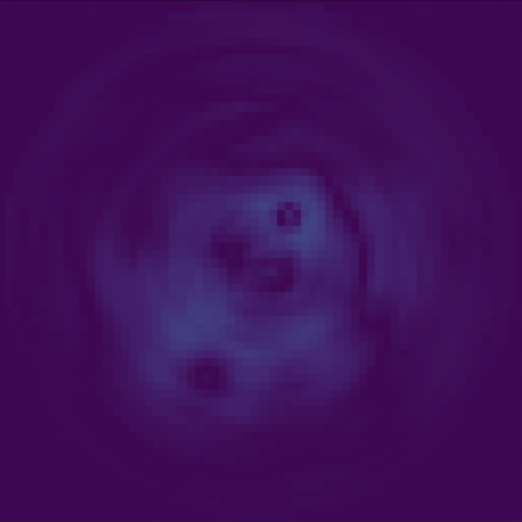}
    \end{minipage}

    $\eta\;\;$\begin{minipage}{0.15\textwidth}
        \centering
       % \textbf{$\eta_{\rm exact}$}
        \includegraphics[width=\linewidth]{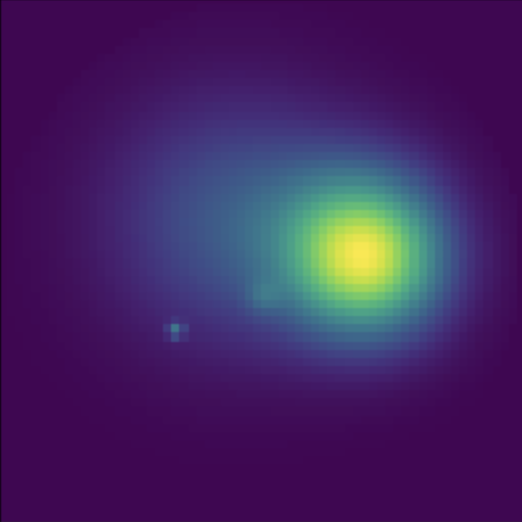}
    \end{minipage}\hfill
    \begin{minipage}{0.15\textwidth}
        \centering
       % \textbf{$\eta_{\rm NN}(2.5\;\&\;5)$}
        \includegraphics[width=\linewidth]{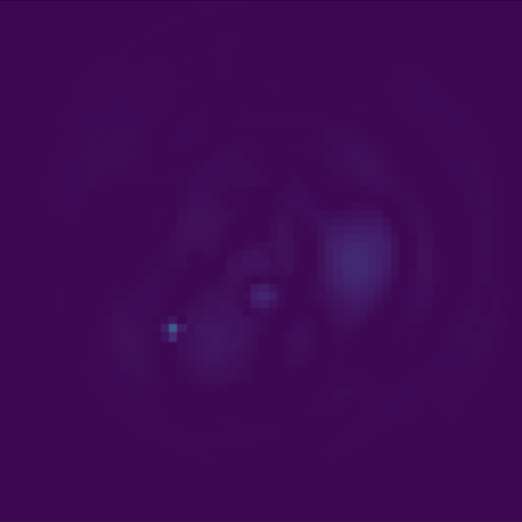}
    \end{minipage}\hfill
    \begin{minipage}{0.15\textwidth}
        \centering
        %\textbf{$\eta_{\rm NN}(2.5\;\&\;10)$}
        \includegraphics[width=\linewidth]{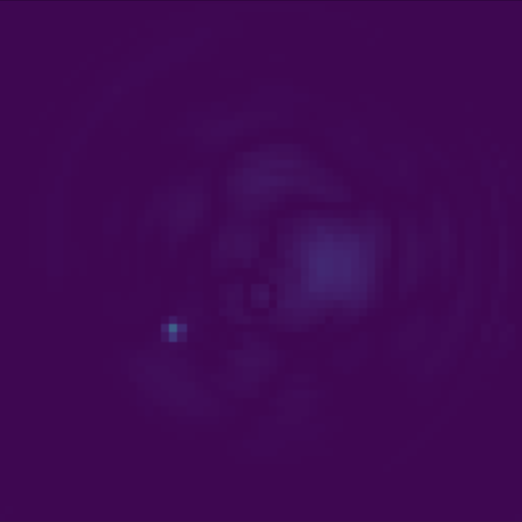}
    \end{minipage}\hfill
    \begin{minipage}{0.15\textwidth}
        \centering
        %\textbf{$\eta_{\rm NN}(5\;\&\;10)$}
        \includegraphics[width=\linewidth]{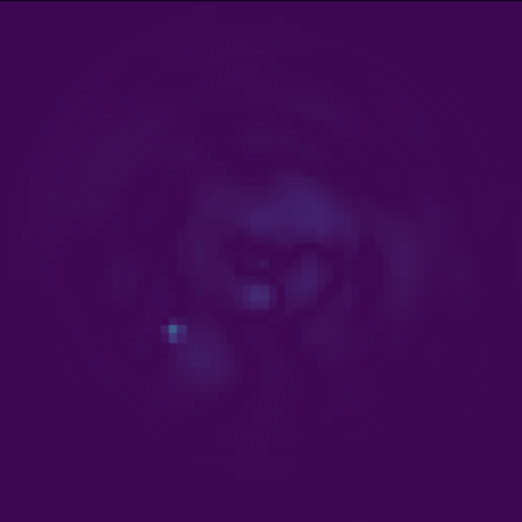}
    \end{minipage}\hfill
    \begin{minipage}{0.15\textwidth}
        \centering
        %\textbf{$\eta_{\rm NN}(5\;\&\;10)$}
        \includegraphics[width=\linewidth]{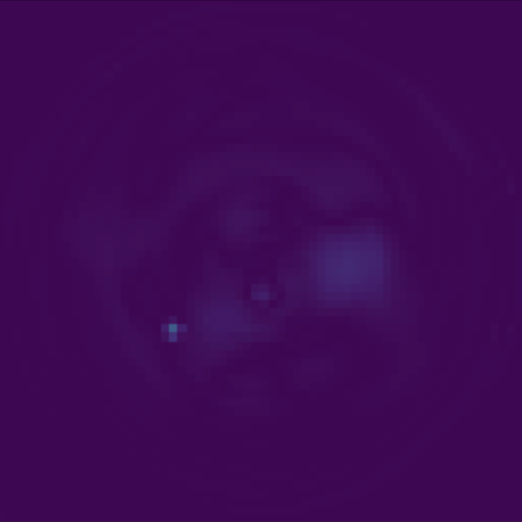}
    \end{minipage}\hfill
    \begin{minipage}{0.15\textwidth}
        \centering
        %\textbf{$\eta_{\rm NN}(10\;\&\;20)$}
        \includegraphics[width=\linewidth]{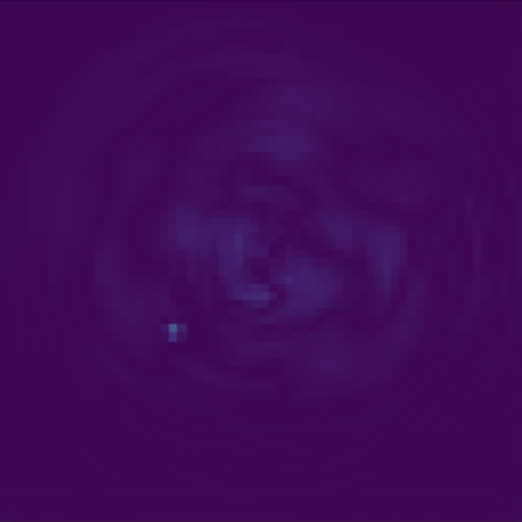}
    \end{minipage}
    
    \bigskip
    
    $\gamma\;\;$\begin{minipage}{0.15\textwidth}
        \centering
        %\textbf{$\gamma_{\rm exact}$}
        \includegraphics[width=\linewidth]{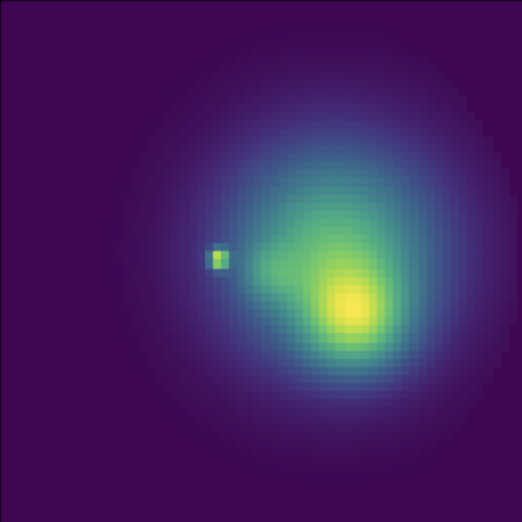}
    \end{minipage}\hfill
    \begin{minipage}{0.15\textwidth}
        \centering
        %\textbf{$\gamma_{\rm NN}(2.5\;\&\;5)$}
        \includegraphics[width=\linewidth]{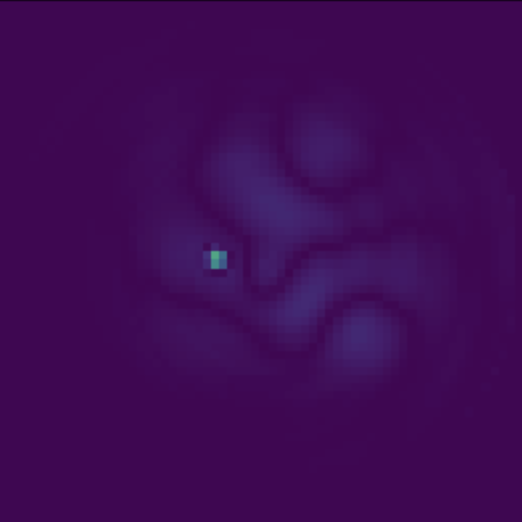}
    \end{minipage}\hfill
    \begin{minipage}{0.15\textwidth}
        \centering
        %\textbf{$\gamma_{\rm NN}(2.5\;\&\;10)$}
        \includegraphics[width=\linewidth]{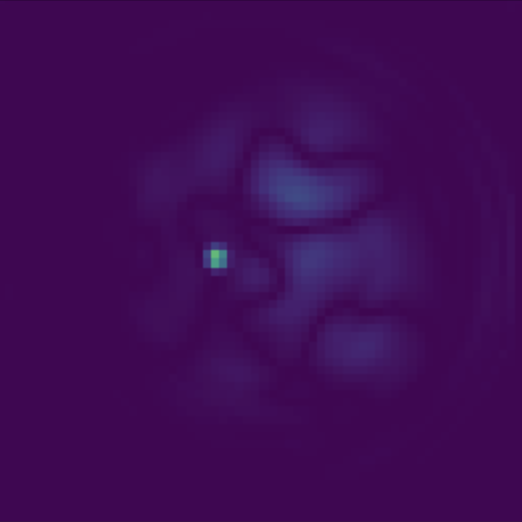}
    \end{minipage}\hfill
    \begin{minipage}{0.15\textwidth}
        \centering
        %\textbf{$\gamma_{\rm NN}(5\;\&\;10)$}
        \includegraphics[width=\linewidth]{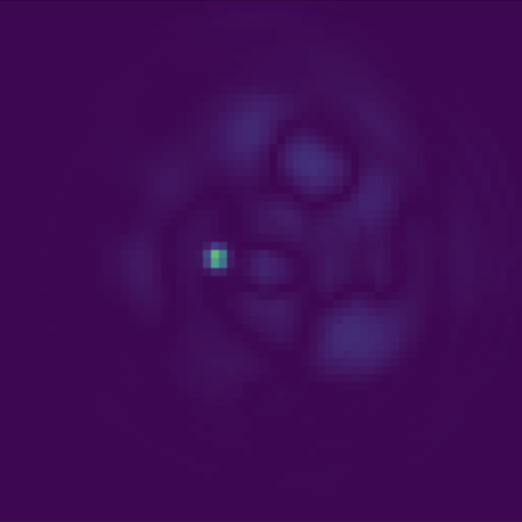}
    \end{minipage}\hfill
    \begin{minipage}{0.15\textwidth}
        \centering
        %\textbf{$\gamma_{\rm NN}(5\;\&\;10)$}
        \includegraphics[width=\linewidth]{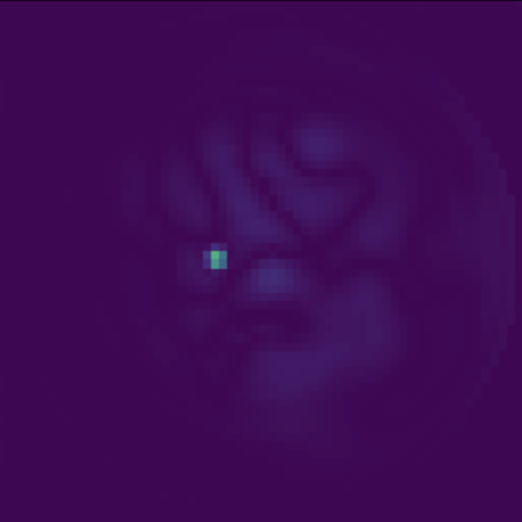}
    \end{minipage}\hfill
    \begin{minipage}{0.15\textwidth}
        \centering
        %\textbf{$\gamma_{\rm NN}(10\;\&\;20)$}
        \includegraphics[width=\linewidth]{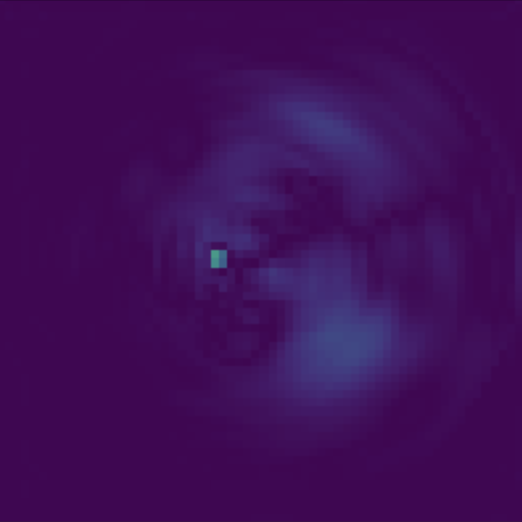}
    \end{minipage}

    $\eta\;\;$\begin{minipage}{0.15\textwidth}
        \centering
        %\textbf{$\eta_{\rm exact}$}
        \includegraphics[width=\linewidth]{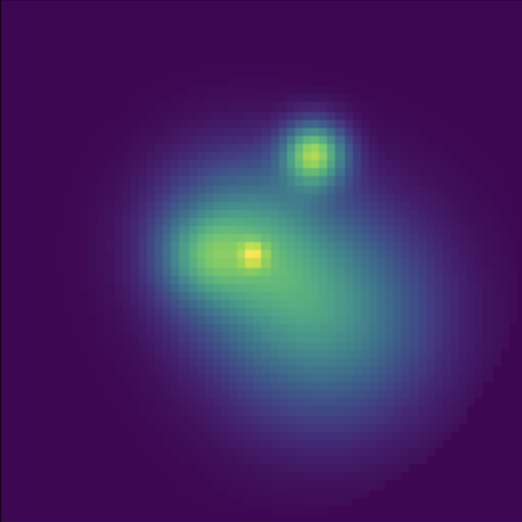}
    \end{minipage}\hfill
    \begin{minipage}{0.15\textwidth}
        \centering
        %\textbf{$\eta_{\rm NN}(2.5\;\&\;5)$}
        \includegraphics[width=\linewidth]{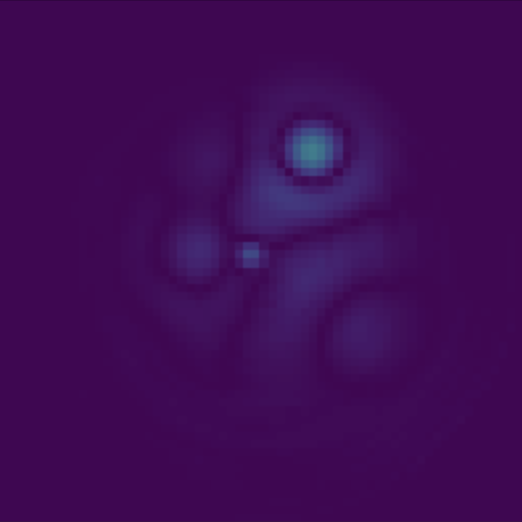}
    \end{minipage}\hfill
    \begin{minipage}{0.15\textwidth}
        \centering
        %\textbf{$\eta_{\rm NN}(2.5\;\&\;10)$}
        \includegraphics[width=\linewidth]{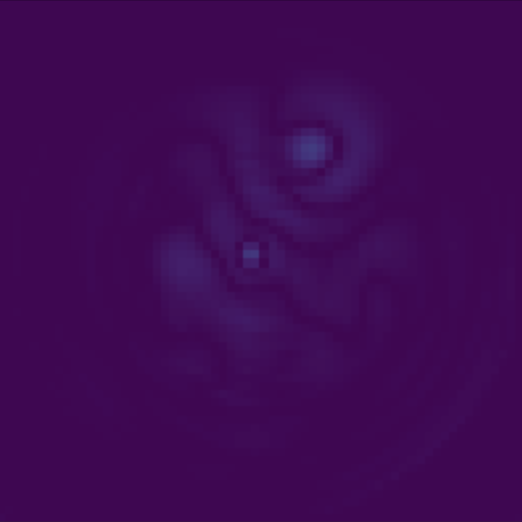}
    \end{minipage}\hfill
    \begin{minipage}{0.15\textwidth}
        \centering
        %\textbf{$\eta_{\rm NN}(5\;\&\;10)$}
        \includegraphics[width=\linewidth]{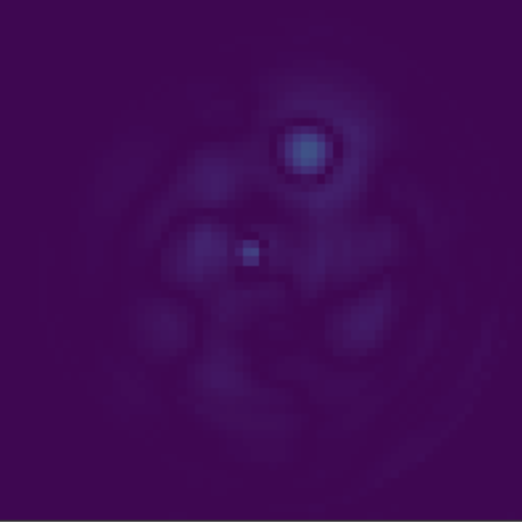}
    \end{minipage}\hfill
    \begin{minipage}{0.15\textwidth}
        \centering
        %\textbf{$\eta_{\rm NN}(5\;\&\;10)$}
        \includegraphics[width=\linewidth]{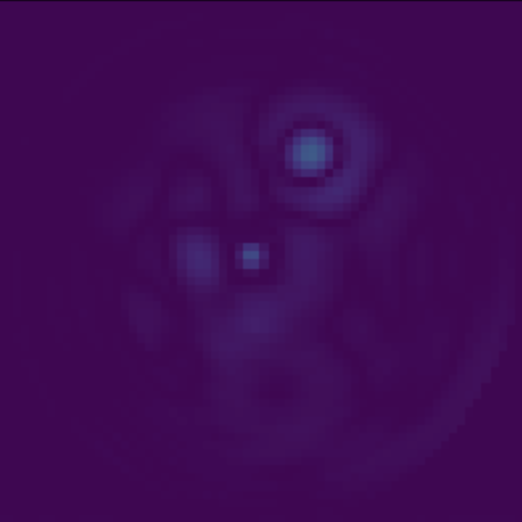}
    \end{minipage}\hfill
    \begin{minipage}{0.15\textwidth}
        \centering
        %\textbf{$\eta_{\rm NN}(10\;\&\;20)$}
        \includegraphics[width=\linewidth]{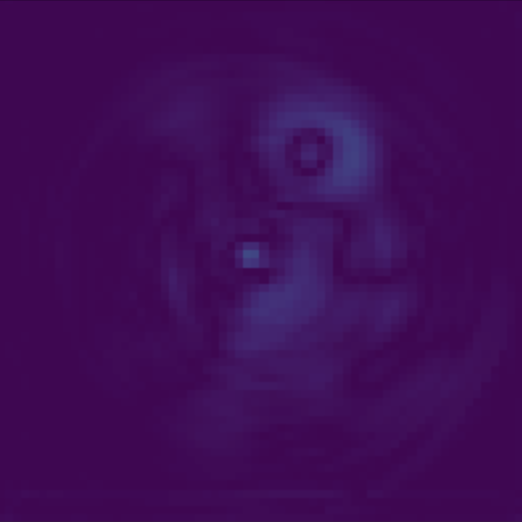}
    \end{minipage}
\caption{The exact perturbations and the reconstructed errors from the proposed neural networks with different wave frequencies $\omega_1,\omega_2$. 
\label{fig:smallinhom}}   
\end{figure}

Based on the observation in Fig. \ref{fig:err}, the recovery capability of the proposed neural networks is expected to be closely related to the frequencies of the perturbations $\gamma,\eta$ and the measurement $\Lambda^\omega$ used in the training process. We fix the size of training dataset and the resolution, the results for the networks with different wave frequencies $\omega_1,\omega_2$ are shown in Table \ref{tab:freq} and {Figs. \ref{fig:freq} and \ref{fig:smallinhom}}.  

These results indicate that a suitable wave {frequency pair} for the training does improve the accuracy of the reconstructions. In particular, higher frequencies allow the neural networks to capture more local details and perturbations with smaller support (high-frequency components). {However, they may reduce} its capability to capture the general features (low-frequency components). {Thus}, with a prior knowledge of the exact perturbations, a proper wave frequency pair can largely improve the capability of the proposed neural networks.

\begin{table}[htp!]
  \centering
  \begin{threeparttable}
  \caption{Relative errors of the proposed neural networks with different wave frequencies on the smooth perturbations training datasets where $N_{\rm tr}=900$ and $n_\theta=n_c=64$.\label{tab:freq}}
    \begin{tabular}{ccccccc}
    \toprule
    \multicolumn{1}{c}{Wave frequencies}&
    \multicolumn{1}{c}{$\;\;\;\;\omega_1=2.5\;\;$}&\multicolumn{1}{c}{$\;\;\omega_1=2.5\;\;$}&\multicolumn{1}{c}{$\;\;\omega_1=2.5\;\;$}&\multicolumn{1}{c}{$\;\;\omega_1=5\;\;$}&\multicolumn{1}{c}{$\;\;\omega_1=5\;\;$}&\multicolumn{1}{c}{$\;\;\omega_1=10\;\;$}\\
    &{$\;\;\;\;\omega_2=5\;\;$}&{$\;\;\omega_2=10\;\;$}&{$\;\;\omega_2=20\;\;$}&{$\;\;\omega_2=10\;\;$}&{$\;\;\omega_2=20\;\;$}&{$\;\;\omega_2=20\;\;$}\\
    \midrule
     $e_{g}$ & 12.57\% & 11.77\% & 12.77\% & 10.40\% & 12.92\% & 17.64\%\\
    $e_{a}$ & 11.98\% & 11.41\% & 12.25\% & 10.21\% & 12.50\% &18.16\% \\
    \bottomrule
    \end{tabular}
    \end{threeparttable}
\end{table}

\subsection{Dependence on the resolution}\label{sec:num_resolu}

By noting the main results in {Theorems \ref{thm:BF_0_app}, \ref{thm:CNN_app} and \ref{thm:gen_0}} on the approximation and generalization analyses of the proposed uncompressed neural networks where the errors are bounded in terms of the resolutions $n_\theta$ and $n_c$ (and $n_\rho=(\omega_2/\omega_1) n_\theta$), we discuss the influence of the resolutions on the approximation and generalization accuracy.
As shown in Table \ref{tab:resolu}, higher resolutions of the inputs, i.e., $n_\theta$, give higher accuracy, while the resolutions of the outputs, i.e., $n_c$, have little influence on the behavior of the proposed neural networks.

\begin{table}[htp!]
  \centering
  \begin{threeparttable}
  \caption{Relative errors of the proposed neural networks with different resolution on the smooth perturbations training datasets where $N_{\rm tr}=900$, $\omega_1=2.5$ and $\omega_2=5$.\label{tab:resolu}}
    \begin{tabular}{ccccccc}
    \toprule
    \multicolumn{1}{c}{Resolutions}&
    \multicolumn{1}{c}{$\;\;\;\;n_\theta=32\;\;$}&\multicolumn{1}{c}{$\;\;n_\theta=32\;\;$}&\multicolumn{1}{c}{$\;\;n_\theta=64\;\;$}&\multicolumn{1}{c}{$\;\;n_\theta=64\;\;$}&\multicolumn{1}{c}{$\;\;n_\theta=64\;\;$}&\multicolumn{1}{c}{$\;\;n_\theta=128\;\;$}\\
    &{$\;\;\;\;n_c=32\;\;$}&{$\;\;n_c=64\;\;$}&{$\;\;n_c=64\;\;$}&{$\;\;n_c=128\;\;$}&{$\;\;n_c=256\;\;$}&{$\;\;n_c=128\;\;$}\\
    \midrule
     $e_{g}$ & 13.11\% & 13.11\% & 12.57\% & 12.35\% & 12.32\% & 11.96\%\\
    $e_{a}$ & 12.45\% & 12.43\% & 11.98\% & 11.79\% & 11.75\% & 11.41\% \\
    \bottomrule
    \end{tabular}
    \end{threeparttable}
\end{table}

\subsection{Dependence on the size of training dataset}

Now, we consider the impact of the training data size on the proposed neural network. The numerical results in Tables \ref{tab:size} and \ref{tab:size2} indicate that the training data size
$N_{\rm tr}=300$ is sufficient for both smooth and nonsmooth perturbations defined in \eqref{eqn:dataset_s} and \eqref{eqn:dataset_nons} to achieve reasonable accuracy. A smaller training set yields lower approximation errors, but performs worse on the generalization; Larger training datasets lead to comparable performance in terms of both approximation and generalization.

\begin{table}[htp!]
  \centering
  \begin{threeparttable}
  \caption{Relative errors of the proposed neural networks on the smooth perturbations training datasets of different size where $n_\theta=n_c=64$, $\omega_1=2.5$ and $\omega_2=5$.\label{tab:size}}
    \begin{tabular}{cccccc}
    \toprule
    \multicolumn{1}{c}{Training size}&
    \multicolumn{1}{c}{$\;\;\;\;\;\;N_{\rm tr}=50\;\;\;\;$}&\multicolumn{1}{c}{$\;\;\;\;N_{\rm tr}=100\;\;\;\;$}&\multicolumn{1}{c}{$\;\;\;\;N_{\rm tr}=300\;\;\;\;$}&\multicolumn{1}{c}{$\;\;\;\;N_{\rm tr}=900\;\;\;\;$}\\
    \midrule
     $e_{g}$ & 16.04\% & 13.57\% & 12.79\%  & 12.57\% \\
    $e_{a}$ & 8.94\% & 9.96\% & 11.86\%  & 11.98\% \\
    \bottomrule
    \end{tabular}
    \end{threeparttable}
\end{table}

\begin{table}[htp!]
  \centering
  \begin{threeparttable}
  \caption{Relative errors of the proposed neural networks on the nonsmooth perturbations training datasets of different size where $n_\theta=n_c=64$, $\omega_1=2.5$ and $\omega_2=5$.\label{tab:size2}}
    \begin{tabular}{cccccc}
    \toprule
    \multicolumn{1}{c}{Training size}&
    \multicolumn{1}{c}{$\;\;\;\;\;\;N_{\rm tr}=50\;\;\;\;$}&\multicolumn{1}{c}{$\;\;\;\;N_{\rm tr}=100\;\;\;\;$}&\multicolumn{1}{c}{$\;\;\;\;N_{\rm tr}=300\;\;\;\;$}&\multicolumn{1}{c}{$\;\;\;\;N_{\rm tr}=900\;\;\;\;$}\\
    \midrule
    $e_{g}$ & 20.55\% & 17.88\% & 16.79\%  & 15.81\% \\
    $e_{a}$ & 13.15\% & 14.57\% & 15.39\%  & 15.39\% \\
    \bottomrule
    \end{tabular}
    \end{threeparttable}
\end{table}

\subsection{Test for the anisotropic media}

It is well-known that neural networks rely heavily on the training and testing datasets. The desirable or even acceptable accuracy is not achievable when the test data falls outside the distribution of the training data. 
We shall discuss the feasibility of the proposed neural networks on test data other than the training-data type.
In particular, we consider the problem of recovering the perturbations from the scattering data generated by an anisotropic scatter, that is solving the problem 
\begin{align*}
\nabla\cdot(A\nabla u^s)+\omega^2 nu^s=- \nabla\cdot(\Gamma\nabla u^i)-\omega^2 \eta u^i, \quad \mbox{where}\quad A=I+\Gamma=I+\begin{bmatrix}
    \gamma_{11} &\gamma_{12}\\\gamma_{12} & \gamma_{22}
\end{bmatrix}
\end{align*}
by the proposed neural network {trained on} the data of isotropic media.
To minimize the loss of features captured from the training dataset and the disruption of the neural network by the unpredictable structure of the test data, we consider a smooth isotropic training dataset without too small inhomogeneities, which is relatively difficult to capture (see Sections \ref{sec:num_capability} and \ref{sec:num_resolu}), generated from \eqref{eqn:dataset_s} with $R_{a,j}\in  \big(1-\max(|x_{a,j}|,|y_{a,j}|)\big)[0.3,1)$ and $R_{n,j}\in \big(1-\max(|x_{n,j}|,|y_{n,j}|)\big)[0.3,1)$. 
The relative approximation and generalization errors of the trained neural network are $g_a=7.40\%$ and $e_g=7.25\%$.
The perturbations $\gamma_{11},\gamma_{12},\gamma_{22}$ and $\eta$ in the testing dataset are generated from the same type of smooth perturbation as that used in the isotropic training dataset. 
The outputs of the trained neural network are the isotropic (i.e., scalar) coefficients. 

\begin{figure}[hbt!]
\centering
\setlength{\tabcolsep}{4pt}
{\tiny\begin{tabular}{c}
\includegraphics[width=0.66\textwidth]{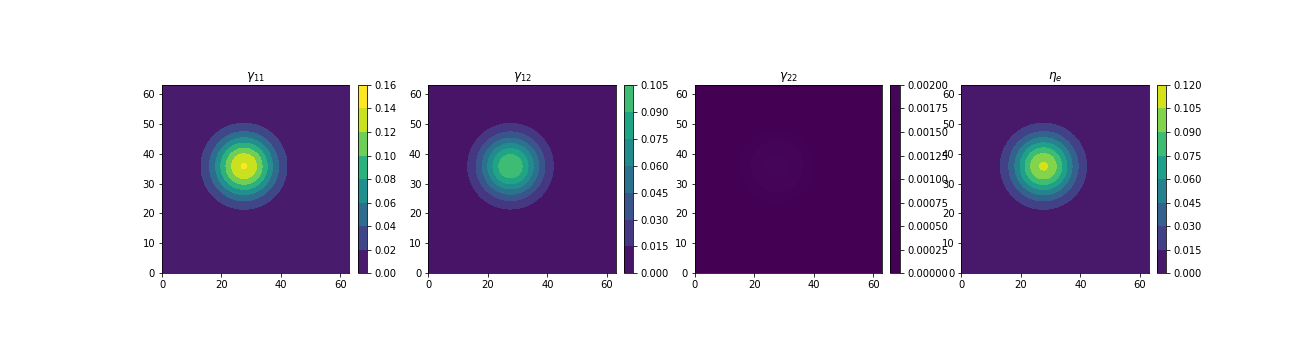} \\
\includegraphics[width=0.66\textwidth]{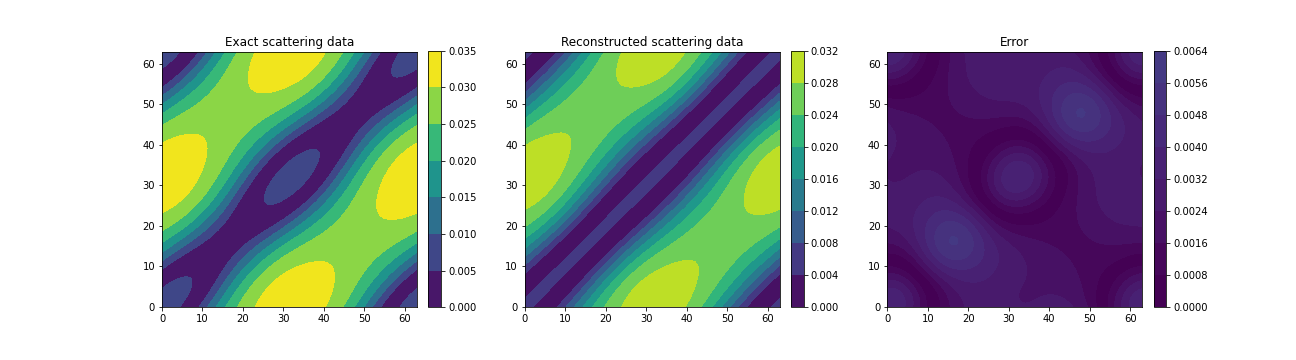}\\
\end{tabular}}
\caption{Exact perturbations $\gamma_{11}$, $\gamma_{12}$, $\gamma_{22}$ and $\eta_e$, and the exact and reconstructed scattering data.\label{fig:aniso_per1_1}}
\end{figure}

\begin{figure}[hbt!]
\centering
\setlength{\tabcolsep}{4pt}
{\tiny\begin{tabular}{c}
\includegraphics[width=0.66\textwidth]{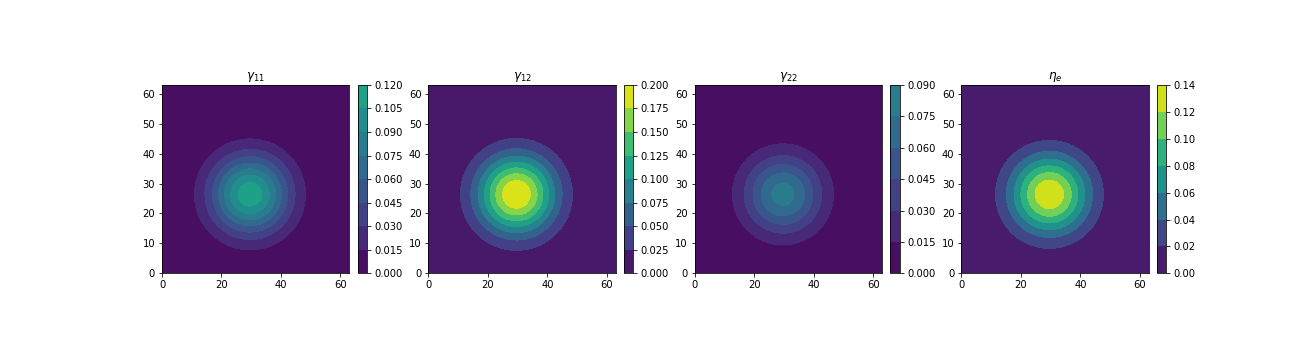} \\
\includegraphics[width=0.66\textwidth]{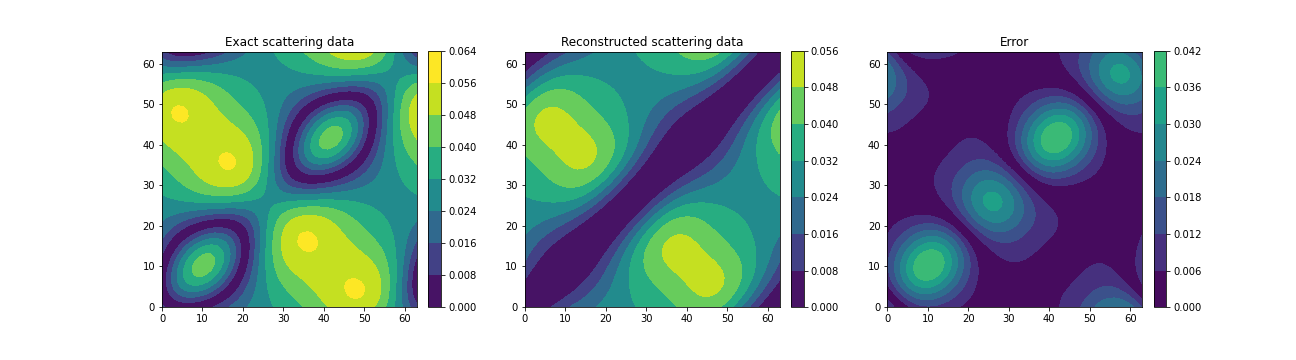}\\
\end{tabular}}
\caption{Exact perturbations $\gamma_{11}$, $\gamma_{12}$, $\gamma_{22}$ and $\eta_e$, and the exact and reconstructed scattering data.\label{fig:aniso_per1_2}}
\end{figure}

\begin{figure}[hbt!]
\centering
\setlength{\tabcolsep}{4pt}
{\tiny\begin{tabular}{c}
\includegraphics[width=0.66\textwidth]{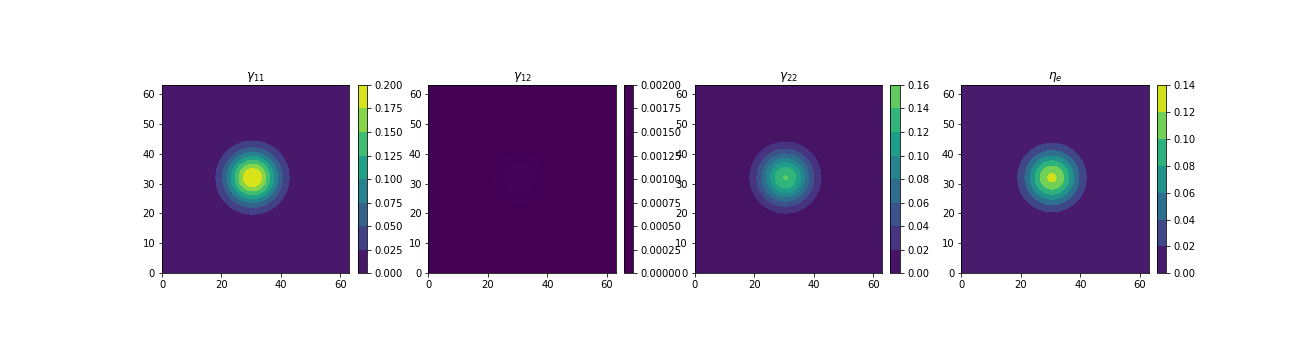} \\
\includegraphics[width=0.66\textwidth]{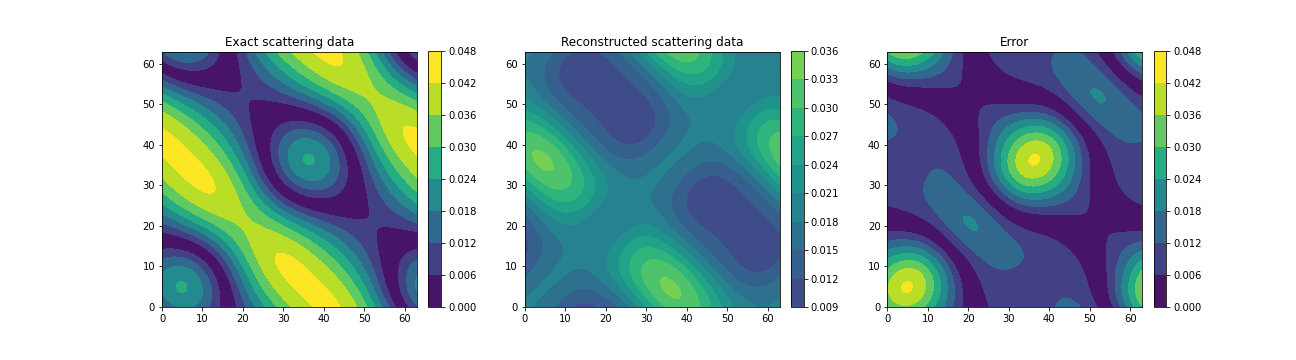}\\
\end{tabular}}
\caption{Exact perturbations $\gamma_{11}$, $\gamma_{12}$, $\gamma_{22}$ and $\eta_e$, and the exact and reconstructed scattering data.\label{fig:aniso_per1_3}}
\end{figure}

\begin{figure}[hbt!]
\centering
\setlength{\tabcolsep}{4pt}
{\tiny\begin{tabular}{c}
\includegraphics[width=0.66\textwidth]{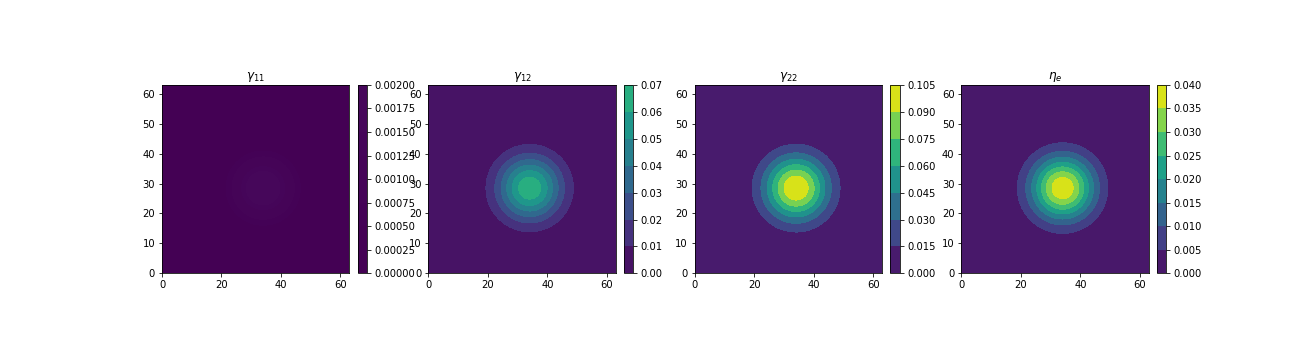} \\
\includegraphics[width=0.66\textwidth]{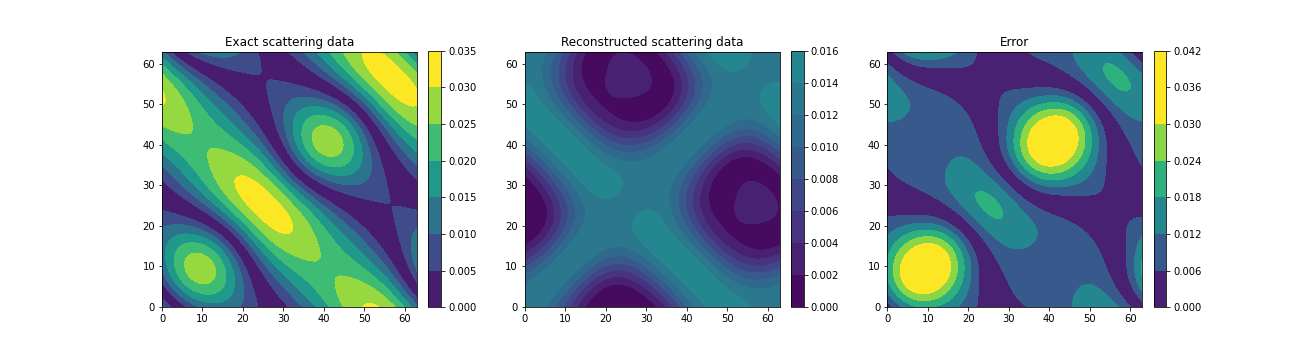}\\
\end{tabular}}
\caption{Exact perturbations $\gamma_{11}$, $\gamma_{12}$, $\gamma_{22}$ and $\eta_e$, and the exact and reconstructed scattering data.\label{fig:aniso_per1_4}}
\end{figure}

As discussed in Section \ref{sec:intro}, in $\mathbb{R}^2$, there is a unique isotropic medium (with scalar coefficients) which produces the same scattering data as the true anisotropic medium. 
We refer to this isotropic medium as the isotropic representation of the anisotropic medium and measure the capability of the trained neural network by the accuracy of constructing this isotropic representation.
However, the isotropic representation cannot be directly computed from a given anisotropic medium.
As a result, we compare the scattering data generated by the reconstructed coefficients with the exact scattering data and measure the distance between them.
The numerical results for the anisotropic medium where $\gamma_{11},\gamma_{12},\gamma_{22}$ and $\eta$ are generated by \eqref{eqn:dataset_s} with $J=1$ (one Gaussian function) and the same support (the region of significant parts) are presented in Figs. \ref{fig:aniso_per1_1}--\ref{fig:aniso_per1_4}.

\begin{figure}[hbt!]
\centering
\setlength{\tabcolsep}{4pt}
{\tiny\begin{tabular}{c}
\includegraphics[width=0.66\textwidth]{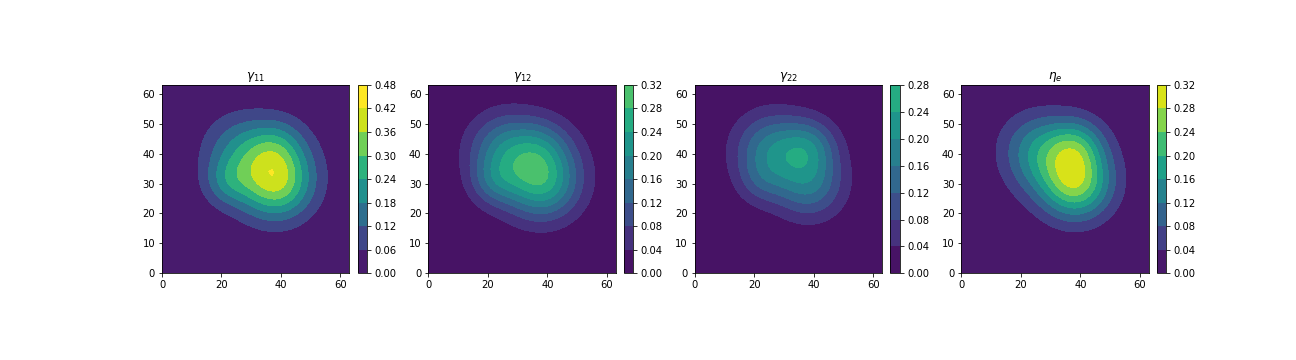} \\
\includegraphics[width=0.66\textwidth]{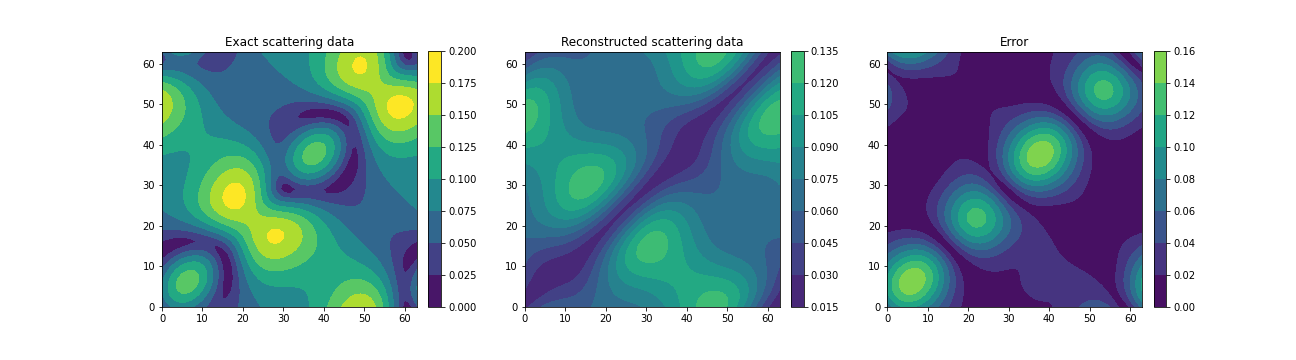}\\
\end{tabular}}
\caption{Exact perturbations $\gamma_{11}$, $\gamma_{12}$, $\gamma_{22}$ and $\eta_e$, and the exact and reconstructed scattering data.\label{fig:aniso_per5_1}}
\end{figure}

\begin{figure}[hbt!]
\centering
\setlength{\tabcolsep}{4pt}
{\tiny\begin{tabular}{c}
\includegraphics[width=0.66\textwidth]{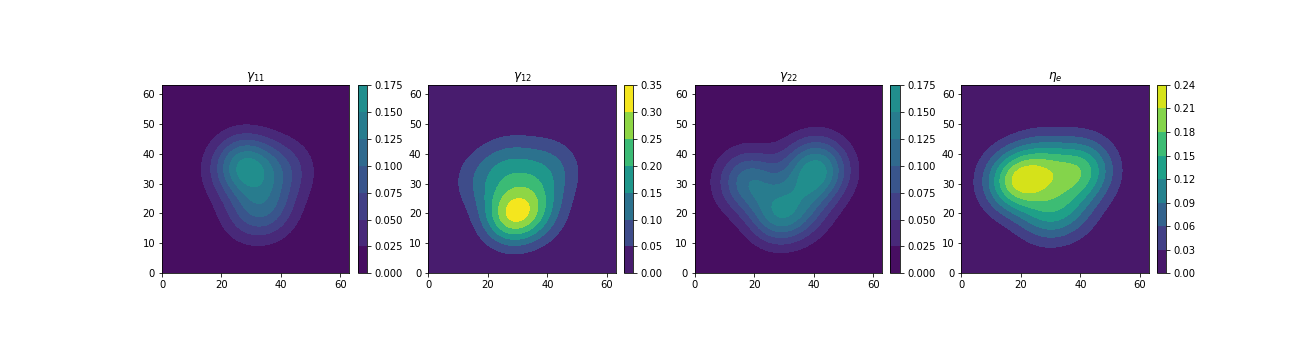} \\
\includegraphics[width=0.66\textwidth]{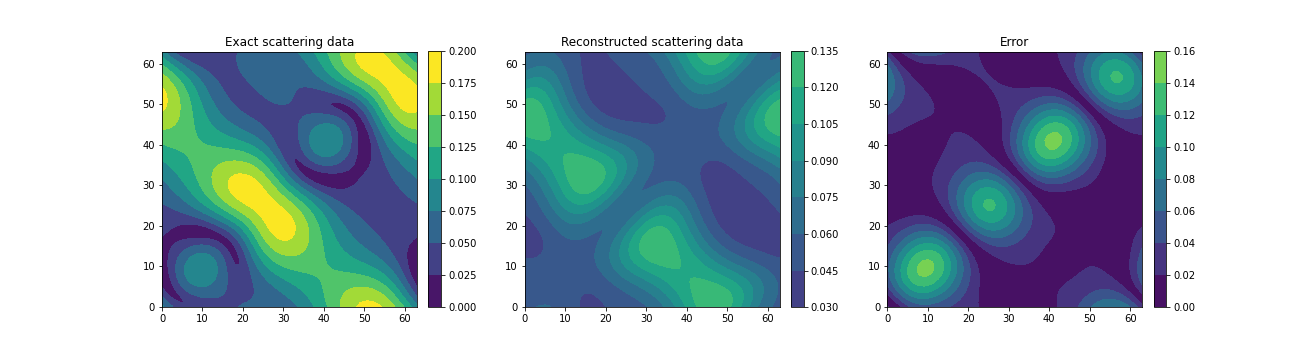}\\
\end{tabular}}
\caption{Exact perturbations $\gamma_{11}$, $\gamma_{12}$, $\gamma_{22}$ and $\eta_e$, and the exact and reconstructed scattering data.\label{fig:aniso_per5_2}}
\end{figure}

\begin{figure}[hbt!]
\centering
\setlength{\tabcolsep}{4pt}
{\tiny\begin{tabular}{c}
\includegraphics[width=0.66\textwidth]{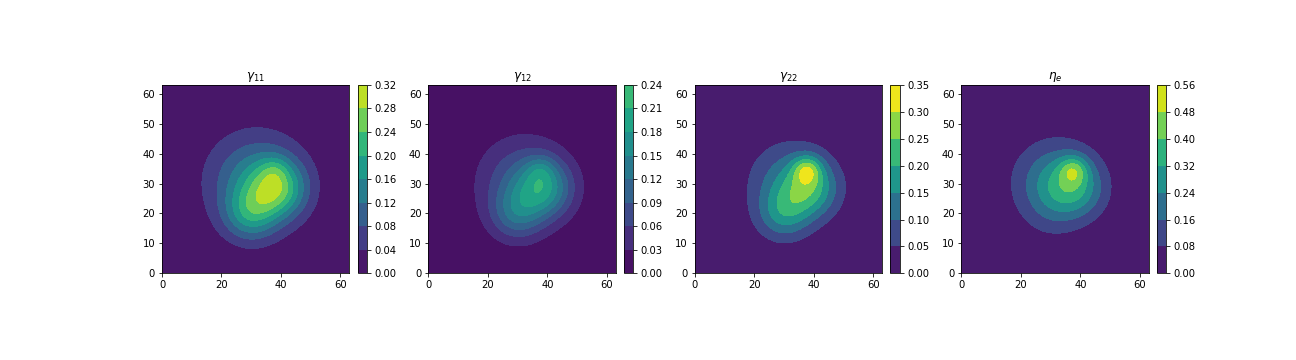} \\
\includegraphics[width=0.66\textwidth]{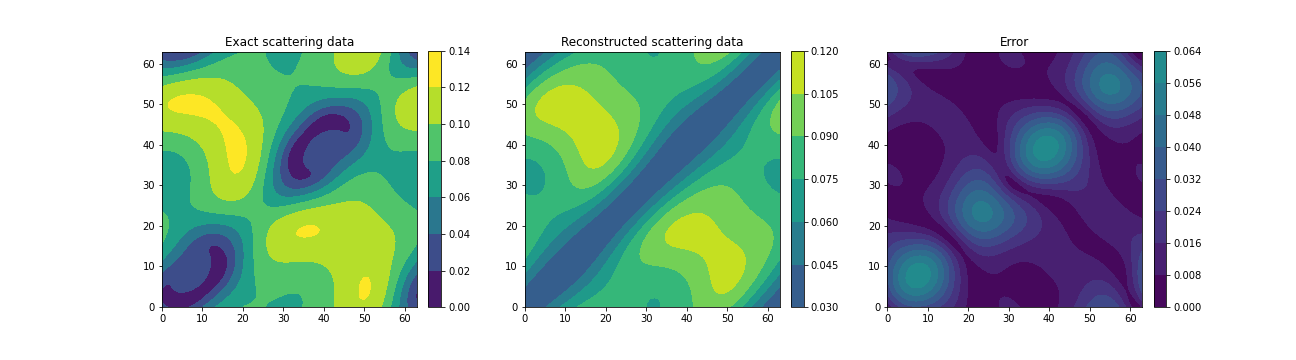}\\
\end{tabular}}
\caption{Exact perturbations $\gamma_{11}$, $\gamma_{12}$, $\gamma_{22}$ and $\eta_e$, and the exact and reconstructed scattering data.\label{fig:aniso_per5_3}}
\end{figure}

\begin{figure}[hbt!]
\centering
\setlength{\tabcolsep}{4pt}
{\tiny\begin{tabular}{c}
\includegraphics[width=0.66\textwidth]{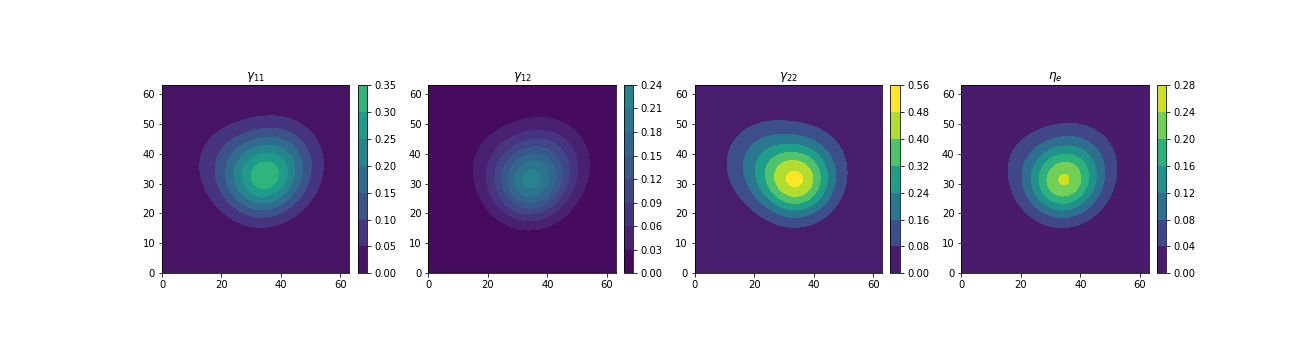} \\
\includegraphics[width=0.66\textwidth]{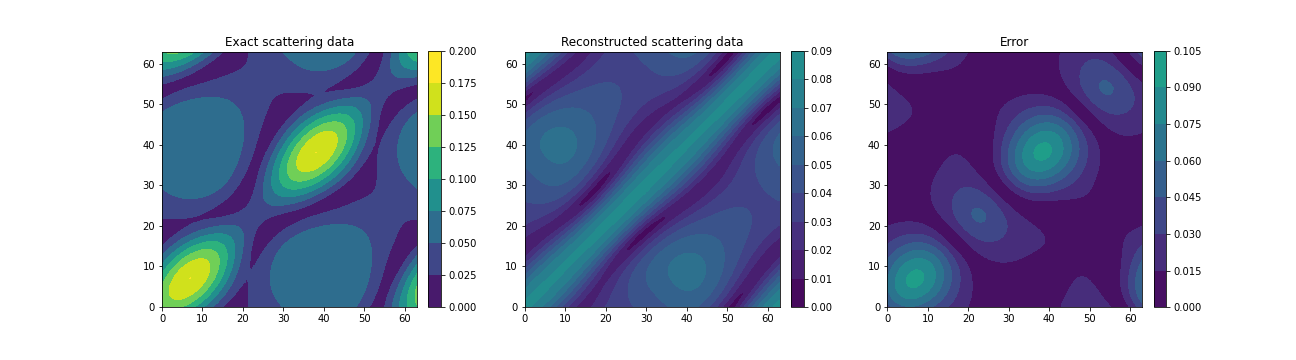}\\
\end{tabular}}
\caption{Exact perturbations $\gamma_{11}$, $\gamma_{12}$, $\gamma_{22}$ and $\eta_e$, and the exact and reconstructed scattering data.\label{fig:aniso_per5_4}}
\end{figure}

The results indicate that the trained neural network can construct the isotropic representations in some cases (e.g., Figs \ref{fig:aniso_per1_1} and \ref{fig:aniso_per1_2}). However, for other cases (e.g., Figs \ref{fig:aniso_per1_3} and \ref{fig:aniso_per1_4}), it will introduce an error to the corresponding scattering data, especially on the forward incident scattering parts. Similar observations hold for the perturbations generated by \eqref{eqn:dataset_s} with $J=5$ (see Figs. \ref{fig:aniso_per5_1}--\ref{fig:aniso_per5_4}).

%%%%%%%%%%%%%%%%%%%%%%%%%%%%%%%%%%%%%%%%%%%%%%%%%%%%%%%%%%%%%%%%%%

\section{Concluding remarks}\label{sec:conc}

In this work, we designed two types of combined neural networks (uncompressed and compressed) to reconstruct {two function-valued coefficients} in the Helmholtz equation for the inverse scattering problem from the scattering data at two different frequencies by approximating the regularized pseudo-inverse of the linearized forward operator for the direct scattering problem.
We first decomposed the regularized pseudo-inverse into two components, i.e., the Fourier integral operator and convolution-type inverse operator, and studied the properties of each component separately. 
We constructed a fully connected neural network of shared weights layers which can be further approximated by the butterfly factorization for approximating the Fourier integral operator, and then combined it with a convolutional neural network with scaled residual layers for approximating the convolution-type inverse operator.
We also discussed the approximation and generalization capabilities of the proposed neural networks, where the results indicate that the proposed neural networks
can effectively approximate the inverse process with desirable generalization when sufficient training data and parameters are utilized.
Finally, we presented several numerical experiments which show the feasibility of the proposed neural networks for recovering different types of isotropic inhomogeneous media and that the trained neural networks can also be used to reconstruct the isotropic representation of certain types of anisotropic media.

The neural networks proposed in this work are designed for simulating the linearized model, which is considered as a limitation of the method (even if their expressiveness can be improved by increasing the nonlinearity of the structures).
Improving the neural networks to simulate the inverse process of the true model (or better approximators) with the support of theoretical analysis is an important topic that desires to be investigated.
We leave this interesting question for future work.

%%%%%%%%%%%%%%%%%%%%%%%%%%%%%%%%%%%%%%%%%%%%%%%%%%%%%%%%%%%%%%%%%%

\section*{Acknowledgements}

The author is grateful to Professor Fioralba Cakoni and Professor Michael Vogelius for {stimulating} discussions on this work, and to Professor Peter Monk for providing parts of the code for generating the scattering data. 
The author also thanks Professor Qin Li's group for sharing the code from their previous work \cite{ZhangZepeda-NunezLi:2024}.
Part of the work was completed during the visit to Hong Kong and Beijing supported by AMS-Simons Travel Grant.

%%%%%%%%%%%%%%%%%%%%%%%%%%%%%%%%%%%%%%%%%%%%%%%%%%%%%%%%%%%%%%%%%%

\section*{Declarations}
\textbf{Conflict of interest}  \;\;The author declares no competing interests.

%%%%%%%%%%%%%%%%%%%%%%%%%%%%%%%%%%%%%%%%%%%%%%%%%%%%%%%%%%%%%%%%%%

\appendix
\section{Auxiliary estimates}\label{app:estimate}
In this appendix, we give several auxiliary estimates that have been used in the approximation and generalization analyses.
First, we state the estimates on the modulus of continuity and $L^\infty$ norms of some functions in Lemmas \ref{lem:est1}, \ref{lem:filter} and \ref{lem:det}. We define the modulus of continuity $w_{f}$ of any function $f$ by $w_{f}(s)=\sup_{|x-y|_{\infty}\leq s}|f(x)-f(y)|$ where $|\cdot|_{\infty}$ denotes the maximum magnitude of the entries in the vector or matrix.
\begin{lemma}\label{lem:est1}
Let the integral kernel $k(t,\rho)=e^{-i\omega_1 \rho \cos t}$ and $\Lambda^{\omega}$ be the approximate far-field pattern defined in \eqref{eqn:F1F2}. 
Then for any $\rho_n\geq 0$, there hold
\begin{align*}
w_{\overline{k(\cdot,\rho_n)}}(s)=w_{k(\cdot,\rho_n)}(s)\leq&\omega_1 \rho_n s, \quad w_{\Lambda^{\omega}}(s)\leq 2c_{\omega} (\omega R\|\eta(y)\|_{L^1(\mathcal{B}_R)}+(\omega R+1)\|\gamma(y)\|_{L^1(\mathcal{B}_R)})s,\\
\|\bar{k}\|_{L^{\infty}}=\|k\|_{L^{\infty}} \leq& 1,\qquad \;\;\|\Lambda^{\omega}\|_{L^{\infty}}=c_{\omega}(\|\eta(y)\|_{L^1(\mathcal{B}_R)}+\|\gamma(y)\|_{L^1(\mathcal{B}_R)}).
\end{align*}
\end{lemma}
\begin{proof}
First, from the definition of the integral kernel $k(t,\rho_n)=e^{-i\omega_1 \rho_n \cos t}$, we can derive that $\|\bar{k}\|_{L^{\infty}}=\|k\|_{L^{\infty}} \leq 1$ and
\begin{align*}
w_{\overline{k(\cdot,\rho_n)}}(s)=w_{k(\cdot,\rho_n)}(s)=&\sup_{|t-t'|\leq s}|e^{-i\omega_1 \rho_n \cos t}-e^{-i\omega_1 \rho_n \cos t'}|\leq\sup_{|t-t'|\leq s}|e^{-i\omega_1 \rho_n (\cos t-\cos t')}-1|\\
\leq&\sup_{|t-t'|\leq s}|-i\omega_1 \rho_n (\cos t-\cos t')|\leq \sup_{|t-t'|\leq s}\omega_1 \rho_n |t-t'|=\omega_1 \rho_n s.
\end{align*}

Similarly, for the far-field pattern $\Lambda^{\omega}$ given by \eqref{eqn:F1F2}, 
we have
\begin{align*}
\|\Lambda^{\omega}\|_{L^{\infty}}=&\sup_{\hat{x},\zeta\in\mathcal{S}}\big|c_{\omega}\int_{\mathcal{B}_R}(\eta(y)-(\hat{x}\cdot \zeta)   \gamma(y))e^{i\omega (\zeta-\hat{x})\cdot y}{\rm d}y\big|\\
\leq &c_{\omega}\sup_{\hat{x},\zeta\in\mathcal{S}}\int_{\mathcal{B}_R}|(\eta(y)-(\hat{x}\cdot \zeta)   \gamma(y))e^{i\omega (\zeta-\hat{x})\cdot y}|{\rm d}y\\
\leq &c_{\omega}\sup_{\hat{x},\zeta\in\mathcal{S}}\int_{\mathcal{B}_R}(|\eta(y)|+|\gamma(y)|){\rm d}y
=c_{\omega}(\|\eta(y)\|_{L^1(\mathcal{B}_R)}+\|\gamma(y)\|_{L^1(\mathcal{B}_R)}).
\end{align*}
Then, with the assumption $\bar{\Omega}\subset [0,2\pi]\times[0,R]$, we rewrite $\Lambda^{\omega}$ on the polar coordinates as
\begin{align*}
\Lambda^{\omega}(\theta_{\hat{x}},\theta_{\zeta})
%=&c_{\omega}\int_{\mathcal{B}_R}(\eta(y)-(\hat{x}\cdot \zeta)   \gamma(y))e^{i\omega (\zeta-\hat{x})\cdot y}{\rm d}y\\
=&c_{\omega}\int_0^{2\pi}\int_0^R(\eta(y)-\cos(\theta_{\hat{x}}-\theta_{\zeta})  \gamma(y))e^{i{\omega}\rho(\cos(\theta_{\zeta}-\theta)-\cos(\theta_{\hat{x}}-\theta))}\rho{\rm d}\rho{\rm d}\theta,
\end{align*}
and then estimate $w_{\Lambda^{\omega}}(s)$ by
\begin{align*}
w_{\Lambda^{\omega}}(s)=&\sup_{|\theta_{\hat{x}}-\theta_{\hat{x}}'|\leq s,\;|\theta_{\zeta}-\theta_{\zeta}'|\leq s}|\Lambda^{\omega}(\theta_{\hat{x}},\theta_{\zeta})-\Lambda^{\omega}(\theta_{\hat{x}}',\theta_{\zeta}')|\\
\leq& \sup_{\theta_{\zeta}}\sup_{|\theta_{\hat{x}}-\theta_{\hat{x}}'|\leq s}|\Lambda^{\omega}(\theta_{\hat{x}},\theta_{\zeta})-\Lambda^{\omega}(\theta_{\hat{x}}',\theta_{\zeta})|+\sup_{\theta_{\hat{x}}'}\sup_{|\theta_{\zeta}-\theta_{\zeta}'|\leq s}|\Lambda^{\omega}(\theta_{\hat{x}}',\theta_{\zeta})-\Lambda^{\omega}(\theta_{\hat{x}}',\theta_{\zeta}')|,
\end{align*}
where
\begin{align*}
&\sup_{|\theta_{\hat{x}}-\theta_{\hat{x}}'|\leq s}|\Lambda^{\omega}(\theta_{\hat{x}},\theta_{\zeta})-\Lambda^{\omega}(\theta_{\hat{x}}',\theta_{\zeta})|\\
=&c_{\omega}\sup_{|\theta_{\hat{x}}-\theta_{\hat{x}}'|\leq s}|\int_0^{2\pi}\int_0^R\big((\eta(y)-\cos(\theta_{\hat{x}}-\theta_{\zeta})  \gamma(y))e^{-i{\omega}\rho\cos(\theta_{\hat{x}}-\theta)}\\
&\qquad\qquad\qquad\qquad\quad-(\eta(y)-\cos(\theta_{\hat{x}}'-\theta_{\zeta})  \gamma(y))e^{-i{\omega}\rho\cos(\theta_{\hat{x}}'-\theta)}\big)e^{i\rho\cos(\theta_{\zeta}-\theta)}\rho{\rm d}\rho{\rm d}\theta|\\
%\leq&c_{\omega}\sup_{|\theta_{\hat{x}}-\theta_{\hat{x}}'|\leq s}|\int_0^{2\pi}\int_0^R\big((\eta(y)-\cos(\theta_{\hat{x}}-\theta_{\zeta})  \gamma(y))(e^{-i\omega\rho\cos(\theta_{\hat{x}}-\theta)}-e^{-i\omega\rho\cos(\theta_{\hat{x}}'-\theta)})\\
%&\qquad\qquad\qquad\qquad\quad+(\cos(\theta_{\hat{x}}'-\theta_{\zeta})-\cos(\theta_{\hat{x}}-\theta_{\zeta}))  \gamma(y)e^{-i\omega\rho\cos(\theta_{\hat{x}}'-\theta)}\big)e^{i\omega\rho\cos(\theta_{\zeta}-\theta)}\rho{\rm d}\rho{\rm d}\theta|
\leq&c_{\omega}\sup_{|\theta_{\hat{x}}-\theta_{\hat{x}}'|\leq s}\int_0^{2\pi}\int_0^R\big(|(\eta(y)-\cos(\theta_{\hat{x}}-\theta_{\zeta})  \gamma(y))(e^{-i\omega\rho\cos(\theta_{\hat{x}}-\theta)}-e^{-i\omega\rho\cos(\theta_{\hat{x}}'-\theta)})|\\
&\qquad\qquad\qquad\qquad\qquad\qquad\quad+|(\cos(\theta_{\hat{x}}'-\theta_{\zeta})-\cos(\theta_{\hat{x}}-\theta_{\zeta}))  \gamma(y)e^{-i\omega\rho\cos(\theta_{\hat{x}}'-\theta)}|\big)\rho{\rm d}\rho{\rm d}\theta\\
\leq&c_{\omega}\sup_{|\theta_{\hat{x}}-\theta_{\hat{x}}'|\leq s}\int_0^{2\pi}\int_0^R\big((|\eta(y)|+|\gamma(y)|)|e^{-i\omega\rho(\cos(\theta_{\hat{x}}-\theta)-\cos(\theta_{\hat{x}}'-\theta))}-1|\\
&\qquad\qquad\qquad\qquad\qquad\qquad\quad+|\cos(\theta_{\hat{x}}'-\theta_{\zeta})-\cos(\theta_{\hat{x}}-\theta_{\zeta})||\gamma(y)|\big)\rho{\rm d}\rho{\rm d}\theta\\
%\leq&c_{\omega}\sup_{|\theta_{\hat{x}}-\theta_{\hat{x}}'|\leq s}\int_0^{2\pi}\int_0^R\big((|\eta(y)|+|\gamma(y)|)\omega\rho|\theta_{\hat{x}}-\theta_{\hat{x}}'|+|\theta_{\hat{x}}'-\theta_{\hat{x}}||\gamma(y)|\big)\rho{\rm d}\rho{\rm d}\theta\\
\leq&c_{\omega}\sup_{|\theta_{\hat{x}}-\theta_{\hat{x}}'|\leq s}\int_0^{2\pi}\int_0^R\big((|\eta(y)|+|\gamma(y)|)\omega\rho+|\gamma(y)|\big)|\theta_{\hat{x}}'-\theta_{\hat{x}}|\rho{\rm d}\rho{\rm d}\theta\\
\leq&c_{\omega} s \int_0^{2\pi}\int_0^R(\omega R|\eta(y)|+(\omega R+1)|\gamma(y)|)\rho{\rm d}\rho{\rm d}\theta
=c_{\omega} (\omega R\|\eta(y)\|_{L^1(\mathcal{B}_R)}+(\omega R+1)\|\gamma(y)\|_{L^1(\mathcal{B}_R)})s.
\end{align*}
Similarly, there holds
\begin{align*}
\sup_{|\theta_{\zeta}-\theta_{\zeta}'|\leq s}|\Lambda^{\omega}(\theta_{\hat{x}}',\theta_{\zeta})-\Lambda^{\omega}(\theta_{\hat{x}}',\theta_{\zeta}')|\leq c_{\omega} (\omega R\|\eta(y)\|_{L^1(\mathcal{B}_R)}+(\omega R+1)\|\gamma(y)\|_{L^1(\mathcal{B}_R)})s.
\end{align*}
This completes the proof.
\end{proof}

We give the estimates and properties of the filter $g_{j{j'}}^\omega(y)$ and the kernels $g_{de}$, $g_2$ and ${\bf g}$ defined in Proposition \ref{prop:F*F_formula} in the following two lemmas.
\begin{lemma}\label{lem:filter}
For any $j,j'=1,2$, let $g_{j{j'}}^\omega(y)$ be the filter given by 
$$g_{j{j'}}^\omega(y)=(-1)^{j+{j'}}c_\omega^2 \int_{\mathbb{S}\times\mathbb{S}} (\hat{x}\cdot \zeta)^{4-j-{j'}} \cos(\omega (\hat{x}-\zeta)\cdot y) {\rm d}\hat{x}\, {\rm d}\zeta:=(-1)^{j+{j'}}c_\omega^2{\rm I}_\omega^{j,j'}(y),$$
then $g_{j{j'}}^\omega(y)$ is a real-valued radially symmetric function satisfying
\begin{align*}
w_{g_{j{j'}}^\omega}(s)\leq&8\sqrt{2} \pi^2\omega c_\omega^2 s \quad \mbox{and}\quad \|g_{j{j'}}^{\omega}\|_{L^{\infty}} \leq 4\pi^2 c_\omega^2.
\end{align*}
Furthermore, the convolutions $g_{11}^{\omega}*g_{22}^{\omega'}$ and $g_{12}^{\omega}*g_{12}^{\omega'}$ are also radially symmetric. 
\end{lemma}
\begin{proof}
It is clear that
$\|g_{j{j'}}^{\omega}\|_{L^{\infty}}\leq c_\omega^2\sup_{y}|{\rm I}_\omega^{j,j'}(y)|\leq c_\omega^2 \int_{\mathbb{S}\times\mathbb{S}} 1{\rm d}\hat{x}\, {\rm d}\zeta\leq 4\pi^2 c_\omega^2$.
Then, we bound the modulus of continuity of $g_{j{j'}}$ by
\begin{align*}
w_{g_{j{j'}}^\omega}(s)=&\sup_{|y-y'|_\infty\leq s}|(-1)^{j+{j'}}c_\omega^2 ({\rm I}_\omega^{j,j'}(y)-{\rm I}_\omega^{j,j'}(y'))|
\leq c_\omega^2 \sup_{|y-y'|_\infty\leq s}|{\rm I}_\omega^{j,j'}(y)-{\rm I}_\omega^{j,j'}(y')|\\
\leq&c_\omega^2 \sup_{|y-y'|_\infty\leq s}\int_{\mathbb{S}\times\mathbb{S}} |\hat{x}\cdot \zeta|^{4-j-{j'}} |\cos(\omega (\hat{x}-\zeta)\cdot y)-\cos(\omega (\hat{x}-\zeta)\cdot y')|{\rm d}\hat{x}\, {\rm d}\zeta\\
\leq&c_\omega^2 \sup_{|y-y'|_\infty\leq s}\int_{\mathbb{S}\times\mathbb{S}}|\omega (\hat{x}-\zeta)\cdot (y-y')|{\rm d}\hat{x}\, {\rm d}\zeta
\leq\omega c_\omega^2 \sup_{|y-y'|_\infty\leq s}\int_{\mathbb{S}\times\mathbb{S}}  \|\hat{x}-\zeta\|_2\|y-y'\|_2{\rm d}\hat{x}\, {\rm d}\zeta\\
\leq&\sqrt{2}\omega c_\omega^2 s \int_{\mathbb{S}\times\mathbb{S}}  2{\rm d}\hat{x}\, {\rm d}\zeta
\leq 8\sqrt{2} \pi^2\omega c_\omega^2 s .
\end{align*}

Next, we shed light on the symmetry of $g_{j{j'}}^\omega$. To this end, we convert the Cartesian coordinates to the polar coordinates via $\hat{x}=(\cos\theta_{\hat{x}},\sin\theta_{\hat{x}})$, $\zeta=(\cos\theta_{\zeta},\sin\theta_{\zeta})$ and $y=(\rho \cos\theta,\rho \sin\theta)$. 
By \eqref{eqn:tri}, we have
\begin{align*}
{\rm I}_\omega^{j,j'}=&\int_{0}^{2\pi}\int_{0}^{2\pi} (\cos(\theta_{\hat{x}}-\theta_{\zeta}))^{4-j-{j'}}  \cos\Big(\omega \rho \big(\cos(\theta_{\hat{x}}-\theta)-\cos(\theta_{\zeta}-\theta)\big)\Big){\rm d}\theta_{\hat{x}}{\rm d}\theta_{\zeta}\\
=&\int_{0}^{2\pi}\int_{0}^{2\pi} (\cos(\theta_{\hat{x}}-\theta_{\zeta}))^{4-j-{j'}}  \cos\Big(\omega \rho \big(\cos(\theta_{\hat{x}})-\cos(\theta_{\zeta})\big)\Big){\rm d}\theta_{\hat{x}}{\rm d}\theta_{\zeta}
\end{align*}
which is independent of $\theta$. This indicates the radially symmetry of the filter $g_{j{j'}}^\omega(y)$ and the convolutions $g_{11}^{\omega}*g_{22}^{\omega'}$ and $g_{12}^{\omega}*g_{12}^{\omega'}$. 
This completes the proof.
\end{proof}

\begin{lemma}\label{lem:det}
Let $g_{de}=\sum_{i,i'=1}^2(g_{11}^{\omega_i}+\frac\alpha2\delta)*(g_{22}^{\omega_{i'}}+\frac\alpha2\delta) -g_{12}^{\omega_i}*g_{12}^{\omega_{i'}}$, $g_2=\mathcal{F}^{-1}\big(\mathcal{F}(g_{de})^{-1}\big)$ and ${\bf g}$ be the kernels defined in Proposition \ref{prop:F*F_formula} with 
$$g_{j{j'}}^\omega(y)=(-1)^{j+{j'}}c_\omega^2 \int_{\mathbb{S}\times\mathbb{S}} (\hat{x}\cdot \zeta)^{4-j-{j'}} \cos(\omega (\hat{x}-\zeta)\cdot y) {\rm d}\hat{x}\, {\rm d}\zeta:=(-1)^{j+{j'}}c_\omega^2 {\rm I}_\omega^{j,j'}(y).$$ 
Then, $\det(\mathcal{F}({\bf g}))=\mathcal{F}(g_{de})\geq\alpha^2$ and there exists a decomposition $g_2=\alpha^{-2}\delta-g_{2,\alpha}$ with real-valued function $g_{2,\alpha}$ such that %$\mathcal{F}(g_{2,\alpha})\in[0, \alpha^{-2}]$ and $\mathcal{F}(g_{2,\alpha})$ vanishes outside the disc $\overline{B_{2\omega_2}}$. 
\begin{align*}
w_{g_{2,\alpha}}(s)\leq&\frac{4\sqrt{2}\omega_2^3}{3\pi\alpha^{2}}s \quad \mbox{and}\quad \|g_{2,\alpha}\|_{L^{\infty}} \leq \frac{\omega_2^2}{\pi\alpha^{2}}.
\end{align*}
\end{lemma}

\begin{proof}
By the definition of $g_{de}$, we can decompose $g_{de}$ into 
\begin{align}
g_{de}=&\sum_{i,i'=1}^2\big((g_{11}^{\omega_i}+\frac\alpha2\delta)*(g_{22}^{\omega_{i'}}+\frac\alpha2\delta) -g_{12}^{\omega_i}*g_{12}^{\omega_{i'}}\big)\nonumber\\
=&(\sum_{i=1}^2 g_{11}^{\omega_i})*(\sum_{i=1}^2 g_{22}^{\omega_{i}})-(\sum_{i=1}^2 g_{12}^{\omega_i})*(\sum_{i=1}^2 g_{12}^{\omega_{i}})+\alpha\sum_{i=1}^2(g_{11}^{\omega_i}+g_{22}^{\omega_{i}})+\alpha^2\delta\nonumber
:=g_{de,\alpha}+\alpha^2\delta,\nonumber%\label{eqn:g_de_decom}
\end{align}
where $g_{jj'}^{\omega}=(-1)^{j+{j'}}c_\omega^2{\rm I}_\omega^{j,j'}(y)$, for any $j,j'=1,2$, with
\begin{align*}
{\rm I}_\omega^{j,j'}(y)=&\int_{\mathbb{S}\times\mathbb{S}} (\hat{x}\cdot \zeta)^{4-j-{j'}} \cos(\omega (\hat{x}-\zeta)\cdot y){\rm d}\hat{x}\, {\rm d}\zeta=\int_{\mathbb{S}\times\mathbb{S}} (\hat{x}\cdot \zeta)^{4-j-{j'}} e^{i\omega (\hat{x}-\zeta)\cdot y}{\rm d}\hat{x}{\rm d}\zeta.
\end{align*}
The Fourier transform of ${\rm I}_\omega^{j,j'}$ is given by
\begin{align*}
[\mathcal{F}(I_\omega^{j,j'}(y))](\xi)
=&\int_{\mathbb{S}\times\mathbb{S}} (\hat{x}\cdot \zeta)^{4-j-{j'}} [\mathcal{F}(e^{i\omega (\hat{x}-\zeta)\cdot y})](\xi){\rm d}\hat{x}{\rm d}\zeta=(2\pi)^2\int_{\mathbb{S}\times\mathbb{S}} (\hat{x}\cdot \zeta)^{4-j-{j'}} \delta(\omega (\hat{x}-\zeta)-\xi){\rm d}\hat{x}{\rm d}\zeta\\
=&(2\pi)^2\int_{\mathbb{S}\times\mathbb{S}} (\hat{x}\cdot \zeta)^{4-j-{j'}} \delta(\omega (\hat{x}-\zeta)+\xi){\rm d}\hat{x}{\rm d}\zeta
=[\mathcal{F}(I_\omega^{j,j'}(y))](-\xi)%\\
%=&(\frac{2\pi}{\omega})^2\int_{\mathbb{S}\times\mathbb{S}} (\hat{x}\cdot \zeta)^{4-j-{j'}} \delta( \hat{x}-\zeta-\frac{\xi}{\omega}){\rm d}\hat{x}{\rm d}\zeta,
\end{align*}
and thus
\begin{align*}
\sum_{i=1}^2[\mathcal{F}(g_{jj'}^{\omega_i})](-\xi)=\sum_{i=1}^2[\mathcal{F}(g_{jj'}^{\omega_i})](\xi)=&(-1)^{j+{j'}}\sum_{i=1}^2 c_{\omega_i}^2[\mathcal{F}({\rm I}_{\omega_i}^{j,j'})](\xi)\\
=&(-1)^{j+{j'}}(2\pi)^2\int_{\mathbb{S}\times\mathbb{S}} (\hat{x}\cdot \zeta)^{4-j-{j'}} \sum_{i=1}^2 c_{\omega_i}^2 \delta(\omega_i (\hat{x}-\zeta)-\xi){\rm d}\hat{x}{\rm d}\zeta.
\end{align*}
Then, the Cauchy-Schwarz inequality indicates that
\begin{align*}
\mathcal{F}(g_{de,\alpha})=&\big(\sum_{i=1}^2\mathcal{F}(g_{11}^{\omega_i})\big)\big(\sum_{i=1}^2\mathcal{F}(g_{22}^{\omega_{i}})\big)-\big(\sum_{i=1}^2\mathcal{F}(g_{12}^{\omega_i})\big)^2+\alpha\sum_{i=1}^2\big(\mathcal{F}(g_{11}^{\omega_i})+\mathcal{F}(g_{22}^{\omega_{i}})\big)\geq 0,
\end{align*}
and thus $[\mathcal{F}(g_{de})](\xi)=[\mathcal{F}(g_{de})](-\xi)\geq\alpha^2$.
%When $\xi=0$, there holds 
%\begin{align*}
%[\mathcal{F}(I_\omega^{j,j'}(y))](0)
%=&(\frac{2\pi}{\omega})^2\int_{\mathbb{S}\times\mathbb{S}} (\hat{x}\cdot \zeta)^{4-j-{j'}} \delta( \hat{x}-\zeta){\rm d}\hat{x}{\rm d}\zeta
%=(\frac{2\pi}{\omega})^2\int_{\mathbb{S}} (\zeta\cdot \zeta)^{4-j-{j'}} {\rm d}\zeta=\frac{(2\pi)^3}{\omega^2},
%\end{align*}
Furthermore, when $\|\xi\|_2>2\max(\omega_1,\omega_2)=2\omega_2$, for any $i=1,2$, there hold 
\begin{align*}
[\mathcal{F}(I_{\omega_i}^{j,j'})](\xi)
=&(2\pi)^2\int_{\mathbb{S}\times\mathbb{S}} (\hat{x}\cdot \zeta)^{4-j-{j'}} 0{\rm d}\hat{x}{\rm d}\zeta
=0, \quad [\mathcal{F}(g_{de,\alpha})](\xi)=0 \quad\mbox{and}\quad [\mathcal{F}(g_{de})](\xi)=\alpha^2.
\end{align*}
Next, by the definition $g_2=\mathcal{F}^{-1}\big(\det(\mathcal{F}({\bf g}))^{-1}\big)=\mathcal{F}^{-1}(\mathcal{F}(g_{de})^{-1})$ and the decompositions $$g_{de}=\alpha^2\delta+g_{de,\alpha}\quad \mbox{and}\quad g_2=\alpha^{-2}\delta-g_{2,\alpha},$$  we have
\begin{align*}
[\mathcal{F}(g_{2,\alpha})](-\xi)=[\mathcal{F}(g_{2,\alpha})](\xi)=\alpha^{-2}-[\mathcal{F}(g_2)](\xi)=\alpha^{-2}-\big([\mathcal{F}(g_{de})](\xi)\big)^{-1}\in[0, \alpha^{-2}],
\end{align*}
and $\mathcal{F}(g_{2,\alpha})=\alpha^{-2}-\alpha^{-2}=0$ when $\|\xi\|_2>2\omega_2$.

Now, we shall estimate the modulus of continuity $w_{g_{2,\alpha}}$ and the maximum magnitude $\|g_{2,\alpha}\|_{L^\infty}$ of $g_{2,\alpha}$.
By computing the inverse Fourier transform of $\mathcal{F}(g_{2,\alpha})$, we derive that 
\begin{align*}
g_{2,\alpha}(y)=(2\pi)^{-2}\int_{\mathbb{R}^2}e^{i\xi\cdot y}[\mathcal{F}(g_{2,\alpha})](\xi){\rm d}\xi
%=(2\pi)^{-2}\int_{B_{2\omega_2}}e^{i\xi\cdot y}[\mathcal{F}(g_{2,\alpha})](\xi){\rm d}\xi
=(2\pi)^{-2}\int_{B_{2\omega_2}}\cos(\xi\cdot y)[\mathcal{F}(g_{2,\alpha})](\xi){\rm d}\xi\in \mathbb{R},
\end{align*}
which implies the estimates 
\begin{align*}
\|g_{2,\alpha}\|_{L^\infty}= &\sup_y|g_{2,\alpha}(y)|\leq (2\pi)^{-2}\int_{B_{2\omega_2}}\alpha^{-2}{\rm d}\xi\leq \frac{\pi(2\omega_2)^2}{(2\pi\alpha)^{2}}=\frac{\omega_2^2}{\pi\alpha^{2}},\\
w_{g_{2,\alpha}}(s)=&\sup_{|y-y'|_\infty\leq s}|g_{2,\alpha}(y)-g_{2,\alpha}(y')|
\leq (2\pi)^{-2}\sup_{|y-y'|_\infty\leq s}\int_{B_{2\omega_2}}|e^{i\xi\cdot y}-e^{i\xi\cdot y'}||[\mathcal{F}(g_{2,\alpha})](\xi)|{\rm d}\xi\\
\leq& (2\pi\alpha)^{-2}\sup_{|y-y'|_\infty\leq s}\int_{B_{2\omega_2}}|e^{i\xi\cdot (y-y')}-1|{\rm d}\xi
\leq(2\pi\alpha)^{-2}\sup_{|y-y'|_\infty\leq s}\int_{B_{2\omega_2}}|i\xi\cdot (y-y')|{\rm d}\xi\\
\leq& (2\pi\alpha)^{-2}\sup_{|y-y'|_\infty\leq s}\int_{B_{2\omega_2}}\|\xi\|_2\|y-y'\|_2{\rm d}\xi
\leq \frac{\sqrt{2}s}{(2\pi\alpha)^{2}}\int_{B_{2\omega_2}}\|\xi\|_2{\rm d}\xi\\
=&\frac{\sqrt{2}s}{(2\pi\alpha)^{2}}\int_{0}^{2\pi}\int_{0}^{2\omega_2}r^2{\rm d}r{\rm d}\theta
%=\frac{\sqrt{2}s}{2\pi\alpha^{2}}\frac{2^3\omega_2^3}{3}
=\frac{4\sqrt{2}\omega_2^3}{3\pi\alpha^{2}}s.
\end{align*}
This completes the proof.
\end{proof}

Next, we establish the Lipschitz estimates of the hypothesis classes $\mathcal{H}^0_{BFNN}$ and $\mathcal{H}^1_{BFNN}$ with respect to the parameters which are used for deriving the generalization errors of the proposed neural networks in Section \ref{sec:gene}.
\begin{lemma}\label{lem:cover_lip_0}
Let $\phi^0,\tilde{\phi^0}\in \mathcal{H}^0_{BFNN}$ be the $0$-level BFNN associated with parameters $K,D$ and $\tilde{K},\tilde{D}$ respectively. Then
there holds
\begin{align*}
\|\phi^0-\tilde{\phi^0}\|:=(\sum_{s=1}^{N}\|\phi^0(\Lambda_s)-\tilde{\phi^0}(\Lambda_s)\|_F^2)^{\frac12}\leq (c_{K}\|K-\tilde{K}\|_2+c_{D}\|D-\tilde{D}\|_2)n_c^{-\frac32}\sqrt{N}
\end{align*}
with the constants $c_{K}=2^3\pi^2 (c_{\omega_1}+c_{\omega_2})(B_D^4+1)^\frac12 B_K B_{in}$ and $c_{D}=2^3\pi^2 (c_{\omega_1}+c_{\omega_2})B_DB_K^2 B_{in}$.
\end{lemma}
\begin{proof}
By the definition of the hypothesis class $\mathcal{H}^0_{BFNN}$ in Section \ref{sec:hypothesis}, 
we can estimate the distance between any two hypothesis operators $\phi^0,\tilde{\phi^0}\in\mathcal{H}^0_{BFNN}$ applied to any input $\Lambda_s\in\mathcal{S}$ with the shifted inputs $\Lambda_{\theta_m}^\omega:=R_m^{\omega}(\Lambda_s)$ by
\begin{align*}
\|\phi^0(\Lambda_s)-\tilde{\phi^0}(\Lambda_s)\|_F^2=&\sum_{m=1}^{n_\theta}\big(\|\sum_{i=1}^2(\phi_1^0 \Lambda_{\theta_m}^{\omega_i}-\tilde{\phi_1^0} \Lambda_{\theta_m}^{\omega_i}) P_{\omega_i}\|_2^2+\|\sum_{i=1}^2(\phi_2^0 \Lambda_{\theta_m}^{\omega_i}-\tilde{\phi_2^0} \Lambda_{\theta_m}^{\omega_i}) P_{\omega_i}\|_2^2\big)\\
:=&\sum_{m=1}^{n_\theta}({\rm I}_{1,m}^2+{\rm I}_{2,m}^2),
\end{align*}
where 
\begin{align*}
{\rm I}_{1,m}\leq& \sum_{i=1}^2\|\phi_1^0 \Lambda_{\theta_m}^{\omega_i}-\tilde{\phi_1^0} \Lambda_{\theta_m}^{\omega_i}\|_2
\\
=&\sum_{i=1}^2\|-\frac{2\pi^2 c_{\omega_i}}{n_\theta^2}{\rm diag}\big(K^* ( D^*\Lambda_{\theta_m}^{\omega_i}D+D\Lambda_{\theta_m}^{\omega_i}D^*) K-\tilde{K}^* ( \tilde{D}^*\Lambda_{\theta_m}^{\omega_i}\tilde{D}+\tilde{D}\Lambda_{\theta_m}^{\omega_i}\tilde{D}^*) \tilde{K}\big)\|_2\\
\leq &\sum_{i=1}^2\frac{2\pi^2 c_{\omega_i}}{n_\theta^2}\|K^* ( D^*\Lambda_{\theta_m}^{\omega_i}D+D\Lambda_{\theta_m}^{\omega_i}D^*) K-\tilde{K}^* ( \tilde{D}^*\Lambda_{\theta_m}^{\omega_i}\tilde{D}+\tilde{D}\Lambda_{\theta_m}^{\omega_i}\tilde{D}^*) \tilde{K}\|_F:=\sum_{i=1}^2\frac{2\pi^2 c_{\omega_i}}{n_\theta^2}{\rm I}_{1,m,i},\\
{\rm I}_{2,m}\leq &\sum_{i=1}^2\|\phi_2^0 \Lambda_{\theta_m}^{\omega_i}-\tilde{\phi_2^0} \Lambda_{\theta_m}^{\omega_i}\|_2
=\sum_{i=1}^2\|\frac{4\pi^2 c_{\omega_i}}{n_\theta^2}{\rm diag}\big(K^* \Lambda_{\theta_m}^{\omega_i} K-\tilde{K}^* \Lambda_{\theta_m}^{\omega_i} \tilde{K}\big)\|_2\\
\leq &\sum_{i=1}^2\frac{4\pi^2 c_{\omega_i}}{n_\theta^2}\|K^* \Lambda_{\theta_m}^{\omega_i} K-\tilde{K}^* \Lambda_{\theta_m}^{\omega_i} \tilde{K}\|_F:=\sum_{i=1}^2\frac{4\pi^2 c_{\omega_i}}{n_\theta^2}{\rm I}_{2,m,i}.
\end{align*}
We discuss the terms ${\rm I}_{1,m,i}$ and ${\rm I}_{2,m,i}$ one by one. First, we decompose the above upper bound of ${\rm I}_{1,m,i}$ into several parts
\begin{align*}
{\rm I}_{1,m,i}
\leq &\|(DK)^*\Lambda_{\theta_m}^{\omega_i}DK-(\tilde{D}\tilde{K})^*\Lambda_{\theta_m}^{\omega_i}\tilde{D}\tilde{K}\|_F+\|(D^*K)^*\Lambda_{\theta_m}^{\omega_i}D^* K-(\tilde{D}^*\tilde{K})^* \Lambda_{\theta_m}^{\omega_i}\tilde{D}^*\tilde{K}
\|_F\\
\leq &\|\big((DK)^*-(\tilde{D}\tilde{K})^*\big)\Lambda_{\theta_m}^{\omega_i}DK\|_F+\|(\tilde{D}\tilde{K})^*\Lambda_{\theta_m}^{\omega_i}\big(DK-\tilde{D}\tilde{K}\big)\|_F\\
&+\|\big((D^*K)^*-(\tilde{D}^*\tilde{K})^*\big)\Lambda_{\theta_m}^{\omega_i}D^*K\|_F+\|(\tilde{D}^*\tilde{K})^*\Lambda_{\theta_m}^{\omega_i}\big(D^*K-\tilde{D}^*\tilde{K}\big)\|_F\\
\leq &\|\Lambda_{\theta_m}^{\omega_i}\|_F\Big(\|(DK)^*-(\tilde{D}\tilde{K})^*\|_2\|DK\|_2+\|(\tilde{D}\tilde{K})^*\|_2\|DK-\tilde{D}\tilde{K}\|_2\\
&\qquad\quad\;\;+\|(D^*K)^*-(\tilde{D}^*\tilde{K})^*\|_2\|D^*K\|_2+\|(\tilde{D}^*\tilde{K})^*\|_2\|D^*K-\tilde{D}^*\tilde{K}\|_2\Big)\\
\leq &2B_DB_K B_{in}\Big(\|DK-\tilde{D}\tilde{K}\|_2+\|D^*K-\tilde{D}^*\tilde{K}\|_2\Big),
\end{align*}
where $\|DK-\tilde{D}\tilde{K}\|_2$ and $\|D^*K-\tilde{D}^*\tilde{K}\|_2$ can be further decomposed into
\begin{align*}
\|DK-\tilde{D}\tilde{K}\|_2= &\|D(K-\tilde{K})+(D-\tilde{D})\tilde{K}\|_2\leq B_D\|K-\tilde{K}\|_2+B_K\|D-\tilde{D}\|_2,\\
\|D^*K-\tilde{D}^*\tilde{K}\|_2=&\|D^*(K-\tilde{K})+(D^*-\tilde{D}^*)\tilde{K}\|_2\leq B_D\|K-\tilde{K}\|_2+B_K\|D-\tilde{D}\|_2,
\end{align*}
which imply immediately that
\begin{align*}
{\rm I}_{1,m,i}
\leq &4B_DB_K B_{in}\Big(B_D\|K-\tilde{K}\|_2+B_K\|D-\tilde{D}\|_2\Big):={\rm I}_1.
\end{align*}
Similarly, we can bound the second term ${\rm I}_{2,m,i}$ by
\begin{align*}
{\rm I}_{2,m,i}\leq \|(K-\tilde{K})^* \Lambda_{\theta_m}^{\omega_i} K\|_F+\|\tilde{K}^*\Lambda_{\theta_m}^{\omega_i} (K-\tilde{K})\|_F\leq 2B_K B_{in}\|K-\tilde{K}\|_2:={\rm I}_2.
\end{align*}
Finally, combining all the above estimates and the assumption $n_\theta=n_c$ yields that 
\begin{align*}
&\|\phi^0(\Lambda_s)-\tilde{\phi^0}(\Lambda_s)\|_F=\big(\sum_{m=1}^{n_\theta}({\rm I}_{1,m}^2+{\rm I}_{2,m}^2)\big)^{\frac12}\leq \Big(\sum_{m=1}^{n_\theta}\big((\sum_{i=1}^2\frac{2\pi^2 c_{\omega_i}}{n_\theta^2}{\rm I}_{1,m,i})^2+(\sum_{i=1}^2\frac{4\pi^2 c_{\omega_i}}{n_\theta^2}{\rm I}_{2,m,i})^2\big)\Big)^\frac12\\
\leq&\frac{2\pi^2 (c_{\omega_1}+c_{\omega_2})}{n_\theta^2}\Big(\sum_{m=1}^{n_\theta}\big(({\rm I}_1)^2+(2{\rm I}_2)^2\big)\Big)^\frac12
\leq 2\pi^2 (c_{\omega_1}+c_{\omega_2})n_\theta^{-\frac32}\big(({\rm I}_1)^2+(2{\rm I}_2)^2\big)^\frac12\\
%\leq& 2\pi^2 (c_{\omega_1}+c_{\omega_2})n_\theta^{-\frac32}\Big(\big(4B_DB_K B_{in}(B_D\|K-\tilde{K}\|_2+B_K\|D-\tilde{D}\|_2)\big)^2+(4B_K B_{in}\|K-\tilde{K}\|_2)^2\big)^\frac12\\
\leq& 2^3\pi^2 (c_{\omega_1}+c_{\omega_2})B_K B_{in}n_\theta^{-\frac32}\Big((B_D^2\|K-\tilde{K}\|_2+B_DB_K\|D-\tilde{D}\|_2)^2+\|K-\tilde{K}\|_2^2\big)^\frac12\\
\leq& 2^3\pi^2 (c_{\omega_1}+c_{\omega_2})B_K B_{in}n_c^{-\frac32}\big((B_D^4+1)^\frac12\|K-\tilde{K}\|_2+B_DB_K\|D-\tilde{D}\|_2\big)
\end{align*}
and completes the proof.
\end{proof}

\begin{lemma}\label{lem:cover_lip_1}
Let $\phi^1,\tilde{\phi^1}\in \mathcal{H}^1_{BFNN}$ be the $1$-level BFNN associated with parameters $U,M,V,D$ and $\tilde{U},\tilde{M},\tilde{V},\tilde{D}$ respectively. Then
there holds
\begin{align*}
\|\phi^1-\tilde{\phi^1}\|:=&(\sum_{s=1}^{N}\|\phi^1(\Lambda_s)-\tilde{\phi^1}(\Lambda_s)\|_F^2)^{\frac12}\\
\leq& (c_{U}  \|U-\tilde{U}\|_2+c_{M}\|M-\tilde{M}\|_2+c_{V} \|V-\tilde{V}\|_2+c_{D_r}\|D-\tilde{D}\|_2)n_c^{-\frac32}\sqrt{N}
\end{align*}
where the constants $c_{U}=c_{K_r}B_M B_V$,  $c_{M}=c_{K_r}B_U B_V$, $c_{V}=c_{K_r}B_U B_M$ and $c_{D_r}=2^3\pi^2 (c_{\omega_1}+c_{\omega_2})B_D B_U^2 B_M^2 B_V^2 B_{in}$ with $c_{K_r}=2^3\pi^2 (c_{\omega_1}+c_{\omega_2})(B_D^4+1)^\frac12 B_U B_M B_V B_{in}$.
\end{lemma}
\begin{proof}
Let $K_r=UMV$, then there holds $\|K_r\|_2\leq B_U B_M B_V$.
By the definition of the hypothesis class $\mathcal{H}^1_{BFNN}$ in Section \ref{sec:hypothesis}, it follows from Lemma \ref{lem:cover_lip_0} that
\begin{align*}
\|\phi^1-\tilde{\phi^1}\|\leq (c_{K_r}\|K_r-\tilde{K_r}\|_2+c_{D_r}\|D-\tilde{D}\|_2)n_c^{-\frac32}\sqrt{N}
\end{align*}
which can be further decomposed as
\begin{align*}
\|K_r-\tilde{K}_r\|_2=& \|UMV-\tilde{U}\tilde{M}\tilde{V}\|_2\leq \|(U-\tilde{U})MV\|_2+\|\tilde{U}(M-\tilde{M})V\|_2+\|\tilde{U}\tilde{M}(V-\tilde{V})\|_2\\
\leq& \|U-\tilde{U}\|_2\|M\|_2\|V\|_2+\|\tilde{U}\|_2\|M-\tilde{M}\|_2\|V\|_2+\|\tilde{U}\|_2\|\tilde{M}\|_2\|V-\tilde{V}\|_2\\
\leq& B_M B_V  \|U-\tilde{U}\|_2+B_U B_V\|M-\tilde{M}\|_2+B_U B_M \|V-\tilde{V}\|_2.
\end{align*}
It completes the proof.
\end{proof}

\bibliographystyle{abbrv}
\bibliography{NN}

\begin{thebibliography}{10}

\bibitem{uni2}
K.~Astala and L.~P\"aiv\"arinta.
\newblock Calder\'on's inverse conductivity problem in the plane.
\newblock {\em Ann. of Math. (2)}, 163(1):265--299, 2006.

\bibitem{john2}
K.~Astala, L.~P\"aiv\"arinta, and M.~Lassas.
\newblock Calder\'on's inverse problem for anisotropic conductivity in the plane.
\newblock {\em Comm. Partial Differential Equations}, 30(1-3):207--224, 2005.

\bibitem{uni1}
A.~L. Bukhgeim.
\newblock Recovering a potential from {C}auchy data in the two-dimensional case.
\newblock {\em J. Inverse Ill-Posed Probl.}, 16(1):19--33, 2008.

\bibitem{CakoniColton:2005}
F.~Cakoni and D.~Colton.
\newblock {\em Qualitative Methods in Inverse Scattering Theory}.
\newblock The MIT Press, Cambridge, Massachusetts, London, England, 2005.

\bibitem{CandesDemanetYing:2009}
E.~Cand\`{e}s, L.~Demanet, and L.~Ying.
\newblock A fast butterfly algorithm for the computation of {F}ourier integral operators.
\newblock {\em Multiscale Modeling \& Simulation}, 7(4):1727--1750, 2009.

\bibitem{colton-kress}
D.~Colton and R.~Kress.
\newblock {\em Inverse Acoustic and Electromagnetic Scattering Theory}, volume~93 of {\em Applied Mathematical Sciences}.
\newblock Springer, Cham, fourth edition, 2019.

\bibitem{Costabel:2015}
M.~Costabel.
\newblock On the spectrum of volume integral operators in acoustic scattering.
\newblock In C.~Constanda and A.~Kirsch, editors, {\em Integral Methods in Science and Engineering}, pages 119--127, 2015.

\bibitem{DeVoreHaninPetrova:2021}
R.~DeVore, B.~Hanin, and G.~Petrova.
\newblock Neural network approximation.
\newblock {\em Acta Numer.}, 30:327--444, 2021.

\bibitem{EnglHankeNeubauer:1996}
H.~W. Engl, M.~Hanke, and A.~Neubauer.
\newblock {\em Regularization of {I}nverse {P}roblems}.
\newblock Kluwer, Dordrecht, 1996.

\bibitem{FanFeliuLin:2019}
Y.~Fan, J.~Feliu-Fabà, L.~Lin, L.~Ying, and L.~Zepeda-Núñez.
\newblock A multiscale neural network based on hierarchical nested bases.
\newblock {\em Res Math Sci}, 6(21), 2019.

\bibitem{FoucartRauhut2013}
S.~Foucart and H.~Rauhut.
\newblock {\em A Mathematical Introduction to Compressive Sensing.}
\newblock Applied and Numerical Harmonic Analysis. Birkhäuser, 2013.

\bibitem{uniq2}
F.~Gylys-Colwell.
\newblock An inverse problem for the {H}elmholtz equation.
\newblock {\em Inverse Problems}, 12(2):139--156, 1996.

\bibitem{HeZhangRen:2016}
K.~He, X.~Zhang, S.~Ren, and J.~Sun.
\newblock Deep residual learning for image recognition.
\newblock In {\em Proceedings of the IEEE Conference on Computer Vision and Pattern Recognition (CVPR)}, 2016.

\bibitem{KingmaBa:2015}
D.~P. Kingma and J.~Ba.
\newblock Adam: A method for stochastic optimization.
\newblock In Y.~Bengio and Y.~LeCun, editors, {\em 3rd International Conference for Learning Representations, San Diego}, 2015.

\bibitem{KV}
R.~V. Kohn and M.~Vogelius.
\newblock Relaxation of a variational method for impedance computed tomography.
\newblock {\em Comm. Pure Appl. Math.}, 40(6):745--777, 1987.

\bibitem{LiDemanetZepeda:2022}
M.~Li, L.~Demanet, and L.~Zepeda-N\'{u}\~{n}ez.
\newblock Wide-band butterfly network: Stable and efficient inversion via multi-frequency neural networks.
\newblock {\em Multiscale Modeling \& Simulation}, 20(4):1191--1227, 2022.

\bibitem{LiYangMartin:2015}
Y.~Li, H.~Yang, E.~R. Martin, K.~L. Ho, and L.~Ying.
\newblock Butterfly factorization.
\newblock {\em Multiscale Modeling \& Simulation}, 13(2):714--732, 2015.

\bibitem{Maurer2016}
A.~Maurer.
\newblock A vector-contraction inequality for {R}ademacher complexities.
\newblock In R.~Ortner, H.~U. Simon, and S.~Zilles, editors, {\em Algorithmic Learning Theory}, pages 3--17. Springer International Publishing, 2016.

\bibitem{SchnoorBehboodiRauhut2023}
E.~Schnoor, A.~Behboodi, and H.~Rauhut.
\newblock {Generalization error bounds for iterative recovery algorithms unfolded as neural networks}.
\newblock {\em Information and Inference: A Journal of the IMA}, 12(3):2267--2299, 2023.

\bibitem{Shalev-ShwartzBen-David2014}
S.~Shalev-Shwartz and S.~Ben-David.
\newblock {\em Understanding Machine Learning: From Theory to Algorithms}.
\newblock Cambridge University Press, USA, 2014.

\bibitem{john1}
J.~Sylvester.
\newblock An anisotropic inverse boundary value problem.
\newblock {\em Comm. Pure Appl. Math.}, 43(2):201--232, 1990.

\bibitem{ZhangZepeda-NunezLi:2024}
B.~Zhang, L.~Zepeda-Nunez, and Q.~Li.
\newblock Solving the wide-band inverse scattering problem via equivariant neural networks.
\newblock {\em Journal of Computational and Applied Mathematics}, 451:116050, 2024.

\end{thebibliography}
\end{document}